\documentclass[preprint,11pt,oneside,leqno]{imsart}
\RequirePackage[T1]{fontenc}
\RequirePackage{times}
\RequirePackage[a4paper,textwidth=6.25in,textheight=9in,centering]{geometry}
\RequirePackage{amsmath,amssymb,amsthm,mathtools}
\RequirePackage[numbers,sort&compress]{natbib}
\RequirePackage{hyperref}
\hypersetup{hidelinks}
\RequirePackage{graphicx,booktabs}

\startlocaldefs
\numberwithin{equation}{section}
\theoremstyle{plain}
\newtheorem{theorem}{Theorem}[section]
\newtheorem{lemma}[theorem]{Lemma}
\newtheorem{proposition}[theorem]{Proposition}
\newtheorem{corollary}[theorem]{Corollary}
\newtheorem{conjecture}[theorem]{Conjecture}
\theoremstyle{definition}
\newtheorem{definition}[theorem]{Definition}
\newtheorem{assumption}[theorem]{Assumption}
\theoremstyle{definition}
\newtheorem{remark}[theorem]{Remark}
\newtheorem{example}[theorem]{Example}
\newcommand{\R}{\mathbb R}
\newcommand{\E}{\mathbb E}

\newcommand{\cO}{\mathcal O}
\newcommand{\cF}{\mathcal F}
\newcommand{\cI}{\mathcal I}
\newcommand{\cL}{\mathcal L}
\newcommand{\norm}[1]{\lVert #1\rVert}
\newcommand{\ip}[2]{\langle #1,#2\rangle}
\newcommand{\supp}{\operatorname{supp}}
\newcommand{\diver}{\operatorname{div}}
\newcommand{\tr}{\operatorname{tr}}
\newcommand{\KL}{\operatorname{KL}}
\newcommand{\Teff}{T}            
\newcommand{\lmax}{\lambda_{\max}}
\newcommand{\lmin}{\lambda_{\min}}
\newcommand{\rate}{a_N}          
\newcommand{\Risk}{\mathfrak R_N}
\newcommand{\condC}{\textup{(C)}}

\endlocaldefs

\begin{document}
\begin{frontmatter}
\title{Minimax estimation of linear functionals in the Darcy inverse problem}
\runtitle{Linear functionals in the Darcy problem}

\begin{aug}
\author[A]{\fnms{Yi}~\snm{Liu}}
\author[B]{\fnms{Xiaotian}~\snm{Chang}}
\author[C]{\fnms{Qirui}~\snm{Hu}}
\address[A]{York University}
\address[B]{Nanyang Technological University}
\address[C]{Shanghai University of Finance and Economics}
\runauthor{Y. Liu, X. Chang and Q. Hu}
\end{aug}

\begin{abstract}
We study estimation of linear functionals $\int f\psi$ of the conductivity $f$ in the Darcy problem $\nabla\cdot(f\nabla u)=g$ from observations of $u$ in Gaussian white noise of level $N^{-1/2}$. When the potential has a single nondegenerate minimum, we determine the minimax rates over Sobolev balls in the linearised experiment. Unless a compatibility condition along the gradient flow holds, the rate is governed by an effective dimension $T$ of the Hessian at the minimum, not by the ambient dimension: it is $N^{-(2s+T)/(4s+8)}$ for $T<4$, has an elbow at $T=4$ with an extra logarithmic factor, and is parametric for $T>4$. In the nonlinear experiment we show that these lower bounds persist when compatibility fails at the reference and $T\le4$, and we prove the minimax rate on classes with isotropic critical Hessians, and the parametric rate near references with $T>4$. In dimension at most four, no nonzero linear functional is uniformly $\sqrt N$-estimable on a Sobolev neighbourhood of the constant conductivity, even when it is regular there.
\end{abstract}


\end{frontmatter}
\section{Introduction}\label{sec:intro}

The Darcy problem of recovering the conductivity, or permeability, $f$ of a medium from the potential $u_f$ it induces is a basic model for subsurface flow. It has become the standard test case for the statistical theory of nonlinear PDE-constrained inverse problems \cite{richter1981inverse,stuart2010inverse,nickl2020convergence,giordano2020consistency,nickl2023bayesian}. In applications one is often less interested in the whole field $f$ than in weighted averages $\int f\psi$, and such scalar summaries are also the natural targets of uncertainty quantification. This paper determines how well they can be estimated.

Let $\cO\subset\R^d$, $d\ge2$, be a bounded domain with smooth boundary, let $g>0$ be a smooth source, and for a positive conductivity $f$ let $u_f$ solve
\begin{equation}\label{eq:darcy}
 \diver(f\nabla u_f)=g\quad\text{in }\cO,\qquad u_f=0\quad\text{on }\partial\cO.
\end{equation}
We observe
\begin{equation}\label{eq:model}
 Y=u_f+N^{-1/2}\,\mathbb W,
\end{equation}
where $\mathbb W$ is Gaussian white noise on $L^2(\cO)$, and we wish to estimate the linear functional
\begin{equation}\label{eq:functional}
 \Phi(f)=\int_\cO f\psi,\qquad \psi\in C_c^\infty(\cO).
\end{equation}
Estimation of $f$ itself is well understood. Rates of convergence for penalised least squares and Bayesian methods were obtained in \cite{nickl2020convergence,giordano2020consistency}; see also \cite{nickl2023bayesian,agapiou2024laplace,siebel2025convergence,koers2024linear,wang2026global}. Rei\ss\ and Wang \cite{reiss2026darcy} recently identified the minimax polynomial $L^2$ rate, up to logarithmic factors, on local H\"older--Zygmund classes around the constant conductivity on the ball. In their analysis the rate is driven by the degeneracy of the problem at the critical point of the potential.

For linear functionals one might hope for the parametric rate $N^{-1/2}$, as in several other nonlinear inverse problems, where the posterior distribution of such functionals even satisfies a Bernstein--von Mises theorem \cite{monard2019efficient,nickl2020bvm,monard2021statistical}. Nickl and Paternain \cite[Section~3.5]{nickl2023information} showed that in the Darcy problem this can fail: the efficient information for $\int f\psi$ vanishes unless a compatibility condition along the characteristics of a transport equation holds. Their analysis assumes $\psi(x_*)=0$ at the critical point $x_*$ of the potential and $d\le3$, and treats the unit disc in detail. The present paper quantifies this phenomenon. We identify the minimax rate when the information vanishes, in every dimension and without assuming $\psi(x_*)=0$, and we show how the picture changes in the nonlinear experiment.

\subsection{Main results}\label{sec:intro-results}

The linearisation of $f\mapsto u_f$ at $f$ is $h\mapsto A_fh=-L_f^{-1}\diver(h\nabla u_f)$, where $L_f=\diver(f\nabla\cdot)$ with Dirichlet boundary conditions. The functional $h\mapsto\ip h\psi$ has an $L^2$ representer, $\ip h\psi=\ip{A_fh}{z}$, exactly when $z=L_fw$ for a solution $w\in H^2\cap H^1_0$ of the transport equation
\[
 \nabla u_f\cdot\nabla w=\psi,\qquad w|_{\partial\cO}=0 .
\]
Under our assumptions $u_f$ is negative with a single nondegenerate minimum $x_*$. All gradient-flow lines start at $x_*$ and end on $\partial\cO$, and $w$ is the integral of $\psi$ along them. Generically $w$ is singular at $x_*$. It contains $\psi(x_*)$ times the flow time, which grows logarithmically, and a term that is constant on each flow line but varies across flow lines. How singular $w$ is depends on the anisotropy of the Hessian $H=D^2u_f(x_*)$ through the \emph{effective dimension}
\begin{equation}\label{eq:Teff}
 \Teff=\Teff(H)=\frac{\tr H}{\lmax(H)}\in(1,d],\qquad \Teff=d\iff H\text{ is a multiple of the identity}.
\end{equation}
Theorem~\ref{thm:linear} shows that in the linearised experiment the minimax rate is
\begin{equation}\label{eq:rate}
 \rate(s,\Teff)=\begin{cases}
 N^{-(2s+\Teff)/(4s+8)}, & \Teff<4,\\
 \sqrt{\log N/N}, & \Teff=4,\\
 N^{-1/2}, & \Teff>4,
 \end{cases}
\end{equation}
over Sobolev balls of smoothness $s$, unless an explicit, coordinate-free compatibility condition \condC\ holds (Definition~\ref{def:C}), in which case the rate is $N^{-1/2}$. We call these the \emph{irregular regime} ($\Teff<4$ and \condC\ fails), the \emph{critical case} ($\Teff=4$ and \condC\ fails) and the \emph{regular regime} (\condC\ holds or $\Teff>4$).

The exponent can be read off from a bias--variance balance. Estimate $\Phi$ by pairing the data with a representer truncated at frequency $R$. The part of $\psi$ that the truncated representer misses has $H^{-s}$-norm of order $R^{-s-\Teff/2}$; this is the bias. The variance is $N^{-1}\norm{z_R}_2^2$. Here $z=L_{f}w$ involves two derivatives of $w$, which are of order $r^{-2}$ at distance $r$ from $x_*$ in the fastest flow direction. The gradient flow maps flow-adapted shells onto each other, and the shell at scale $r$, of width $r$ in the fastest direction and $r^{\lambda_i/\lmax}$ in the direction of the $i$th eigenvector, has volume $r^{\Teff}$. Hence
\[
 \norm{z_R}_2^2\asymp\int_{1/R}^{1}r^{-4}\,r^{\Teff-1}\,dr
 \asymp\begin{cases}R^{4-\Teff},&\Teff<4,\\ \log R,&\Teff=4,\\ 1,&\Teff>4.\end{cases}
\]
For $\Teff<4$, balancing $R^{-s-\Teff/2}$ against $N^{-1/2}R^{(4-\Teff)/2}$ at $R=N^{1/(2s+4)}$ gives the first case of \eqref{eq:rate}. For $\Teff>4$ the exact representer is square integrable and the rate is parametric, while at $\Teff=4$ the logarithmic divergence produces the factor $\sqrt{\log N}$. The lower bounds reverse this calculation. A perturbation of the potential supported in a flow-adapted box of scale $r=N^{-1/(2s+4)}$ cannot be detected at noise level $N^{-1/2}$, but moves $\Phi$ by the amount in \eqref{eq:rate} (Section~\ref{sec:ideas}).

The rate thus has an elbow at $\Teff=4$, and the location of the elbow does not depend on $s$; Figure~\ref{num:fig:exponents} shows the computed exponents as a function of $\Teff$. Since $\Teff\le d$, every functional violating \condC\ is irregular or critical when $d\le4$. Thus, for functionals violating \condC, $d\ge5$ is necessary but not sufficient for parametric estimation; for instance, $H=\operatorname{diag}(1,1,1,1,2)$ in $d=5$ has $\Teff=3$. In the irregular regime the exponent approaches $1/2$ as $s\to\infty$, with gap $(4-\Teff)/(4s+8)$.

These statements concern a fixed linearisation. In the nonlinear experiment over a Sobolev neighbourhood of a reference conductivity $f_0$, both the Hessian at the critical point and the flow-line integrals vary with $f$, and the picture changes in three ways.
\begin{enumerate}
\item[(i)] The linearised lower bounds persist (Theorem~\ref{thm:lower}). If \condC\ fails at $f_0$ and $\Teff_0\le4$, the minimax risk over the neighbourhood is at least of order $\rate(s,\Teff_0)$. If $\psi(x_*)\ne0$, it is at least of order $N^{-(2s+d)/(4s+8)}$ for every critical Hessian, even when $f$ is known away from a neighbourhood of $x_*$.
\item[(ii)] On classes of conductivities near the constant conductivity on the unit ball whose critical Hessians are all isotropic, the minimax rate is $\rate(s,d)$ for \emph{every} nonzero $\psi$ (Theorem~\ref{thm:isotropic}).
\item[(iii)] If $\Teff_0>4$, the rate on a small full Sobolev neighbourhood is parametric (Theorem~\ref{thm:full-regular}). Near the constant conductivity in $d\le4$, by contrast, no nonzero functional is uniformly $\sqrt N$-estimable, and the local rate exponent converges to that of $\rate(s,d)$ as the neighbourhood shrinks (Proposition~\ref{prop:instability} and Corollary~\ref{cor:local-exponent}).
\end{enumerate}
The distinction between compatibility at a reference and compatibility throughout a class is essential here. On the canonical isotropic class, the uniform minimax rate is $\rate(s,d)$ even for functionals satisfying \condC\ at the reference. The reason is that compatibility can fail at arbitrarily nearby conductivities with the same critical point and Hessian (Proposition~\ref{prop:instability}). By contrast, imposing \condC\ at every member of a class permits a uniform parametric upper bound (Proposition~\ref{up:prop:uniformC}).

\subsection{Relation to previous work and methods}\label{sec:intro-related}

In linear Gaussian models the minimax risk for a linear functional is governed by the modulus of continuity \cite{ibragimov1984nonparametric,donoho1994statistical,cai2004minimax}, and its rate for indirect observations is determined by the smoothness of the representer relative to the parameter class \cite{goldenshluger2000adaptive,cavalier2008nonparametric,butucea2009adaptive}. Theorem~\ref{thm:linear} computes this modulus for the Darcy forward operator, whose representers for smooth functionals may be singular at the critical point. The nonlinear Darcy problem has been studied extensively. Stability estimates go back to Richter \cite{richter1981inverse} and Alessandrini \cite{alessandrini1986identification}, with related internal-data problems treated in \cite{kuchment2012stabilizing,bal2013reconstruction,montalto2013stability}. Statistical guarantees for $f$ are in \cite{nickl2020convergence,giordano2020consistency,nickl2023bayesian,koers2024linear,wang2026global,reiss2026darcy}, and numerical recovery is studied in \cite{bonito2017diffusion}; see \cite{abraham2020calderon,monard2021statistical,nickl2024lowfreq,giordano2024lowfreq,nickl2024navier,nickl2026killed} for other PDE models.

For functionals of the conductivity, R\o mer \cite[Theorem~3.1 and Remark~4]{romer2024valid} constructs asymptotically valid credible ellipsoids for linear functionals of the log conductivity, and Baek and Papagiannouli \cite[Section~4.1 and Theorem~4.3]{baek2024frequentist} obtain contraction and conservative coverage for functionals of a linked conductivity parameter, under strong smoothness and compatibility assumptions. Their credible sets shrink at nearly parametric rates for very smooth truths. Our target is the minimax risk for $\int f\psi$, uniformly over Sobolev classes of conductivities of finite smoothness. This is consistent with their findings, since the polynomial exponent of $\rate(s,\Teff)$ approaches $1/2$ as $s\to\infty$; the critical logarithmic factor remains when $\Teff=4$.

Transitions from nonparametric to parametric rates are classical for \emph{nonlinear} functionals. For the integrated squared density, Bickel and Ritov \cite{bickel1988estimating} found the elbow at smoothness $d/4$; see also \cite{laurent1996efficient,cai2006optimal}. Analogous elbows occur for functionals in semiparametric models with rough nuisance parameters \cite{robins2008higher,robins2017minimax,kennedy2024minimax}, on sparsity classes \cite{collier2017minimax}, and under local differential privacy \cite{butucea2023interactive}. In indirect observation models the location of the elbow moves with the degree of ill-posedness \cite{butucea2007goodness}. For a \emph{linear} functional observed through a fixed linear operator, the rate is parametric exactly when the representer is square integrable. In the Darcy problem this is decided by the degeneracy of the forward operator at the critical point. The elbow therefore occurs at a value of a geometric quantity, $\Teff=4$, rather than at a smoothness threshold.

Our upper bounds are attained by a one-step estimator whose correction term is a regularised transport representer computed from a pilot. The lower bounds use perturbations prescribed through the forward map, namely localised flow-adapted changes of the potential, and this is what makes the effective dimension appear. Section~\ref{sec:ideas} explains both constructions.

\subsection{Organisation and notation}\label{sec:intro-notation}

Section~\ref{sec:setting} introduces the setting, the parameter classes and condition \condC. Section~\ref{sec:results} states the main results, which are summarised in Table~\ref{tab:classes}. Section~\ref{sec:ideas} explains the main ideas of the proofs, Section~\ref{sec:numerics} gives a numerical illustration, and Section~\ref{sec:discussion} discusses open problems. Complete proofs and numerical details are in the Supplement \cite{liu2026supp}, whose sections, results and equations carry the prefix S.

$\norm{\cdot}_2$ is the $L^2(\cO)$ norm, $H^s=H^s(\cO)$ is the ordinary Sobolev space, and $H^1_0$ is the closure of $C_c^\infty(\cO)$ in $H^1$. $L_f v=\diver(f\nabla v)$ with Dirichlet boundary conditions. For a potential $u$ we write $T_uh=\diver(h\nabla u)$. For a source $F$, $S_uF$ denotes the solution $h$ of $T_uh=F$ defined by \eqref{lin:eq:S} through the backward flow towards the critical point (Section~\ref{lin:sec:flow}). The operator $T_u$ should not be confused with the number $\Teff(H)$. $a\lesssim b$ means $a\le Cb$ for a constant independent of $N$ and of the parameter, and $a\asymp b$ means $a\lesssim b\lesssim a$. Rates refer to the expected absolute error $\E_f|\widehat\Phi-\Phi(f)|$. The minimax risk over a class $\cF$ is
\[
 \Risk(\cF)=\inf_{\widehat\Phi}\sup_{f\in\cF}\E_f|\widehat\Phi-\Phi(f)|,
\]
where the infimum extends over all estimators, that is, all measurable functions $\widehat\Phi=\widehat\Phi(Y)$ of the observation. The stated absolute-error lower bounds also imply squared-error lower bounds at the squared rates, and the upper-bound constructions give the same bounds for root mean squared error. Consequently, for classes with matching minimax rates, the minimax mean squared error has the order of the squared stated rate.
\section{Setting, parameter classes and condition (C)}\label{sec:setting}

\begin{assumption}\label{ass:base}
$\cO\subset\R^d$, $d\ge2$, is a bounded connected domain with $C^\infty$ boundary. The source $g\in C^\infty(\overline\cO)$ and the reference conductivity $f_0\in C^\infty(\overline\cO)$ are bounded below by positive constants. The reference solution $u_0=u_{f_0}$ has exactly one critical point $x_*\in\cO$. It is a nondegenerate minimum with Hessian $H=D^2u_0(x_*)>0$ and eigenvalues $0<\lambda_1\le\dots\le\lambda_d=\lmax(H)$.
\end{assumption}

\subsection{Parameter classes and boundary convention}\label{sec:classes}

For the unit ball $B_1$ and $f_0=g=1$ we have $u_0(x)=(|x|^2-1)/(2d)$ and $H=I_d/d$; we call this the \emph{canonical model}. For $s>d/2$ and $\eta>0$ let
\[
 \cF_s(f_0,\eta)=\bigl\{f=f_0+h:\ \norm h_{H^s(\cO)}\le\eta,\ \inf_\cO f>0\bigr\}.
\]
For $\eta$ sufficiently small the positivity constraint is automatic and $\inf f\ge\inf f_0/2$. If moreover $s>d/2+1$, then $u_f\in C^2(\overline\cO)$ uniformly over the class, and $u_f$ has a unique nondegenerate minimum $x_*(f)$ close to $x_*$. Upper bounds are stated for sufficiently small fixed $\eta$. Lower bounds over a small ball also hold over every larger ball.

The class $\cF_s(f_0,\eta)$ does not fix $f$ near $\partial\cO$. Much of the literature fixes the boundary values of $f$ \cite{nickl2023information} or assumes $f-f_0$ compactly supported \cite{nickl2023bayesian}. The distinction matters for one of the two parts of condition \condC\ below. The obstruction $\psi(x_*)\ne0$ is local: the corresponding lower bounds use perturbations supported in a small neighbourhood of $x_*$, and they, as well as the isotropic rates in that case, remain valid under interior support constraints. The flow-profile obstruction, by contrast, requires perturbations that propagate along whole flow lines to the boundary, and it disappears under interior support constraints (Remark~\ref{rem:collar}). In the anisotropic case with $\psi(x_*)\ne0$, the minimax rate under interior support constraints is open.

\subsection{Linearisation and representer}\label{sec:linearisation}

Subtracting the Darcy equations for $f+h$ and $f$ gives the exact identity
\begin{equation}\label{eq:exact}
 u_{f+h}-u_f=-L_{f+h}^{-1}T_{u_f}h .
\end{equation}
The derivative of $f\mapsto u_f$ is $A_fh=-L_f^{-1}T_{u_f}h$, with adjoint $A_f^*z=\nabla u_f\cdot\nabla(L_f^{-1}z)$. In the \emph{linearised experiment} at $f_0$ one observes $Y=A_{f_0}h+N^{-1/2}\mathbb W$ and estimates $\ip h\psi$ over $\norm h_{H^s}\le\eta$. A representer $z\in L^2$ with $A_{f_0}^*z=\psi$ exists if and only if $z=L_{f_0}w$ for some $w\in H^2(\cO)\cap H^1_0(\cO)$ solving
\begin{equation}\label{eq:transport}
 b\cdot\nabla w=\psi,\qquad w|_{\partial\cO}=0,\qquad b=\nabla u_0 .
\end{equation}
We call $\psi$ (or $\Phi$) \emph{regular at $f_0$} if such a representer exists. Since $u_0<0$ in $\cO$ and $u_0=0$ on $\partial\cO$, every integral curve of $b$ other than the stationary point $x_*$ emanates from $x_*$ in backward time and leaves $\cO$ in finite forward time. The characteristic solution of \eqref{eq:transport} is obtained by integrating $\psi$ along these curves backward from the boundary.

\subsection{The flow profile and condition (C)}\label{sec:conditionC}

Let $\cL$ be the set of oriented integral curves of $b$ in $\cO\setminus\{x_*\}$. For each $\ell\in\cL$ choose a parametrisation $\gamma_\ell$ with $\dot\gamma_\ell(t)=b(\gamma_\ell(t))$, and let $t_\partial(\ell)$ be its exit time. If $\psi(x_*)=0$, then $\psi(\gamma_\ell(t))$ decays exponentially as $t\to-\infty$, and the \emph{flow profile}
\begin{equation}\label{eq:profile}
 \Psi(\ell)=-\int_{-\infty}^{t_\partial(\ell)}\psi(\gamma_\ell(t))\,dt,\qquad \ell\in\cL,
\end{equation}
is well defined. It does not depend on the choice of time origin on each trajectory. Labelling the trajectories by their intersections with a transverse section, for example a small level set of $u_0$, identifies $\cL$ with $S^{d-1}$ and leaves the definition unchanged.

\begin{definition}[Condition (C)]\label{def:C}
$\psi$ satisfies \condC\ at $f_0$ if $\psi(x_*)=0$ and $\Psi$ is constant on $\cL$.
\end{definition}

Condition \condC\ is equivalent to continuity of the characteristic solution $w$ of \eqref{eq:transport} at $x_*$, and then $w$ is smooth (Proposition~\ref{prop:C-smooth}). Under \condC, $\psi$ is regular at $f_0$ in every dimension. If $\psi(x_*)\ne0$, $w$ has a logarithmic singularity at $x_*$ proportional to $\psi(x_*)$. If $\psi(x_*)=0$ but $\Psi$ is not constant, $w$ is bounded but discontinuous at $x_*$, with different limits along different flow lines.
\section{Main results}\label{sec:results}

Throughout, $\psi\in C_c^\infty(\cO)$ is nonzero, $\Teff_0=\Teff(H)$ is the effective dimension \eqref{eq:Teff} at the reference, and $\rate(s,t)$ is defined in \eqref{eq:rate}. Constants may depend on $\cO,g,f_0,\psi,s,\eta$, but not on $N$. Table~\ref{tab:classes} summarises the results.

\begin{table}[t]
\caption{Rates of the minimax risk. Columns two to four assume that \condC\ fails at $f_0$. ${\asymp}$: minimax rate; ${\gtrsim}$, ${\lesssim}$: one-sided bound; open: see Section~\ref{sec:discussion}. Rows: linearised experiment over an $H^s$ ball (Theorem~\ref{thm:linear}); nonlinear experiment over $\cF_s(f_0,\eta)$ (Theorems~\ref{thm:lower} and~\ref{thm:full-regular}); isotropic classes at the canonical model (Theorem~\ref{thm:isotropic}). Nonlinear upper bounds require the smoothness conditions of the cited results.}\label{tab:classes}
\centering\footnotesize
\setlength{\tabcolsep}{5pt}
\begin{tabular}{@{}lcccc@{}}
\toprule
 & $\Teff_0<4$ & $\Teff_0=4$ & $\Teff_0>4$ & \condC\ at $f_0$\\
\midrule
Linearised & $\asymp N^{-\frac{2s+\Teff_0}{4s+8}}$ & $\asymp\sqrt{\log N/N}$ & $\asymp N^{-1/2}$ & $\asymp N^{-1/2}$\\[2pt]
Nonlinear, lower & $\gtrsim N^{-\frac{2s+\Teff_0}{4s+8}}$ & $\gtrsim\sqrt{\log N/N}$ & $\gtrsim N^{-1/2}$ & $\gtrsim N^{-1/2}$\\[2pt]
Nonlinear, upper & open$^{a}$ & open$^{a}$ & $\lesssim N^{-1/2}$ & open$^{b}$\\[2pt]
Isotropic classes & $\asymp\rate(s,d)$, $d\le3$ & $\asymp\sqrt{\log N/N}$, $d=4$ & $\asymp N^{-1/2}$, $d\ge5$ & as when \condC\ fails\\
\bottomrule
\end{tabular}

\smallskip
\begin{minipage}{0.97\textwidth}\footnotesize
$^{a}$At the canonical model the local exponent is identified as $\eta\downarrow0$ (Corollary~\ref{cor:local-exponent}); for anisotropic $f_0$ see Conjecture~\ref{conj:anisotropic}.
$^{b}$If $\Teff_0>4$, the rate is $N^{-1/2}$ by Theorem~\ref{thm:full-regular}. At the canonical model with $d\le4$ the risk is $\gtrsim\rate(s,d)\gg N^{-1/2}$ (Theorem~\ref{thm:isotropic}). If \condC\ holds at every member of a class, a uniform $O(N^{-1/2})$ risk bound is attainable over that class (Proposition~\ref{up:prop:uniformC}).
\end{minipage}
\end{table}

\subsection{The linearised experiment}

\begin{theorem}[Minimax rates in the linearised experiment]\label{thm:linear}
Let Assumption~\ref{ass:base} hold, $s\ge0$ and $\eta>0$. In the linearised experiment at $f_0$,
\[
 \inf_{\widehat\ell}\sup_{\norm h_{H^s}\le\eta}\E_h|\widehat\ell-\ip h\psi|
 \asymp\begin{cases}
 N^{-1/2}, & \text{if \condC\ holds},\\
 N^{-(2s+\Teff_0)/(4s+8)}, & \text{if \condC\ fails and }\Teff_0<4,\\
 \sqrt{\log N/N}, & \text{if \condC\ fails and }\Teff_0=4,\\
 N^{-1/2}, & \text{if \condC\ fails and }\Teff_0>4,
 \end{cases}
\]
where the infimum extends over all measurable functions $\widehat\ell$ of the observation. The minimax mean squared error is of the order of the square of these rates, with constants depending on $\eta$, $s$ and the fixed model. Moreover, $\psi$ is regular at $f_0$ if and only if \condC\ holds or $\Teff_0>4$.
\end{theorem}

In particular, the minimax rate has an elbow at $\Teff_0=4$, and its location does not depend on the smoothness $s$. No nonresonance or linearisability assumption on the gradient flow is needed. Since $\Teff_0\le d$, with equality only for isotropic $H$, the theorem has the following consequence for the efficient information.

\begin{corollary}\label{cor:linear}
If $d\le4$, the linearised functional is $\sqrt N$-estimable, equivalently has positive efficient information \cite{vandervaart1998asymptotic}, if and only if \condC\ holds. If $d\ge5$, this holds if and only if \condC\ holds or $\Teff_0>4$. For every $d\ge5$ there are smooth references at which some functionals are not $\sqrt N$-estimable.
\end{corollary}

Theorem~\ref{thm:linear} and Corollary~\ref{cor:linear} are proved in Section~\ref{sec:linear-proofs} of the Supplement \cite{liu2026supp}.

\subsection{Nonlinear lower bounds}

\begin{theorem}[Nonlinear lower bounds]\label{thm:lower}
Let Assumption~\ref{ass:base} hold, $s>d/2$, and let $\eta>0$ be sufficiently small.
\begin{enumerate}
\item[(a)] If $\psi(x_*)\ne0$, then $\Risk(\cF_s(f_0,\eta))\gtrsim N^{-(2s+d)/(4s+8)}$. This only uses $\nabla u_0(x_*)=0$; nondegeneracy is not needed. The bound also holds if $f$ is known outside any fixed neighbourhood of $x_*$.
\item[(b)] If \condC\ fails at $f_0$ and $\Teff_0\le4$, then $\Risk(\cF_s(f_0,\eta))\gtrsim\rate(s,\Teff_0)$.
\item[(c)] Always, $\Risk(\cF_s(f_0,\eta))\gtrsim N^{-1/2}$.
\end{enumerate}
\end{theorem}

Part (b) shows that the linearised rate of Theorem~\ref{thm:linear} cannot be improved in the nonlinear experiment when $\Teff_0\le4$. Part (a) is weaker for anisotropic $H$, since then $\Teff_0<d$, but its perturbations have interior support; it is informative only for $d\le4$, since for $d\ge5$ it is implied by (c). Theorem~\ref{thm:lower} is proved in Section~\ref{sec:lower-proofs}.

\subsection{Minimax rates on isotropic classes}

Consider the canonical model, and let $s>d/2+1$, so that $u_f\in C^2$ uniformly on $\cF_s(1,\eta)$. Define
\[
 \cI_s(\eta)=\bigl\{f\in\cF_s(1,\eta):\ x_*(f)=0,\ D^2u_f(0)=I_d/d\bigr\},
\]
\[
 \cI_s^{\rm iso}(\eta)=\bigl\{f\in\cF_s(1,\eta):\ D^2u_f(x_*(f))=\lambda_fI_d\text{ for some }\lambda_f>0\bigr\}.
\]
The class $\cI_s(\eta)$ fixes both the critical point and its Hessian, whereas $\cI_s^{\rm iso}(\eta)\supset\cI_s(\eta)$ allows their location and scale to vary while preserving isotropy. Neither class requires the conductivity to be radial or the potential to be quadratic. Since $\tr D^2u_f(x_*)=g(x_*)/f(x_*)$, every $f\in\cI_s(\eta)$ satisfies $f(0)=1$.

\begin{theorem}[Minimax rates on isotropic classes]\label{thm:isotropic}
Consider the canonical model and a nonzero $\psi\in C_c^\infty(B_1)$, let $s>d/2+3$, and let $\eta>0$ be sufficiently small.
\begin{enumerate}
\item[(i)] There is an estimator $\widehat\Phi_N$, the one-step estimator of Section~\ref{up:sec:iso}, and a constant $C$ such that
\[
 \sup_{f\in\cI_s^{\rm iso}(\eta)}\E_f|\widehat\Phi_N-\Phi(f)|\le C\,\rate(s,d)\qquad\text{for all }N\ge2.
\]
\item[(ii)] There is a constant $c>0$ such that $\Risk(\cI_s(\eta))\ge c\,\rate(s,d)$ for all $N\ge2$.
\end{enumerate}
Consequently $\Risk(\cI_s(\eta))\asymp\Risk(\cI_s^{\rm iso}(\eta))\asymp\rate(s,d)$, \emph{whether or not \condC\ holds at $f_0=1$}.
\end{theorem}

On isotropic classes, then, the rate depends only on the dimension and not on any property of $\psi$; even functionals that are regular at $f_0=1$ are estimated at the irregular rate when $d\le4$.

\begin{remark}[Smoothness thresholds]\label{rem:thresholds}
The condition $s>d/2+3$ is a simple sufficient condition. The proof gives the weaker sufficient conditions
\[
 s>3\ (d=2),\qquad s\ge s_3=3.8264\ldots\ (d=3),\qquad s\ge\tfrac{5+\sqrt{23}}2\ (d=4),
\]
and $s>d/2+2$ with $2s^2-(d+6)s+1\ge0$ for $d\ge5$. Here $s_3$ is the root above $3.5$ of $4s^3-8s^2-29s+4$. These sufficient thresholds follow from control of the nonlinear remainder of the one-step estimator (Section~\ref{sec:ideas-upper}). The lower bound holds for every $s>d/2+1$ and uses only $C^\infty$ alternatives in $\cI_s(\eta)$.
\end{remark}

The pilot in the one-step estimator is constructed in Lemma~\ref{up:lem:pilot} as a measurable near-minimiser over a countable dense subset of the class. A computationally efficient implementation is a direction for further work. The upper bound is proved in Section~\ref{up:sec:iso}, and the lower bound in Section~\ref{low:sec:iso}. If instead \condC\ holds at \emph{every} member of a class $\mathcal G\subset\cF_s(f_0,\eta)$ and $s>d/2+2$, the same one-step estimator with the exact representer attains a uniform $O(N^{-1/2})$ risk bound over $\mathcal G$, without any isotropy assumption (Proposition~\ref{up:prop:uniformC}). That \condC\ at the reference does not suffice is a nonlinear phenomenon:

\begin{proposition}[Instability of isotropy and of condition (C)]\label{prop:instability}
Consider the canonical model.
\begin{enumerate}
\item[(a)] For every small $\varepsilon>0$ there is a smooth positive $f_\varepsilon$, equal to $1$ near $0$, with $\norm{f_\varepsilon-1}_{H^s}\lesssim_s\varepsilon$ for every $s$, critical point $0$ and $\Teff(D^2u_{f_\varepsilon}(0))=d/(1+2d\varepsilon)<d$.
\item[(b)] If $\psi\ne0$ satisfies \condC\ at $f_0=1$, there are smooth conductivities arbitrarily close to $1$ in every $H^s$ norm, with the same critical point and Hessian $I_d/d$, at which \condC\ fails.
\end{enumerate}
Consequently, for $d\le4$, every nonzero $\psi\in C_c^\infty(B_1)$ and every $\eta>0$, no estimator attains the rate $N^{-1/2}$ uniformly over $\cF_s(1,\eta)$. If $d\in\{2,3\}$ and $\psi(0)\ne0$, the rate over $\cF_s(1,\eta)$ is polynomially slower than $\rate(s,d)$, by part (a) and Theorem~\ref{thm:lower}(b) applied at $f_\varepsilon$.
\end{proposition}

The first consequence already follows from Theorem~\ref{thm:isotropic}(ii), since $\cI_s(\eta)\subset\cF_s(1,\eta)$; part (b) explains why. Proposition~\ref{prop:instability} is proved in Section~\ref{low:sec:instability}.

\subsection{Full neighbourhoods}

\begin{theorem}[Parametric rate in the regular regime]\label{thm:full-regular}
Let Assumption~\ref{ass:base} hold with $\Teff_0>4$. If $s>d/2+2$ and $s^2-(d+1)s+(d-2)/2\ge0$ (in particular if $s>d+1$), then for all sufficiently small $\eta>0$ and every nonzero $\psi$,
\[
 \Risk(\cF_s(f_0,\eta))\asymp N^{-1/2}.
\]
If $H$ is isotropic (so $d\ge5$), $s>d/2+3$ suffices.
\end{theorem}

Above the elbow the regular regime is therefore stable: since $\Teff$ depends continuously on $f$, every conductivity in a small neighbourhood is again regular, and the rate is uniform. The quadratic condition on $s$ again comes from the nonlinear remainder. Theorem~\ref{thm:full-regular} is proved in Section~\ref{up:sec:full-regular}.

Below the elbow, full neighbourhoods contain conductivities with different effective dimensions. Theorem~\ref{thm:full-bracket} in the Supplement gives an explicit upper bound over classes on which the condition number of the critical Hessian is at most $\kappa\in(1,d)$. Its polynomial exponent $\beta_\kappa=(s+d/2)/(2(s+2\kappa))$ tends to $(2s+d)/(4s+8)$ as $\kappa\downarrow1$, provided the nonlinear remainder exponent
\[
 \gamma_1=\frac{2s-2+1/(2s)}{2s+d+2}
\]
exceeds $(s+d/2)/(2s+4)$. Combining this with the lower bounds identifies the local rate exponent at the canonical model. For a class $\cF$ let $e(\cF)=\liminf_{N\to\infty}\log(1/\Risk(\cF))/\log N$.

\begin{corollary}[Local exponent at the canonical model]\label{cor:local-exponent}
Consider the canonical model with $d\in\{2,3,4\}$, and let $s>d/2+2$ satisfy $\gamma_1>(s+d/2)/(2s+4)$. Then for every nonzero $\psi$,
\[
 \lim_{\eta\downarrow0}e\bigl(\cF_s(1,\eta)\bigr)=\frac{2s+d}{4s+8}\quad(d=2,3),\qquad \lim_{\eta\downarrow0}e\bigl(\cF_s(1,\eta)\bigr)=\frac12\quad(d=4).
\]
Here the limit $N\to\infty$ is taken first, for fixed $\eta$. For each fixed $\eta>0$, $d\in\{2,3\}$ and $\psi(0)\ne0$, $e(\cF_s(1,\eta))<(2s+d)/(4s+8)$.
\end{corollary}

The minimax rate on a fixed full neighbourhood of the canonical model is thus strictly slower than $\rate(s,d)$. The gap in the polynomial rate exponent tends to zero as the neighbourhood radius shrinks. Corollary~\ref{cor:local-exponent} is proved in Section~\ref{up:sec:local}. For anisotropic references with $\Teff_0\le4$, the minimax rate on full neighbourhoods is open; see Conjecture~\ref{conj:anisotropic}.

\begin{remark}[Interior perturbations]\label{rem:collar}
Suppose the parameter class only allows perturbations $h$ supported in a fixed compact set $K\subset\cO$. In the linearised experiment, $\psi(x_*)=0$ then suffices for the parametric rate, because the characteristic solution minus a function constant along flow lines can be cut off outside $K$. If $x_*\in\operatorname{int}K$ and $\psi(x_*)\ne0$, Theorem~\ref{thm:lower}(a) still applies (Proposition~\ref{low:collar}). For isotropic Hessians the corresponding radial constructions are supported in a small ball around $x_*$, so the rate $\rate(s,d)$ persists.
\end{remark}
\section{Main ideas of the proofs}\label{sec:ideas}

\subsection{The representer and the effective dimension}\label{sec:ideas-linear}
Consider the quadratic model $u_0(x)=\frac12x^{\top}Hx-c$ with constant $f_0$, $x_*=0$ and $H=\operatorname{diag}(\lambda_1,\dots,\lambda_d)$. Put $a_i=\lambda_i/\lmax$ and $D_r=\operatorname{diag}(r^{a_1},\dots,r^{a_d})$, so that $\det D_r=r^{\Teff}$. The gradient flow maps the flow-adapted boxes $D_r[-1,1]^d$ onto each other. After subtraction of a smooth particular solution, the characteristic solution has a singular part $F$ satisfying
\[
 F(D_rx)=F(x)+\frac{\psi(x_*)}{\lmax}\log r .
\]
Thus $F$ is constant if and only if \condC\ holds. If $F$ is not constant, its second derivatives in the fastest coordinates do not vanish identically. They scale like $r^{-2}$ on the shell $D_r(\overline{B_2}\setminus B_1)$ of volume $r^{\Teff}$, which gives the integral $\int r^{-4}r^{\Teff-1}dr$ of Section~\ref{sec:intro-results} and the representer criterion in Theorem~\ref{thm:linear}.

\begin{example}[The disc]\label{ex:disc}
In the canonical model with $d=2$ we have $\nabla u_0(x)=x/2$. For a radial $\psi$ the profile $\Psi$ is constant, so \condC\ holds if and only if $\psi(0)=0$. Writing $r=|x|$, the transport equation \eqref{eq:transport} becomes $\frac r2w'(r)=\psi(r)$ with $w(1)=0$, so that
\[
 w(x)=2\psi(0)\log|x|+w_1(x),\qquad w_1(x)=-2\int_{|x|}^1\frac{\psi(t)-\psi(0)}{t}\,dt,
\]
with $w_1$ smooth. Since $\Delta\log|x|=2\pi\delta_0$ in $\R^2$, the candidate representer is
\[
 z=\Delta w=4\pi\psi(0)\,\delta_0+\Delta w_1 .
\]
It is square integrable if and only if $\psi(0)=0$. Otherwise we replace $\delta_0$ by a radial mollifier $\delta_\rho=\rho^{-2}\chi(\cdot/\rho)$. By Newton's theorem $\Delta^{-1}(\delta_0-\delta_\rho)$ is supported in $B_\rho$, so the residual $\psi-A^*z_\rho$ is a bounded function supported in $B_\rho$. After cancelling finitely many of its moments, it pairs with $h\in H^s$ to give a bias of order $\rho^{s+1}\norm h_{H^s}$. The variance is $N^{-1}\norm{\delta_\rho}_2^2\asymp N^{-1}\rho^{-2}$. Balancing at $\rho=N^{-1/(2s+4)}$ gives $N^{-(s+1)/(2s+4)}=\rate(s,2)$. Conversely, a conductivity bump of width $\rho$ at the origin changes $\Phi$ by about $\rho^{s+1}$ and $u_f$ by about $\rho^{s+2}$ in $L^2$ (Section~\ref{sec:ideas-lower}), so the rate cannot be improved. In this example the irregularity is entirely due to the point mass that the logarithm of the flow time produces at the minimum.
\end{example}

For the upper bound, the representer is regularised by an isotropic Fourier cutoff at frequency $R$. The key estimate is that the Fourier mass of the localised singular part on the annulus $\{R\le|\xi|\le2R\}$ is of order $R^{-\Teff}$, although $F$ is only anisotropically homogeneous. This yields the residual bound $\norm{\psi-A^*z_R}_{H^{-s}}\lesssim R^{-s-\Teff/2}$ and the variance bound $\norm{z_R}_2^2\lesssim R^{4-\Teff}$ used in the heuristic above. For general references, a resonant normal form \cite{sternberg1957local,chen1963equivalence} reduces the local flow to a triangular system in which the fastest coordinates enter linearly, and this preserves both scalings (Section~\ref{sec:normalform}).

\subsection{Lower bounds through the forward map}\label{sec:ideas-lower}
The natural first attempt is a bump $h(x)=a\,\chi(x/\rho)$ in the conductivity at $x_*$. Its forcing $T_{u_0}h=\diver(h\nabla u_0)$ has vanishing mass. Writing $\nabla u_0(x)=Hx+O(|x|^2)$, the contribution $\diver(hHx)$ also has vanishing first moments when $\chi$ is even. Duality with $H^2\cap H^1_0$ and the Bramble--Hilbert lemma \cite{bramble1970estimation}, with the $O(|x|^2)$ remainder controlled separately in the dual $H^2$ norm, then show that $\norm{A_{f_0}h}_2\lesssim\rho^2\norm{h}_2$, while $\ip h\psi\asymp\rho^{d/2}\norm h_2\,\psi(x_*)$. Balancing these against the $H^s$ constraint yields Theorem~\ref{thm:lower}(a), with $d$ in place of $\Teff$. For anisotropic $H$ this is not optimal: a bump is not adapted to the flow, and its transport along the flow lines is not exploited.

To obtain $\Teff$ we prescribe the change of the potential rather than that of the conductivity. In the quadratic model, let $q\in C_c^\infty(\R^d\setminus\{0\})$, $b>0$, $r\in(0,1)$, and set
\[
 V_r(x)=b\,q(D_r^{-1}x),\qquad h_r=-S_{u_0}L_{f_0}V_r .
\]
Since $T_{u_0}S_{u_0}$ is the identity, the definition of $A_{f_0}$ gives
\[
 A_{f_0}h_r=-L_{f_0}^{-1}T_{u_0}h_r=L_{f_0}^{-1}L_{f_0}V_r=V_r .
\]
The conductivity perturbation $h_r$ is not localised: it has a transport tail along the flow lines leaving the support of $V_r$. Nevertheless $S_{u_0}$ is bounded on $H^s$ (Lemma~\ref{low:S}), and the three relevant quantities scale as follows.
\begin{enumerate}
\item[(1)] \emph{Size.} $L_{f_0}V_r=f_0\Delta V_r$ consists of second derivatives of $q(D_r^{-1}\cdot)$, the largest of which are in the fastest coordinates and of order $r^{-2}$; $s$ further derivatives cost at most $r^{-s}$, and the support has volume $r^{\Teff}$. Hence $\norm{h_r}_{H^s}\lesssim b\,r^{\Teff/2-s-2}$.
\item[(2)] \emph{Signal.} $\norm{V_r}_2=b\,r^{\Teff/2}\norm q_2$.
\item[(3)] \emph{Separation.} Write $\psi=b\cdot\nabla w$ with the characteristic solution $w$ of \eqref{eq:transport} and $b=\nabla u_0$. Integrating by parts twice, using that $w$ vanishes on $\partial\cO$ and $V_r$ is compactly supported,
\[
 \ip{h_r}\psi=-\ip{T_{u_0}h_r}{w}=\ip{L_{f_0}V_r}{w}=f_0\,b\,r^{\Teff}\int q(y)\,(\Delta w)(D_ry)\,dy .
\]
By the scaling of $F$, $(\Delta w)(D_ry)=r^{-2}(\Delta_EF)(y)+O(r^{-2+\delta})$ for some $\delta>0$, where $\Delta_E$ is the Laplacian in the fastest coordinates. Whenever $\Delta_EF\not\equiv0$ away from the origin, we may choose $q$ with $\int q\,\Delta_EF\ne0$. Then $|\ip{h_r}\psi|\asymp b\,r^{\Teff-2}$.
\end{enumerate}
The exceptional case is $d=2$ with isotropic $H$ and $F=c\log|x|+\mathrm{constant}$, $c\ne0$. Then we take $q\in C_c^\infty(\R^2)$ with $q(0)\ne0$ and use the distributional pairing $\langle F,\Delta q\rangle=2\pi c q(0)$; the same separation scaling follows, as in the proof of Theorem~\ref{thm:lower}(b).
Now take $b=\gamma r^{s+2-\Teff/2}$ with $\gamma>0$ small and $r=N^{-1/(2s+4)}$. Then $\norm{h_r}_{H^s}\lesssim\gamma\le\eta$, so $f_0+h_r\in\cF_s(f_0,\eta)$, and the functional moves by $|\ip{h_r}\psi|\asymp\gamma r^{s+\Teff/2}=\gamma N^{-(2s+\Teff)/(4s+8)}$.

The step from the linearised to the nonlinear experiment is exact. By \eqref{eq:exact} and $T_{u_0}h_r=-L_{f_0}V_r$,
\[
 u_{f_0+h_r}-u_0=L_{f_0+h_r}^{-1}L_{f_0}V_r=V_r-L_{f_0+h_r}^{-1}\diver(h_r\nabla V_r).
\]
By duality with $H^2\cap H^1_0$, uniform elliptic regularity \cite{gilbarg2001elliptic} and the multiplier estimate $\norm{\diver(h\nabla\varphi)}_2\lesssim\norm h_{H^s}\norm\varphi_{H^2}$ for $s>d/2$, the last term is bounded in $L^2$ by a constant times $\norm{h_r}_{H^s}\norm{V_r}_2\lesssim\eta\norm{V_r}_2$. There is no Taylor remainder, so
\[
 \KL\bigl(P_{f_0+h_r},P_{f_0}\bigr)=\frac N2\norm{u_{f_0+h_r}-u_0}_2^2\lesssim N\norm{V_r}_2^2\asymp\gamma^2Nr^{2s+4}=\gamma^2 .
\]
For $\gamma$ small the two experiments cannot be distinguished, whereas $\Phi$ differs by a multiple of $N^{-(2s+\Teff)/(4s+8)}$, and Le Cam's two-point argument \cite{tsybakov2009introduction} (Lemma~\ref{lin:lem:twopoint}) gives Theorem~\ref{thm:lower}(b) for $\Teff_0<4$. For $\Teff_0=4$, each flow-adapted shell contributes the same separation and the same KL divergence; about $\log N$ perturbations on disjoint shells are added. Their contributions to the KL divergence add up, while their separations add linearly, which gives the separation $\sqrt{\log N/N}$ at bounded divergence. For general references the same construction is carried out in normal-form coordinates (Proposition~\ref{lin:prop:construction}).

On the isotropic class the potential perturbation is chosen constant near $0$, so that the critical point and the Hessian are preserved exactly and the alternative lies in $\cI_s(\eta)$. When \condC\ holds at the reference, we first move to a nearby conductivity in the class at which \condC\ fails (Proposition~\ref{prop:instability}(b)) and then transplant the construction. Full details are in Section~\ref{sec:lower-proofs} of the Supplement.

\subsection{The one-step upper bound}\label{sec:ideas-upper}
All upper bounds use a prediction pilot and a one-step correction. Given a pilot $p$ from the first half of the observation and $w\in H^2(\cO)\cap H^1_0(\cO)$ determined by $p$, the second half yields
\[
 \widehat\Phi=\ip p\psi+Y_2(L_pw)-\ip{u_p}{L_pw}.
\]
Write $h=f-p$, $v=u_f-u_p$ and $r=\psi-\nabla u_p\cdot\nabla w$. Lemma~\ref{up:lem:identity} gives
\[
 \E_f|\widehat\Phi-\Phi(f)|
 \lesssim
 \E_f|\ip hr|
 +N^{-1/2}\E_f\norm{L_pw}_2
 +\E_f\Bigl|\int_\cO h\nabla v\cdot\nabla w\Bigr| .
\]
These terms are the regularisation bias, the noise contribution and the nonlinear remainder. The geometry enters through the choice of $w$.

On isotropic classes we regularise the characteristic solution of $p$ in a ball of radius $\rho$ around the pilot's critical point. Subtracting the value of the logarithmic term at radius $\rho$ keeps the residual bounded independently of $\rho$. Small support alone does not give the required negative Sobolev estimate when $s>d/2$, so we also cancel finitely many residual moments. This yields $\norm{r}_{H^{-s}}\lesssim\rho^{s+d/2}$ and, for $d\le3$, $\norm{L_pw}_2^2\lesssim\rho^{d-4}$. When the exact representer is uniformly square integrable, we use it directly and $r=0$.

Combining the three terms gives, uniformly over the isotropic class,
\[
 \E_f|\widehat\Phi_\rho-\Phi(f)|\lesssim\rho^{s+d/2}+N^{-1/2}V_d(\rho)^{1/2}+N^{-\gamma(s,d)},
\]
\[
 V_d(\rho)=\begin{cases}\rho^{d-4},&d\le3,\\ 1+\log(1/\rho),&d=4,\\ 1,&d\ge5,\end{cases}
\]
where the last term is the remainder discussed next. With $\rho=N^{-1/(2s+4)}$ the first two terms are of order $\rate(s,d)$. In the radial $d=2$ example, the variance comes from the regularised point mass of Example~\ref{ex:disc}; nonconstant flow profiles also generate angular singularities. In $d=3$ it comes from the $r^{-2}$ singularity of $\Delta w$, and in $d=4$ from its logarithmic divergence.

The remainder is bounded through $|\nabla w|\lesssim|x-x_p|^{-1}$ and a fractional Hardy inequality by $\norm h_{H^{1/2}}\norm v_{H^{3/2}}$. The pilot's prediction rate and local transport stability then give the remainder exponent
\[
 \gamma(s,d)=\frac{(2s-1)^2}{2s(2s+d+2)},
\]
which must be at least the exponent of $\rate(s,d)$. This requirement determines the thresholds in Remark~\ref{rem:thresholds}. The regularisation is centred at the pilot's critical point, and the exact identity accounts for the discrepancy between the flows of $u_p$ and $u_f$, so no separate estimate of the location error is needed.
\section{Numerical illustration}\label{sec:numerics}

We illustrate Theorem~\ref{thm:linear} by deterministic computations of the modulus of continuity
\begin{equation}\label{num:eq:modulus}
 \omega(\varepsilon)=\sup\bigl\{\ip h\psi:\ \norm{A_{f_0}h}_2\le\varepsilon,\ h\in\Theta\bigr\}
\end{equation}
in the linearised experiment at a fixed reference. The minimax risk for $\ip h\psi$ over a convex symmetric set $\Theta$ is of order $\omega(N^{-1/2})$ \cite{ibragimov1984nonparametric,donoho1994statistical}, so $\omega(\varepsilon)\asymp\varepsilon^{2e}$ corresponds to the rate $N^{-e}$. In the irregular and regular regimes Theorem~\ref{thm:linear} predicts
\begin{equation}\label{num:eq:pred}
 e_{\mathrm{pred}}(s,\Teff)=\min\Bigl\{\frac{2s+\Teff}{4s+8},\frac12\Bigr\},
\end{equation}
with an additional factor $\sqrt{\log(1/\varepsilon)}$ in the critical regime $\Teff=4$, and $e_{\rm pred}=1/2$ when \condC\ holds. We compute $\omega$ on balls in $\R^d$, $d=2,\dots,6$ (isotropic Hessian, $\Teff=d$), and on ellipsoids with $f_0=g=1$, whose Hessian has eigenvalues proportional to $(\mu,\dots,\mu,1)$ in $d=2,3,5$, so that $\Teff$ ranges over $(1,6]$. Here $\Theta$ is the unit ball of a spectral or polynomial Hilbert norm of order $s$. We report one selected local exponent of $\omega$ per case, read off at the smallest scales that are free of truncation effects. The discretisations, the selection rule and the complete list of values are given in Section~\ref{num:app} of the supplement.

\begin{figure}[t]
\centering
\includegraphics[width=0.9\textwidth]{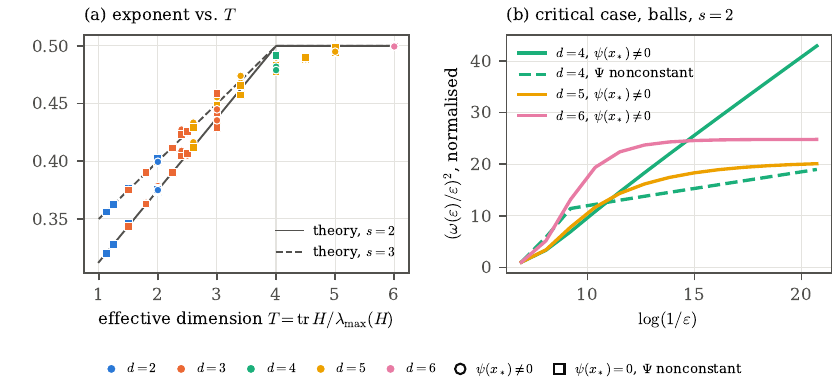}
\caption{(a) Selected local exponents against the effective dimension $\Teff$, with the prediction \eqref{num:eq:pred} for $s=2$ (solid) and $s=3$ (dashed). Colour gives the dimension $d$, the marker the type of target. (b) The critical regime: $(\omega(\varepsilon)/\varepsilon)^2$ for the isotropic balls, $s=2$, each curve normalised by its value at $\varepsilon=10^{-3}$, so that only the shapes are compared.}
\label{num:fig:exponents}
\end{figure}

Figure~\ref{num:fig:exponents}(a) shows that the selected exponents follow \eqref{num:eq:pred} as a function of $\Teff$, rather than of $d$. For strongly anisotropic Hessians the difference from \eqref{num:eq:pred} is at most $0.002$, while replacing $\Teff$ by $d$ would predict exponents larger by up to $0.13$. In $d=5$ with $\Teff<4$ the exponents stay below $1/2$, in line with Corollary~\ref{cor:linear}. In the critical regime $\Teff=4$, $(\omega/\varepsilon)^2$ grows linearly in $\log(1/\varepsilon)$, whereas it saturates for $d=5,6$ (Figure~\ref{num:fig:exponents}(b)). Targets satisfying \condC\ give exponents of at least $0.4965$. The parameter sets are Hilbert balls with boundary behaviour different from that of ordinary Sobolev balls, and the selected exponents are finite-resolution slopes.

\section{Discussion}\label{sec:discussion}

The elbow of the linearised rate is at $\Teff=4$ for every smoothness $s$. In the linearised experiment the rate of a smooth linear functional is determined by condition \condC\ and by one number attached to the critical point of the potential, the effective dimension $\Teff=\tr H/\lmax(H)$. The relevant length scale is that of the fastest flow direction, while the volume available for hiding a perturbation is governed by all eigenvalues. In the irregular regime the exponent $(2s+\Teff)/(4s+8)$ lies strictly between that of the ill-posed problem of estimating $f$ and the parametric exponent $1/2$.

In linear inverse problems the classification of functionals depends only on the fixed forward operator. In the Darcy problem the forward operator, and with it the effective dimension and condition \condC, change with $f$, and Proposition~\ref{prop:instability} shows that this matters for uniform statements. A semiparametric theory based on \condC\ therefore needs a parameter class on which \condC\ holds at every member (Proposition~\ref{up:prop:uniformC}), or a local asymptotic formulation with shrinking neighbourhoods; fixed Sobolev balls in $d\le4$ do not qualify. Establishing efficiency, asymptotic normality and Bernstein--von Mises limits requires further distributional analysis beyond the parametric-rate bounds.

The main open question is the minimax rate on full neighbourhoods of anisotropic references below the elbow. Theorem~\ref{thm:lower}(b) and Theorem~\ref{thm:linear} suggest the following.

\begin{conjecture}\label{conj:anisotropic}
Under Assumption~\ref{ass:base} with $\Teff_0\le4$, if \condC\ fails at $f_0$, then for all sufficiently large $s$,
\[
 \lim_{\eta\downarrow0}e\bigl(\cF_s(f_0,\eta)\bigr)=\min\Bigl\{\frac{2s+\Teff_0}{4s+8},\frac12\Bigr\}.
\]
\end{conjecture}

The upper bound of Theorem~\ref{thm:full-bracket} uses isotropic truncation. A flow-adapted truncation, as in the linearised analysis, would have to be made uniform over the pilot and compatible with the Hardy remainder. Further open problems include:
\begin{itemize}
\item the minimal smoothness, since the conditions in Theorems~\ref{thm:isotropic} and~\ref{thm:full-regular} are only sufficient for our estimators;
\item several, degenerate or saddle-type critical points, which arise for instance for $g$ of variable sign;
\item the minimax rate under interior support constraints or known boundary values, for anisotropic references with $\psi(x_*)\ne0$;
\item computationally efficient pilots, for instance based on the method of \cite{wang2026global};
\item Bernstein--von Mises theorems at the nonparametric rate $\rate(s,\Teff)$.
\end{itemize}

\clearpage
\begin{supplement}
\stitle{Supplement to ``Minimax estimation of linear functionals in the Darcy inverse problem''}
\sdescription{The supplement \cite{liu2026supp} contains the proofs of all results, auxiliary lemmas, an upper bound under controlled anisotropy (Theorem~\ref{thm:full-bracket}) and details of the numerical computations. Sections, results and equations whose numbers carry the prefix S refer to the supplement.}
\end{supplement}

\section*{Supplementary proofs and details}
This supplement contains the proofs of all results stated in the main paper, the auxiliary lemmas they use, an explicit upper bound under controlled anisotropy, and the details of the numerical computations. Equations, theorems and sections without the prefix S refer to the main paper.

\setcounter{section}{0}
\setcounter{figure}{0}
\setcounter{table}{0}
\renewcommand{\thesection}{S\arabic{section}}
\renewcommand{\thefigure}{S\arabic{figure}}
\renewcommand{\thetable}{S\arabic{table}}
\renewcommand{\theHsection}{supp.\arabic{section}}
\renewcommand{\theHfigure}{supp.\arabic{figure}}
\renewcommand{\theHtable}{supp.\arabic{table}}

\section{Overview and an upper bound under controlled anisotropy}\label{sup:sec:overview}

Section~\ref{sec:linear-proofs} proves the results on the linearised experiment (Theorem~\ref{thm:linear}, Corollary~\ref{cor:linear} and Proposition~\ref{prop:C-smooth}). Section~\ref{sec:lower-proofs} proves the nonlinear lower bounds (Theorem~\ref{thm:lower}, the lower bound in Theorem~\ref{thm:isotropic}, Proposition~\ref{prop:instability} and Remark~\ref{rem:collar}). Section~\ref{sec:upper-proofs} proves the nonlinear upper bounds (the upper bound in Theorem~\ref{thm:isotropic}, Theorem~\ref{thm:full-regular}, Theorem~\ref{thm:full-bracket} below and Corollary~\ref{cor:local-exponent}). The remaining sections contain auxiliary lemmas and the details of the numerical computations of Section~\ref{sec:numerics}. Notation is as in Sections~\ref{sec:intro} and~\ref{sec:setting} of the main paper.

The local exponent in Corollary~\ref{cor:local-exponent} is derived from the following explicit bound. Its exponents are not rate-optimal for $\kappa>1$, and its constants are not uniform as $\kappa\downarrow1$.

\begin{theorem}[Upper bound under controlled anisotropy]\label{thm:full-bracket}
Let Assumption~\ref{ass:base} hold, let $d\in\{2,3,4\}$, $s>d/2+2$ and $1<\kappa<d$, and let $\eta>0$ be sufficiently small. Let $\mathcal H\subset\cF_s(f_0,\eta)$ be a class on which, for some $m>0$ and a fixed radius $r_0>0$,
\[
 mI\le D^2u_f\le\kappa mI\qquad\text{on }B_{r_0}(x_*(f)),
\]
uniformly in $f\in\mathcal H$. Then there is an estimator, depending on $\kappa$, with
\[
 \sup_{f\in\mathcal H}\E_f|\widehat\Phi-\Phi(f)|\lesssim N^{-\beta_\kappa}(\log N)^{\vartheta_\kappa}+N^{-\gamma_\kappa},
\]
where
\[
 \beta_\kappa=\frac{s+d/2}{2(s+2\kappa)},\qquad \vartheta_\kappa=\frac{2\kappa-d/2}{s+2\kappa},\qquad \gamma_\kappa=\frac{2s-1-\kappa+\kappa/(2s)}{2s+d+2}.
\]
The constant depends on $\kappa$, $m$, $r_0$, the geometry of the reference and the Sobolev radius, but not on $f$ or $N$.
\end{theorem}

Theorem~\ref{thm:full-bracket} is proved in Section~\ref{up:sec:bracket}. As $\kappa\downarrow1$, $\beta_\kappa\to(2s+d)/(4s+8)$ and $\gamma_\kappa\to\gamma_1=(2s-2+1/(2s))/(2s+d+2)$.
\section{Proofs for the linearised experiment}\label{sec:linear-proofs}

This section proves Theorem~\ref{thm:linear}, Corollary~\ref{cor:linear} and the equivalence between condition \condC\ and regularity of the transport solution. In the linearised experiment the risk of the linear estimator $Y(z)$ is the sum of a squared bias $\ip h{A_{f_0}^*z-\psi}^2$ and a variance $N^{-1}\norm z_2^2$; see \eqref{lin:eq:linest}. Everything therefore reduces to the singular part $F$ of the characteristic solution $w$ near $x_*$, and the proof has three quantitative ingredients.
\begin{enumerate}
\item[(i)] A scaling law. Second derivatives of $F$ in the fastest flow direction scale like $r^{-2}$ on flow shells of volume $r^{\Teff}$ (Lemma~\ref{lin:lem:fast} and \eqref{lin:eq:hom-deriv}). This decides whether $w\in H^2$, hence whether an $L^2$ representer exists (Proposition~\ref{lin:prop:range-quad}).
\item[(ii)] A Fourier shell bound $\int_{R\le|\xi|\le2R}|\widehat{\zeta F}|^2\lesssim R^{-\Teff}$ (Lemmas~\ref{lin:lem:band} and~\ref{lin:lem:general-band}). Truncating $F$ in frequency at $R$ gives bias $R^{-s-\Teff/2}$, through a commutator (Lemma~\ref{lin:lem:comm}), and variance $R^{4-\Teff}$; with $R=N^{1/(2s+4)}$ this is the upper bound.
\item[(iii)] Flow-adapted alternatives $h=-S_{u_0}L_{f_0}V$ with $A_{f_0}h=V$ exactly (Proposition~\ref{lin:prop:construction}). Their $H^s$ norm, $L^2$ image and pairing with $\psi$ scale as $e^{(s+2)\lambda_*t-Qt/2}$, $e^{-Qt/2}$ and $e^{(2\lambda_*-Q)t}$, which gives the matching lower bound.
\end{enumerate}
Section~\ref{lin:sec:flow} collects facts on the gradient flow and the transport inverse. Section~\ref{lin:sec:gauss} records two elementary facts on Gaussian experiments. Section~\ref{lin:sec:quadratic} gives the complete argument for an exactly quadratic reference, where $F$ is quasi-homogeneous. Section~\ref{sec:normalform} extends it to general references, including resonant ones. Routine flow estimates are deferred to Section~\ref{lin:app}.

Throughout this section $f_0,g,u_0,x_*,H$ are as in Assumption~\ref{ass:base}, $b=\nabla u_0$, $Q=\tr H$, $\lambda_*=\lmax(H)$, and
\[
 a_i=\lambda_i/\lambda_*\in(0,1],\qquad \Teff=\Teff_0=\sum_{i=1}^d a_i=Q/\lambda_* .
\]
The linear target is $\ell(h)=\ip h\psi$.
\subsection{Gradient flow, transport inverse and the characteristic solution}\label{lin:sec:flow}

Since $g>0$, the maximum principle gives $u_0<0$ in $\cO$, and the Hopf lemma gives $\partial_\nu u_0>0$ on $\partial\cO$, where $\nu$ is the outer normal. Let $\phi_t$ denote the flow of $b$ and $X_t=\phi_{-t}$ the backward flow. Along backward trajectories $u_0$ decreases at rate $|b|^2$. Since $x_*$ is the only critical point and it is a nondegenerate minimum, we obtain the following facts; the proof is in Section~\ref{lin:app}.

\begin{lemma}\label{lin:lem:flow}
\begin{enumerate}
\item[(i)] $X_t$ maps $\overline\cO$ into itself for $t\ge0$, and $X_tx\to x_*$ as $t\to\infty$, uniformly and exponentially: $|X_tx-x_*|\le Ce^{-\mu t}$ for some $\mu>0$. Every derivative of fixed order $k\ge1$ satisfies $\norm{D^kX_t}_\infty\le C_ke^{-\mu t}$.
\item[(ii)] For $x\ne x_*$ the forward trajectory leaves $\cO$ transversally at a finite time $\tau(x)$, with exit point $\vartheta(x)=\phi_{\tau(x)}x\in\partial\cO$. The functions $\tau,\vartheta$ are smooth on $\overline\cO\setminus\{x_*\}$, and $\tau(x)\to\infty$ as $x\to x_*$.
\item[(iii)] The Jacobian $J_t=\det DX_t$ satisfies $J_t(x)=\exp(-\int_0^t\Delta u_0(X_rx)\,dr)$, $J_t\le Ce^{-ct}$ for some $c>0$, and $|D^kJ_t|\le C_kJ_t$ for every $k$.
\end{enumerate}
\end{lemma}

For a smooth potential $u$ with these properties, let $T_uh=\diver(h\nabla u)$ and define the \emph{transport inverse}
\begin{equation}\label{lin:eq:S}
 S_uq(x)=\int_0^\infty J_t(x)\,q(X_tx)\,dt .
\end{equation}

\begin{lemma}\label{lin:lem:S}
For $q\in L^2(\cO)$ we have $T_uS_uq=q$ in the sense of distributions. For every $s\ge0$, $S_u$ is bounded on $H^s(\cO)$. If $q$ is supported in a compact set $K$, then $S_uq$ vanishes outside the forward orbit $\bigcup_{t\ge0}\phi_t(K)$.
\end{lemma}

The function $S_uq$ is generally nonzero on the whole forward orbit of $\supp q$, which reaches $\partial\cO$. This ``transport tail'' is part of every perturbation built from $S_u$, and Lemma~\ref{lin:lem:S} controls it in $H^s$. The proof is in Section~\ref{lin:app}.

\paragraph*{The characteristic solution}
Put $\psi_*=\psi(x_*)$ and define
\begin{equation}\label{lin:eq:vpsi}
 v_\psi(x)=\int_0^\infty\{\psi(X_tx)-\psi_*\}\,dt .
\end{equation}
By Lemma~\ref{lin:lem:flow}(i), the integrand and all its spatial derivatives decay exponentially, so $v_\psi\in C^\infty(\overline\cO)$, $v_\psi(x_*)=0$, and $b\cdot\nabla v_\psi=\psi-\psi_*$; the last identity follows by differentiating $v_\psi(X_sx)$ in $s$. The function
\begin{equation}\label{lin:eq:w}
 w(x)=-\int_0^{\tau(x)}\psi(\phi_tx)\,dt=v_\psi(x)-v_\psi(\vartheta(x))-\psi_*\,\tau(x),\qquad x\ne x_*,
\end{equation}
is smooth on $\overline\cO\setminus\{x_*\}$, vanishes on $\partial\cO$, and satisfies $b\cdot\nabla w=\psi$ there. Indeed $b\cdot\nabla\tau=-1$, and $v_\psi\circ\vartheta$ is constant along trajectories. We call $w$ the \emph{characteristic solution} of \eqref{eq:transport}. Its singular part
\begin{equation}\label{lin:eq:F}
 F=w-v_\psi=-v_\psi\circ\vartheta-\psi_*\tau
\end{equation}
is defined only on $\overline\cO\setminus\{x_*\}$; as it is locally bounded by $C(1+|\log|x-x_*||)$, we regard it as an element of $L^2_{\rm loc}(\cO)$, and its value at the single point $x_*$ is irrelevant. It satisfies $b\cdot\nabla F=\psi_*$ off $x_*$. Since $\tau(x)\le C(1+|\log|x-x_*||)$ and $|b(x)|\le C|x-x_*|$, removing a ball $B_\epsilon(x_*)$ produces boundary terms of order $\epsilon^d(1+|\log\epsilon|)$. Hence $b\cdot\nabla F=\psi_*$ and $b\cdot\nabla w=\psi$ hold in the sense of distributions on $\cO$, and $w\in L^2(\cO)$.

When $\psi_*=0$, the profile \eqref{eq:profile} is $\Psi(\ell)=-v_\psi(\vartheta(x))$ for any $x$ on $\ell$. Indeed, letting $x\to x_*$ backward along $\ell$ in \eqref{lin:eq:w} gives $w(x)\to\Psi(\ell)$ and $v_\psi(x)\to0$. Thus
\begin{equation}\label{lin:eq:C-F}
 \text{\condC\ holds}\iff\psi_*=0\text{ and }v_\psi|_{\partial\cO}\text{ is constant}\iff F\text{ is constant}.
\end{equation}
For the last equivalence, note that $\tau$ is unbounded near $x_*$, while $v_\psi\circ\vartheta$ is bounded; so $F$ constant forces $\psi_*=0$.

\begin{remark}\label{lin:rem:clock}
When $\psi_*\ne0$, any regularised profile of the form $-\int_{-\infty}^{0}\{\psi(\gamma(t))-\psi_*\}\,dt-\int_0^{t_\partial}\psi(\gamma(t))\,dt$ depends on the choice of a time origin on each line. Shifting that origin by $a(\ell)$ changes the profile by $\psi_*a(\ell)$. This is why Definition~\ref{def:C} requires $\psi_*=0$ first. With this convention, \condC\ is invariant under changes of coordinates and of transverse sections.
\end{remark}

\begin{proposition}\label{prop:C-smooth}
\begin{enumerate}
\item[(i)] Condition \condC\ holds if and only if the characteristic solution $w$ extends continuously to $x_*$. In that case $w=v_\psi-c$, where $c$ is the constant boundary value of $v_\psi$, so $w\in C^\infty(\overline\cO)$.
\item[(ii)] For a smooth positive $f_0$, the equation $A_{f_0}^*z=\psi$ has a solution $z\in L^2(\cO)$ if and only if
\begin{equation}\label{lin:eq:range}
 \psi\in b\cdot\nabla\bigl(H^2(\cO)\cap H^1_0(\cO)\bigr),
\end{equation}
and this holds if and only if $w\in H^2(\cO)$. In that case $z=L_{f_0}w$, and this solution is unique.
\item[(iii)] If $d\le3$, an $L^2$ representer exists if and only if \condC\ holds.
\end{enumerate}
\end{proposition}

\begin{proof}
(i) Along a trajectory approaching $x_*$ in backward time, \eqref{lin:eq:w} diverges linearly in the time variable if $\psi_*\ne0$. If $\psi_*=0$, it converges to $\Psi(\ell)$. Continuity at $x_*$ therefore forces $\psi_*=0$ and $\Psi$ constant. Conversely, under \condC\ the function $v_\psi$ takes a constant value $c$ on $\partial\cO$, namely $c=-\Psi$. Hence $v_\psi\circ\vartheta\equiv c$ and $w=v_\psi-c$ is smooth.

(ii) We have $A_{f_0}^*z=b\cdot\nabla(L_{f_0}^{-1}z)$, and $L_{f_0}^{-1}$ maps $L^2(\cO)$ isomorphically onto $H^2(\cO)\cap H^1_0(\cO)$. This gives the first equivalence, with $z=L_{f_0}\tilde w$ for a solution $\tilde w\in H^2\cap H^1_0$ of $b\cdot\nabla\tilde w=\psi$. Any such $\tilde w$ coincides with $w$ almost everywhere. To see this, use smooth flow-box coordinates $(t,\theta)$ on $\cO\setminus\{x_*\}$, in which $b\cdot\nabla=\partial_t$. Since $\tilde w\in H^1$, it is absolutely continuous on almost every coordinate line in these smooth coordinates (the ACL property). Thus for almost every $\theta$, the function $t\mapsto\tilde w(t,\theta)$ is absolutely continuous with derivative $\psi$, and its trace at the exit time is zero. So it agrees with \eqref{lin:eq:w} on almost every trajectory. Uniqueness follows.

(iii) If $d\le3$, then $H^2(\cO)\subset C(\overline\cO)$, so a representer forces $w$ to be continuous, and (i) applies. Conversely, (i) shows that \condC\ gives a smooth $w$.
\end{proof}

Continuity of $w$ is sufficient for \eqref{lin:eq:range} in every dimension. The Sobolev embedding argument above proves necessity for $d\le3$; Theorem~\ref{thm:linear} extends necessity to $d=4$. In dimensions $d\ge5$, an $L^2$ representer may exist even when the characteristic solution is discontinuous or unbounded.
\subsection{Two facts on linear Gaussian experiments}\label{lin:sec:gauss}

In the linearised experiment $Y=A_{f_0}h+N^{-1/2}\mathbb W$, the estimator $\widehat\ell_z=Y(z)=\ip{A_{f_0}h}z+N^{-1/2}\mathbb W(z)$, for $z\in L^2$, satisfies
\begin{equation}\label{lin:eq:linest}
 \E_h|\widehat\ell_z-\ell(h)|^2=\ip h{A_{f_0}^*z-\psi}^2+N^{-1}\norm z_2^2 .
\end{equation}
All lower bounds in this paper, in the linearised and in the nonlinear experiment, rest on the following elementary two-point bound.

\begin{lemma}[Two-point bound]\label{lin:lem:twopoint}
Let $P_0,P_1$ be probability measures with $\KL(P_0,P_1)\le1/2$, and let $\theta_0,\theta_1\in\R$. Then every estimator $\widehat\theta$ satisfies
\begin{equation}\label{lin:eq:twopoint}
 \max_{i=0,1}\E_{P_i}|\widehat\theta-\theta_i|\ge\tfrac18|\theta_1-\theta_0| .
\end{equation}
\end{lemma}

\begin{proof}
Let $\Delta=|\theta_1-\theta_0|$ and let $\phi$ be the indicator that $\widehat\theta$ is closer to $\theta_1$ than to $\theta_0$. On $\{\phi=1\}$ we have $|\widehat\theta-\theta_0|\ge\Delta/2$, and on $\{\phi=0\}$ we have $|\widehat\theta-\theta_1|\ge\Delta/2$. Hence $\max_i\E_{P_i}|\widehat\theta-\theta_i|\ge\frac\Delta4\{P_0(\phi=1)+P_1(\phi=0)\}\ge\frac\Delta4(1-\norm{P_0-P_1}_{TV})$, and Pinsker's inequality gives $\norm{P_0-P_1}_{TV}\le\sqrt{\KL/2}\le1/2$.
\end{proof}

In the linearised experiment, the Kullback--Leibler divergence between the laws of $Y$ under $h_0$ and $h_1$ is $\frac N2\norm{A_{f_0}(h_1-h_0)}_2^2$; in the nonlinear experiment \eqref{eq:model} it is $\frac N2\norm{u_{f_1}-u_{f_0}}_2^2$. If the divergence is bounded by a fixed constant rather than by $1/2$, the same argument gives \eqref{lin:eq:twopoint} with a smaller constant. Choosing $h_1=\epsilon h^\circ$ for a fixed smooth compactly supported $h^\circ$ with $\ell(h^\circ)\ne0$ and $\epsilon=N^{-1/2}$ gives the \emph{parametric lower bound}: the minimax risk is at least $cN^{-1/2}$ for every nonzero $\psi$. Upper bounds for the expected absolute error follow from \eqref{lin:eq:linest} by the Cauchy--Schwarz inequality.

We use Fourier analysis on $\R^d$ with $\widehat q(\xi)=\int e^{-ix\cdot\xi}q(x)\,dx$, and a fixed bounded extension operator $\mathcal E:H^s(\cO)\to H^s(\R^d)$. If $r\in H^{-s}(\R^d)$ is supported in a compact subset of $\cO$, then $|\ip hr|\le\norm{\mathcal Eh}_{H^s(\R^d)}\norm r_{H^{-s}(\R^d)}\le C\norm h_{H^s(\cO)}\norm r_{H^{-s}(\R^d)}$. Let $p\in C_c^\infty(\R^d)$ be radial, real, equal to one on $\{|\xi|\le1\}$ and supported in $\{|\xi|\le2\}$, and put $P_R=p(D/R)$.
\subsection{The quadratic model}\label{lin:sec:quadratic}

In this subsection $H=\operatorname{diag}(\lambda_1,\dots,\lambda_d)$, $\bar f>0$ and $c_\cO>0$ are constants, and
\[
 \cO=\{x:\ x^{\top}Hx<c_\cO\},\qquad u_0(x)=\tfrac12(x^{\top}Hx-c_\cO),\qquad f_0=\bar f,\qquad g=\bar fQ .
\]
This reference satisfies Assumption~\ref{ass:base}, with $x_*=0$ and $b(x)=Hx$. Write $\mathsf X=Hx\cdot\nabla$ and $\mathsf T_0=\diver(\cdot\,Hx)=\mathsf X+Q$. The linearised forward map is $A=-(\bar f\Delta_D)^{-1}\mathsf T_0$, where $\Delta_D$ is the Dirichlet Laplacian. The backward flow is $X_tx=e^{-tH}x$, and $J_t=e^{-Qt}$.

Let $\lambda_*=\lambda_d$ and $D_r=\operatorname{diag}(r^{a_1},\dots,r^{a_d})=e^{(\log r/\lambda_*)H}$ for $r>0$, so that $\det D_r=r^{\Teff}$. Since $\tau(e^{sH}x)=\tau(x)-s$ and $\vartheta$ is invariant along the flow, the singular part \eqref{lin:eq:F} satisfies
\begin{equation}\label{lin:eq:hom}
 F(D_rx)=F(x)+c_0\log r,\qquad c_0=\psi_*/\lambda_* .
\end{equation}
We use \eqref{lin:eq:hom} to extend $F$ to a smooth function on $\R^d\setminus\{0\}$. Write $x=(y,z)$, where $z$ collects the \emph{fast} coordinates $E=\{i:a_i=1\}$ and $y$ the remaining \emph{slow} ones, and let $\Delta_E=\sum_{i\in E}\partial_{ii}$. Note that $F\in L^2_{\rm loc}(\R^d)$, since $|F(x)|\le C(1+|\log|x||)$.

\begin{lemma}\label{lin:lem:fast}
Suppose $F$ is not constant. Then:
\begin{enumerate}
\item[(i)] for every $i\in E$, $\partial_{ii}F$ does not vanish identically on $\R^d\setminus\{0\}$;
\item[(ii)] $\Delta_EF$ is a nonzero distribution on $\R^d$; if $E\ne\{1,\dots,d\}$ or $d\ge3$, it is nonzero on $\R^d\setminus\{0\}$.
\end{enumerate}
\end{lemma}

\begin{proof}
(i) Fix $i\in E$, write $x=(x',x_i)$, and suppose $\partial_{ii}F=0$ off the origin. For every $x'\ne0$ the line $\{(x',x_i):x_i\in\R\}$ avoids the origin, so $F(x',x_i)=\alpha(x')x_i+\beta(x')$. Evaluating at $x_i=1$ and $x_i=2$ shows that $\alpha,\beta$ extend smoothly to $x'=0$, since $F$ is smooth near $(0,1)$ and $(0,2)$. Writing $D_r=\operatorname{diag}(D'_r,r)$, \eqref{lin:eq:hom} gives
\[
 \alpha(D'_rx')=r^{-1}\alpha(x'),\qquad \beta(D'_rx')=\beta(x')+c_0\log r .
\]
As $r\downarrow0$ the left-hand sides converge to $\alpha(0)$ and $\beta(0)$. Hence $\alpha\equiv0$, then $c_0=0$, and then $\beta\equiv\beta(0)$. So $F$ is constant on $\{x'\ne0\}$, and by continuity on $\R^d\setminus\{0\}$, a contradiction.

(ii) Suppose first that there is a slow coordinate, and that $\Delta_EF=0$ on $\R^d\setminus\{0\}$. Fix a slow vector $y\ne0$. Then $z\mapsto F(y,z)$ is an entire harmonic function of the fast variables. It grows at most logarithmically: for $|z|$ large choose $r$ with $r|z|=1$. Then $D_r(y,z)$ lies in a fixed compact set avoiding the origin, because its slow part tends to zero and its fast part has unit length. So \eqref{lin:eq:hom} gives $|F(y,z)|\le C+|c_0|\log|z|$. By the gradient estimate for harmonic functions, an entire harmonic function of logarithmic growth is constant. Hence $F(y,z)=G(y)$ for $y\ne0$. Since $G(y)=F(y,z_0)$ for a fixed $z_0\ne0$, $G$ extends smoothly to $y=0$. Then \eqref{lin:eq:hom} gives $G(D'_ry)=G(y)+c_0\log r$, and letting $r\downarrow0$ yields $c_0=0$ and $G\equiv G(0)$, a contradiction.

If all coordinates are fast, $H=\lambda_*I$ and $D_r=rI$. If $\Delta F=0$ as a distribution on $\R^d$, then $F$ is harmonic on $\R^d$ by Weyl's lemma and grows at most logarithmically at infinity, hence is constant. If $d\ge3$ and $\Delta F=0$ only on $\R^d\setminus\{0\}$, the singularity is removable, because $|F(x)|=O(|\log|x||)=o(|x|^{2-d})$; so again $F$ is constant. In $d=2$ with $H$ isotropic, $\Delta(\log|x|)=2\pi\delta_0$ shows that the distribution $\Delta F$ may be supported at the origin, which is why (ii) is stated on $\R^d$.
\end{proof}

\begin{proposition}[Representer criterion]\label{lin:prop:range-quad}
In the quadratic model, $A^*z=\psi$ has a solution $z\in L^2(\cO)$ if and only if \condC\ holds or $\Teff>4$.
\end{proposition}

\begin{proof}
By Proposition~\ref{prop:C-smooth}(ii) we must decide when $w=v_\psi+F\in H^2(\cO)$. Since $v_\psi$ is smooth, this is the question whether $F\in H^2$ near $0$. Use flow-polar coordinates $x=D_r\theta$, with $r\in(0,1]$ and $\theta\in\partial\cO$. In these coordinates
\begin{equation}\label{lin:eq:volume}
 dx=r^{\Teff-1}\mathsf j(\theta)\,dr\,d\sigma(\theta),
\end{equation}
where $\mathsf j=\lambda_*^{-1}\,H\theta\cdot\nu(\theta)$ is smooth and positive. Differentiating \eqref{lin:eq:hom},
\begin{equation}\label{lin:eq:hom-deriv}
 (\partial^\alpha F)(D_r\theta)=r^{-\sum_i\alpha_ia_i}(\partial^\alpha F)(\theta),\qquad |\alpha|\ge1 .
\end{equation}
\emph{Sufficiency.} If $\Teff>4$, then $\int_0^1r^{\Teff-1-2(a_i+a_j)}dr\le\int_0^1r^{\Teff-5}dr<\infty$, so all classical second derivatives are square integrable near $0$, and so are the first derivatives. To identify them with weak derivatives, let $\chi_\epsilon(x)=\chi(\rho(x)/\epsilon)$, where $\chi$ is a smooth cutoff and $\rho$ is the quasi-radius $\rho(D_r\theta)=r$. By \eqref{lin:eq:hom-deriv}, the $H^2$ norm of $\chi_\epsilon F$ is at most $C\epsilon^{\Teff/2-2}(1+|\log\epsilon|)\to0$. So the smooth functions $(1-\chi_\epsilon)F$ converge to $F$ in $H^2$ near $0$. If \condC\ holds, $F$ is constant by \eqref{lin:eq:C-F}.

\emph{Necessity.} Suppose \condC\ fails, so $F$ is not constant, and $\Teff\le4$. By Lemma~\ref{lin:lem:fast}(i) there is an open set $U\subset\partial\cO$ on which $\partial_{dd}F\ne0$. If $w\in H^2$, its weak second derivatives coincide on $\cO\setminus\{0\}$ with the classical ones. By \eqref{lin:eq:volume} and \eqref{lin:eq:hom-deriv}, $\int|\partial_{dd}F|^2$ over $\{D_r\theta:\theta\in U,\ 0<r<1\}$ is a positive multiple of $\int_0^1r^{\Teff-5}dr=\infty$. Since $\partial_{dd}v_\psi$ is bounded, this contradicts $w\in H^2$.
\end{proof}

The criterion cannot be circumvented by cancellations in the classical Laplacian off the origin. If the distributional Laplacian of $w$ were in $L^2$, elliptic regularity, the smooth boundary behaviour of $w$, and $w\in L^2$ would force $w\in H^2$. For example, in $d=2$ with $H$ isotropic, the term $\log|x|$ in $F$ is harmonic off the origin but contributes a point mass.

\begin{lemma}[Fourier shells]\label{lin:lem:band}
Let $\zeta\in C_c^\infty(\R^d)$ equal one near the origin and put $q=\zeta F$. Then
\[
 \int_{R\le|\xi|\le2R}|\widehat q(\xi)|^2\,d\xi\le CR^{-\Teff},\qquad R\ge1 .
\]
\end{lemma}

\begin{proof}
Let $\rho$ be the quasi-radius above; then $\rho(D_rx)=r\rho(x)$, $\rho(x)\asymp\max_i|x_i|^{1/a_i}$, and $|\{\rho<r\}|\asymp r^{\Teff}$. By \eqref{lin:eq:hom-deriv} and $a_i\le1$, for $0<\rho(x)\le1$ and every integer $k\ge1$,
\begin{equation}\label{lin:eq:Dk}
 |D^kF(x)|\le C_k\rho(x)^{-k}.
\end{equation}
Fix an integer $k>\Teff/2$ and a vector $h$ with $a:=|h|$ small, and let $\Delta_h^k$ denote the $k$-th forward difference.

\emph{Outer region $\rho(x)>Ca$.} Let $m=\rho(x)\le1$. Some coordinate satisfies $|x_i|\asymp m^{a_i}\ge m\ge Ca$. If $C\ge4k$, then every point $x+jh$, $0\le j\le k$, has quasi-radius comparable to $m$. By the mean value theorem and \eqref{lin:eq:Dk}, $|\Delta_h^kq(x)|\le Ca^km^{-k}$. Integrating with \eqref{lin:eq:volume},
\[
 \int_{\rho>Ca}|\Delta_h^kq|^2\le Ca^{2k}\int_{Ca}^1r^{\Teff-1-2k}dr+Ca^{2k}\le Ca^{\Teff},
\]
where the second term accounts for the region where $\zeta$ is not constant.

\emph{Inner region $\rho(x)\le Ca$.} Substitute $x=D_ay$. Then $x+jh=D_a(y+jh')$ with $h'=D_a^{-1}h$, and $|h'_i|=a^{-a_i}|h_i|\le a^{1-a_i}\le1$. For $a$ small, $\zeta=1$ at all points involved. By \eqref{lin:eq:hom}, $F(D_a(y+jh'))=F(y+jh')+c_0\log a$. The constant cancels in the difference because $k\ge1$. Hence
\[
 \int_{\rho\le Ca}|\Delta_h^kq|^2\,dx=a^{\Teff}\int_{\rho(y)\le C}|\Delta_{h'}^kF(y)|^2\,dy\le Ca^{\Teff}\int_{\{\rho\le C'\}}|F|^2\le Ca^{\Teff},
\]
since $F\in L^2_{\rm loc}$. Altogether $\norm{\Delta_h^kq}_2^2\le C|h|^{\Teff}$.

By Plancherel, $\norm{\Delta_h^kq}_2^2=(2\pi)^{-d}\int|e^{ih\cdot\xi}-1|^{2k}|\widehat q(\xi)|^2d\xi$. Average over $h$ uniformly distributed on the ball $\{|h|\le c_1/R\}$. For $R\le|\xi|\le2R$, the average of $|e^{ih\cdot\xi}-1|^{2k}$ is bounded below by a constant $c_2>0$ depending only on $c_1,k,d$; we take $c_1=1$. Therefore
\[
 c_2\int_{R\le|\xi|\le2R}|\widehat q|^2\le(2\pi)^d\sup_{|h|\le1/R}\norm{\Delta_h^kq}_2^2\le CR^{-\Teff}. \qedhere
\]
\end{proof}

\begin{proposition}[Rates in the quadratic model]\label{lin:prop:rates-quad}
The conclusion of Theorem~\ref{thm:linear} holds in the quadratic model. In mean squared error, the minimax risk is of order $N^{-1}$ if \condC\ holds or $\Teff>4$; of order $(\log N)/N$ if \condC\ fails and $\Teff=4$; and of order $N^{-(2s+\Teff)/(2s+4)}$ if \condC\ fails and $\Teff<4$.
\end{proposition}

\begin{proof}
\emph{Upper bound.} If \condC\ holds or $\Teff>4$, Proposition~\ref{lin:prop:range-quad} gives a representer, and \eqref{lin:eq:linest} gives the rate $N^{-1}$. Otherwise choose $\zeta\in C_c^\infty(\cO)$, equal to one near $0$, and put $q=\zeta F$ and $w_{\rm reg}=w-q$. Then $w_{\rm reg}\in C^\infty(\overline\cO)$ vanishes on $\partial\cO$. Also $\mathsf Xq=\psi_*\zeta+F\,\mathsf X\zeta$ is smooth and compactly supported, because $\mathsf X\zeta$ vanishes near $0$. Choose $\zeta_2\in C_c^\infty(\cO)$ equal to one on a neighbourhood of $\supp\zeta$, and define
\[
 w_R=w_{\rm reg}+\zeta_2P_Rq,\qquad z_R=\bar f\Delta w_R .
\]
Since $w_R\in H^2\cap H^1_0$, we have $A^*z_R=\mathsf Xw_R$. Since $\psi=\mathsf Xw_{\rm reg}+\mathsf Xq$ in the sense of distributions,
\[
 \psi-A^*z_R=(1-P_R)\mathsf Xq+[P_R,\mathsf X]q-\mathsf X\bigl((\zeta_2-1)P_Rq\bigr).
\]
The first term is $O(R^{-M})$ in $H^{-s}$ for every $M$, because $\mathsf Xq\in C_c^\infty$. The third is $O(R^{-M})$ in every Sobolev norm, because $1-\zeta_2$ vanishes on a neighbourhood of $\supp q$ and the kernel of $P_R$ decays rapidly on the scale $R^{-1}$. For the commutator, $\mathsf X=\sum_j\lambda_jx_j\partial_j$ and the Fourier transform conjugates $x_j$ to $i\partial_{\xi_j}$, so
\[
 \widehat{[P_R,\mathsf X]q}(\xi)=\frac{H\xi\cdot(\nabla p)(\xi/R)}{R}\,\widehat q(\xi).
\]
The multiplier is bounded by $2\lambda_*\norm{\nabla p}_\infty$ and supported in $\{R\le|\xi|\le2R\}$. By Lemma~\ref{lin:lem:band}, $\norm{[P_R,\mathsf X]q}_{H^{-s}}\le CR^{-s-\Teff/2}$. All terms are supported in a fixed compact subset of $\cO$ up to $O(R^{-M})$ tails. Hence the bias in \eqref{lin:eq:linest} is at most $C\eta R^{-s-\Teff/2}$.

For the variance, $\norm{z_R}_2\le\bar f\{\norm{\Delta w_{\rm reg}}_2+\norm{\zeta_2P_Rq}_{H^2}\}$. Summing Lemma~\ref{lin:lem:band} over dyadic shells $2^j\le2R$, and using $\Teff>1$ for the lower-order terms,
\begin{equation}\label{lin:eq:var}
 \norm{z_R}_2^2\le C\begin{cases}R^{4-\Teff},&\Teff<4,\\ 1+\log R,&\Teff=4.\end{cases}
\end{equation}
With $R=N^{1/(2s+4)}$, the squared bias $R^{-2s-\Teff}$ and the variance $N^{-1}R^{4-\Teff}$ are both of order $N^{-(2s+\Teff)/(2s+4)}$ when $\Teff<4$. When $\Teff=4$ the sum is $O((\log N)/N)$.

\emph{Lower bound, irregular regime ($\Teff<4$).} On $\R^d$ define
\[
 (\mathsf Sq)(x)=\int_0^\infty e^{-Qt}q(e^{-tH}x)\,dt .
\]
Since $\partial_t[e^{-Qt}q(e^{-tH}x)]=-\mathsf T_0[e^{-Qt}q(e^{-tH}\cdot)](x)$, we get $\mathsf T_0\mathsf Sq=q$. The Fourier transform of $e^{-Qt}q(e^{-tH}\cdot)$ is $\widehat q(e^{tH}\cdot)$. Substituting $\eta=e^{tH}\xi$ and using $|e^{-tH}\eta|\le|\eta|$ shows that its $H^s$ norm is at most $e^{-Qt/2}\norm q_{H^s}$. By Minkowski's inequality,
\begin{equation}\label{lin:eq:Ssob}
 \norm{\mathsf Sq}_{H^s(\R^d)}\le\frac2Q\norm q_{H^s(\R^d)},\qquad s\ge0 .
\end{equation}
This bound includes the whole transport tail of $\mathsf Sq$ along the flow, so no separate estimate of the tail is needed. On $\cO$, $\mathsf S$ coincides with the transport inverse \eqref{lin:eq:S}.

By Lemma~\ref{lin:lem:fast}(ii), there is $v\in C_c^\infty(\R^d)$ with
\[
 \Lambda:=\ip F{\Delta_Ev}\ne0 .
\]
For $r\in(0,1)$ small and $\beta>0$, put $v_r(x)=\beta\,v(D_r^{-1}x)$, so $\supp v_r\Subset\cO$, and
\[
 h_r=-\bar f\,\mathsf S\Delta v_r\big|_{\cO}.
\]
Then $\mathsf T_0h_r=-\bar f\Delta v_r$, so $Ah_r=\Delta_D^{-1}\Delta v_r=v_r$, since $v_r$ is compactly supported in $\cO$. Each derivative of $v(D_r^{-1}\cdot)$ costs at most a factor $r^{-1}$, because $a_i\le1$. With $\det D_r=r^{\Teff}$, interpolation and \eqref{lin:eq:Ssob} give
\begin{equation}\label{lin:eq:hr}
 \norm{h_r}_{H^s(\cO)}\le C\beta r^{\Teff/2-s-2},\qquad \norm{Ah_r}_2=\beta r^{\Teff/2}\norm v_2 .
\end{equation}
For the functional, use $\psi=\mathsf Xw$ and integrate by parts. The boundary term on $\partial\cO$ vanishes because $w=0$ there. The term on $\partial B_\epsilon(0)$ is $O(\epsilon^d|\log\epsilon|)$, since $h_r$ is smooth. Hence
\[
 \ell(h_r)=\ip{h_r}{\mathsf Xw}=-\ip{\mathsf T_0h_r}w=\bar f\ip{\Delta v_r}w .
\]
Substitute $x=D_ry$ and use $w(D_ry)=v_\psi(D_ry)+F(y)+c_0\log r$:
\[
 \ip{\Delta v_r}w=\beta r^{\Teff}\sum_i r^{-2a_i}\int(\partial_{ii}v)(y)\{v_\psi(D_ry)+F(y)+c_0\log r\}\,dy .
\]
The constant $c_0\log r$ integrates to zero against $\partial_{ii}v$. After integrating by parts twice, the term with $v_\psi$ equals $\beta r^{\Teff}\int v(y)\,(\Delta v_\psi)(D_ry)\,dy$, which is $O(\beta r^{\Teff})$. For $F$, the fast indices contribute $\beta r^{\Teff-2}\Lambda$, and each slow index contributes $\beta r^{\Teff-2}\cdot O(r^{2(1-a_i)})$. Therefore
\begin{equation}\label{lin:eq:pairing}
 \ell(h_r)=\bar f\beta r^{\Teff-2}\{\Lambda+o(1)\},\qquad r\downarrow0 .
\end{equation}
Take $\beta=\varsigma r^{s+2-\Teff/2}$, with $\varsigma>0$ small enough that $\norm{h_r}_{H^s}\le\eta$, and $r=N^{-1/(2s+4)}$. The divergence between $h=0$ and $h=h_r$ is $\frac N2\varsigma^2r^{2s+4}\norm v_2^2=\frac12\varsigma^2\norm v_2^2$, while $|\ell(h_r)|\ge c\varsigma r^{s+\Teff/2}=c\varsigma N^{-(2s+\Teff)/(4s+8)}$. Now \eqref{lin:eq:twopoint} gives the lower bound.

\emph{Lower bound, critical regime ($\Teff=4$).} We claim that $\Delta_EF\ne0$ on $\R^d\setminus\{0\}$. If $E$ is a proper subset, this is Lemma~\ref{lin:lem:fast}(ii). If all coordinates are fast, then $d=\Teff=4\ge3$, and Lemma~\ref{lin:lem:fast}(ii) applies again. So we may choose $v$ supported in a flow annulus $\{1<\rho<2\}$ with $\Lambda>0$. Let $B>2$ be such that the sets $D_{B^{-j}}\{1<\rho<2\}$ are disjoint, put $r_j=r_0B^{-j}$, and let
\[
 V_j=r_j^{-\Teff/2}v(D_{r_j}^{-1}\cdot)/\norm v_2,\qquad H_j=-\bar f\,\mathsf S\Delta V_j .
\]
By \eqref{lin:eq:hr} and \eqref{lin:eq:pairing} with $\Teff=4$, and for $r_0$ small enough,
\[
 \norm{V_j}_2=1,\qquad \ell(H_j)\ge c_*>0,\qquad \norm{H_j}_{H^s}\le Cr_j^{-s-2} .
\]
Let $\epsilon=N^{-1/2}$ and $m=\lfloor\log(1/\epsilon)/(2(s+2)\log B)\rfloor$, and put
\[
 V=\frac{\varsigma\epsilon}{\sqrt m}\sum_{j=1}^mV_j,\qquad h=-\bar f\,\mathsf S\Delta V .
\]
The $V_j$ have disjoint supports, so $\norm{Ah}_2=\norm V_2=\varsigma\epsilon$. By the triangle inequality and geometric summation,
\[
 \norm h_{H^s}\le C\frac{\varsigma\epsilon}{\sqrt m}r_m^{-s-2}\le C\frac{\varsigma\epsilon^{1/2}}{\sqrt m}\to0,
\]
while $\ell(h)\ge c_*\varsigma\epsilon\sqrt m\ge c\varsigma\sqrt{(\log N)/N}$. Then \eqref{lin:eq:twopoint} gives the lower bound $\sqrt{(\log N)/N}$.

\emph{Regular regime.} The parametric lower bound of Section~\ref{lin:sec:gauss} completes the proof.
\end{proof}

\begin{remark}\label{lin:rem:shape}
The alternative $h_r$ is built from a localised change $v_r$ of the observation, of width $r$ in the fast directions and $r^{a_i}$ in the others, followed by the transport inverse. The conductivity perturbation therefore consists of a core of the same shape plus a transport tail along the flow lines leaving the core. Simple conductivity bumps, also anisotropic ones, give only the exponent with $d$ in place of $\Teff$, because the tail is what makes $Ah_r$ small.
\end{remark}
\subsection{General references and resonant normal forms}\label{sec:normalform}

We now remove the quadratic assumption. The upper bound needs no normal form: it only uses flow estimates near $x_*$. The lower bound and the necessity part use a smooth conjugacy to a polynomial normal form. In it, the fast variables enter linearly, so resonances do not change the scaling.

Translate $x_*$ to the origin and rotate so that $H=Db(0)=\operatorname{diag}(\lambda_i)$.

\begin{lemma}[Fourier shells for a general flow]\label{lin:lem:general-band}
Let $F$ be smooth on a punctured neighbourhood of $0$, with $b\cdot\nabla F=c$ there and smooth values on a small sphere around $0$. Let $\zeta$ be a smooth cutoff supported in that neighbourhood and equal to one near $0$, and put $q=\zeta F$. Then $\int_{R\le|\xi|\le2R}|\widehat q|^2\le CR^{-\Teff}$ for $R\ge1$.
\end{lemma}

\begin{proof}
Choose a ball $B$ around $0$ on which $Db\ge mI$ for some $m>0$. Then $X_t(B)\subset B$, $|X_tx|\le e^{-mt}|x|$, and since $Db(x)=H+O(|x|)$ and $\diver b=Q+O(|x|)$, integration along the flow gives
\begin{equation}\label{lin:eq:capJ}
 C^{-1}e^{-Qt}\le J_t(x)\le Ce^{-Qt},\qquad x\in B,\ t\ge0 .
\end{equation}
For $x,y\in\overline B$ we have $\frac d{dt}|X_tx-X_ty|^2\ge-2\{\lambda_*+C(|X_tx|+|X_ty|)\}|X_tx-X_ty|^2$. The error term has bounded time integral, so
\begin{equation}\label{lin:eq:colip}
 |X_tx-X_ty|\ge C^{-1}e^{-\lambda_*t}|x-y| .
\end{equation}
Let $E_t=X_t(B)$, and fix $C_0$ with $X_{C_0}(\overline B)\subset\frac12B$. Applying \eqref{lin:eq:colip} to points of $X_{C_0}(\overline B)$ and $\partial B$, and using that $\partial E_t=X_t(\partial B)$, gives
\begin{equation}\label{lin:eq:gap}
 \operatorname{dist}(\overline E_{t+C_0},\partial E_t)\ge c_1e^{-\lambda_*t} .
\end{equation}
The finite differences $\Delta_h^kq$ used below compare $F$ at points up to $k|h|$ apart. Since $F$ is controlled through the exit time, we need such points to lie at comparable exit times. This is what \eqref{lin:eq:gap} provides: at the scale $|h|\asymp e^{-\lambda_*t}$, a point of $E_t$ shifted by $k|h|$ stays in $E_{t-C_0}$, so the shift costs at most a bounded amount of flow time.

Let $\tau_B$ be the exit time from $B$. Every $x\in B\setminus\{0\}$ is uniquely $x=X_r\theta$ with $\theta\in\partial B$ and $r=\tau_B(x)$. The volume element is $J_r(\theta)\,b(\theta)\cdot\nu(\theta)\,dr\,d\sigma(\theta)\asymp e^{-Qr}dr\,d\sigma$. Along each orbit segment inside $B$, the time integral of $\norm{Db-H}$ is bounded. Hence the forward variational propagator satisfies $\norm{D\phi_t}\le Ce^{\lambda_*t}$, and, by induction in the higher variational equations, $\norm{D^k\phi_t}\le C_ke^{k\lambda_*t}$. Implicit differentiation of the exit condition, which is transversal, gives $|D^k\tau_B|+|D^k\vartheta_B|\le C_ke^{k\lambda_*\tau_B}$ for the exit point $\vartheta_B$. Since $F(x)=F(\vartheta_B(x))-c\,\tau_B(x)$,
\begin{equation}\label{lin:eq:genDk}
 |D^kF(x)|\le C_ke^{k\lambda_*\tau_B(x)},\qquad |F(x)|\le C(1+\tau_B(x)) .
\end{equation}
Fix an integer $k>\Teff/2$. For a small vector $h$, choose $t$ with $|h|=c_2e^{-\lambda_*t}$, where $c_2\le c_1/(2k)$. By \eqref{lin:eq:gap}, all sets $E_t+jh$, $0\le j\le k$, lie in $E_{t-C_0}$. Increasing the threshold on $t$, they also lie where $\zeta=1$. The constant $-ct$ cancels in the difference. Using $y=X_{t-C_0}z$, \eqref{lin:eq:capJ} and $F(X_rz)=F(z)-cr$,
\begin{align*}
 \int_{E_t}|\Delta_h^kq|^2&\le C\sum_{j=0}^k\int_{E_t}|F(x+jh)+ct|^2dx\le C\int_{E_{t-C_0}}|F+ct|^2\\
 &\le Ce^{-Qt}\int_B|F+cC_0|^2\le Ce^{-Qt}.
\end{align*}
If $x\in B\setminus E_t$ has $r=\tau_B(x)\le t$, then by \eqref{lin:eq:gap} every point used in the difference has exit time at most $r+C_0$. So \eqref{lin:eq:genDk} bounds $|\Delta_h^kq(x)|$ by $C|h|^ke^{k\lambda_*r}$, and
\[
 \int_{B\setminus E_t}|\Delta_h^kq|^2\le C|h|^{2k}\int_0^te^{(2k\lambda_*-Q)r}dr\le Ce^{-Qt}\le C|h|^{\Teff}.
\]
The part outside $B$ is smooth and contributes $O(|h|^{2k})$. The Plancherel averaging argument of Lemma~\ref{lin:lem:band} concludes the proof.
\end{proof}

In the upper bound, the truncated representer leaves the residual $\psi-A_{f_0}^*z_R$. Of its terms, only a commutator of $P_R$ with the transport field sees the singularity of $F$. The commutator acts on frequencies $|\xi|\asymp R$ only, so the shell bound converts into the bias $R^{-s-\Teff/2}$. For the quadratic model this was an explicit Fourier multiplier; for a variable field we use the following lemma.

\begin{lemma}[Commutator with a variable field]\label{lin:lem:comm}
Let $\tilde b\in C_c^\infty(\R^d;\R^d)$ and let $q\in L^2(\R^d)$ satisfy $\int_{R\le|\xi|\le2R}|\widehat q|^2\le CR^{-2\beta}$ for all $R\ge1$, for some $\beta>0$. Then for every $s\ge0$,
\[
 \norm{[P_R,\tilde b\cdot\nabla]q}_{H^{-s}}\le CR^{-s-\beta}.
\]
\end{lemma}

\begin{proof}
With $\omega=\xi-\eta$, the Fourier kernel of the commutator is $i\,\widehat{\tilde b}(\omega)\cdot\eta\,\{p(\xi/R)-p(\eta/R)\}$, up to a factor $(2\pi)^{-d}$. For all $\xi,\eta$,
\begin{equation}\label{lin:eq:symb}
 |\eta|\,|p(\xi/R)-p(\eta/R)|\le C|\xi-\eta| .
\end{equation}
If $|\eta|\le4R$, this follows from the Lipschitz bound on $p$. If $|\eta|>4R$, a nonzero difference requires $|\xi|\le2R$, hence $|\xi-\eta|\ge|\eta|/2$. Split $\tilde b=\tilde b_{\rm lo}+\tilde b_{\rm hi}$ with a smooth Fourier cutoff at $|\omega|\le R/16$. Since $\norm{|\omega|\widehat{\tilde b}_{\rm hi}}_{L^1}=O(R^{-M})$ for every $M$, \eqref{lin:eq:symb} and Young's inequality bound the high part by $O(R^{-M})\norm q_2$ in $L^2$. For the low part the kernel vanishes unless $R/2\le|\eta|\le3R$, and then $|\xi|$ also lies between fixed multiples of $R$. Young's inequality, \eqref{lin:eq:symb} and the shell bound give an $L^2$ bound $CR^{-\beta}$, with output frequencies of order $R$. This gives a further factor $R^{-s}$ in $H^{-s}$.
\end{proof}

\begin{proposition}[Upper bound for a general reference]\label{lin:prop:general-upper}
Under Assumption~\ref{ass:base}, for every $s\ge0$ the linearised minimax mean squared risk is at most $CN^{-1}$ if \condC\ holds or $\Teff>4$; $C(\log N)/N$ if $\Teff=4$; and $CN^{-(2s+\Teff)/(2s+4)}$ if $\Teff<4$. If $\Teff>4$, then $w\in H^2(\cO)$, so an $L^2$ representer exists.
\end{proposition}

\begin{proof}
Use the decomposition $w=v_\psi+F$ from \eqref{lin:eq:F}, where $b\cdot\nabla F=\psi_*$ in the sense of distributions. Choose $\zeta$ as in Lemma~\ref{lin:lem:general-band}, and put $q=\zeta F$ and $w_{\rm reg}=w-q$. Then $b\cdot\nabla q=\psi_*\zeta+F\,b\cdot\nabla\zeta$ is smooth and compactly supported, and $w_{\rm reg}$ is smooth with zero boundary values. Define $w_R=w_{\rm reg}+\zeta_2P_Rq$ and $z_R=L_{f_0}w_R$. Then $A_{f_0}^*z_R=b\cdot\nabla w_R$ exactly, and $\norm{z_R}_2\le C\norm{w_R}_{H^2}$. The residual $\psi-A_{f_0}^*z_R$ decomposes as in the quadratic case. The commutator term is estimated by Lemma~\ref{lin:lem:comm} with $\tilde b$ a compactly supported extension of $b$ and $\beta=\Teff/2$ from Lemma~\ref{lin:lem:general-band}. This gives bias $O(R^{-s-\Teff/2})$, and the variance bounds \eqref{lin:eq:var} follow by summing shells. The choice $R=N^{1/(2s+4)}$ gives the rates. If $\Teff>4$, the shell bound gives $\int|\xi|^4|\widehat q|^2<\infty$, so $w\in H^2$ and Proposition~\ref{prop:C-smooth}(ii) gives the representer. If \condC\ holds, $w$ is smooth.
\end{proof}

For the lower bound we need a normal form. A monomial vector field $y^\alpha e_i$ with $|\alpha|\ge2$ is \emph{resonant} if $\alpha\cdot\lambda=\lambda_i$. Since all $\lambda_j>0$, resonance forces $\lambda_j<\lambda_i$ whenever $\alpha_j\ge1$, and $|\alpha|\le\lambda_*/\lambda_1$. So there are finitely many resonant monomials, and the $i$-th component of a resonant field depends only on coordinates with eigenvalues strictly smaller than $\lambda_i$. In particular, no resonant field depends on a fast coordinate.

\begin{lemma}[Positive-spectrum normal form]\label{lin:lem:normalform}
There is a local $C^\infty$ diffeomorphism $\mathcal K$ near $0$, with $\mathcal K(0)=0$ and $D\mathcal K(0)=I$, which conjugates $b$ to
\[
 \mathcal N(y)=Hy+P(y),
\]
where $P$ is a finite sum of resonant monomials. In particular $\diver P=0$ and $[Hy,P]=0$.
\end{lemma}

This is the Poincar\'e-domain case of the smooth normal form theorem of Sternberg and Chen \cite{sternberg1957local,chen1963equivalence}. A self-contained proof is given in Section~\ref{lin:app}.

Let $\phi^{\mathcal N}_t$ be the flow of $\mathcal N$. Since $[Hy,P]=0$, the flows of $Hy$ and $P$ commute, so
\[
 \phi^{\mathcal N}_t(x)=e^{tH}G_t(x),\qquad G_t=\psi^P_t,
\]
where $\psi^P_t$ is the flow of $P$. Because $P$ is triangular and each of its components depends only on strictly slower coordinates, $\psi^P_t$ is obtained by solving successively in groups of increasing eigenvalue, and $G_t(x)$ is a polynomial in $t$ and $x$. The flow $\phi^{\mathcal N}_t$ itself is not polynomial in $t$. Since the fast components of $\mathcal N$ are $\lambda_*z+P_E(y)$ and the slow components do not involve $z$,
\begin{equation}\label{lin:eq:fast-linear}
 \phi^{\mathcal N}_t(y,z+a)=\phi^{\mathcal N}_t(y,z)+(0,e^{\lambda_*t}a),
\end{equation}
where $z$ are the fast coordinates. It is equivariant, $\psi^P_t(e^{sH}x)=e^{sH}\psi^P_t(x)$, and its $k$-th component is a polynomial in which every monomial $x^\alpha$ satisfies $\sum_i\alpha_i\lambda_i=\lambda_k$. Differentiating the equivariance relation gives derivative bounds at the \emph{shrinking} initial points $x=\phi^{\mathcal N}_{-t}\theta$, which are the only points where they are used. Let $\Theta$ be a fixed compact set and $\theta\in\Theta$. Then, at $x=\phi^{\mathcal N}_{-t}\theta$,
\begin{equation}\label{lin:eq:orbit-derivs}
 \partial_{x_{i_1}}\cdots\partial_{x_{i_m}}\bigl(\phi^{\mathcal N}_t\bigr)_k(x)
 =e^{(\lambda_{i_1}+\dots+\lambda_{i_m})t}\,\pi_{k,i}(t,\theta),
\end{equation}
where $\pi_{k,i}$ is a polynomial, and $\pi_{k,i}=0$ unless $\lambda_{i_1}+\dots+\lambda_{i_m}\le\lambda_k\le\lambda_*$. In particular, first derivatives in a slow coordinate $j$ are $e^{\lambda_jt}$ times a polynomial in $(t,\theta)$, and higher derivatives are $O((1+t)^Ce^{\lambda_*t})$. These bounds hold only along backward orbits of a compact set, not globally. For $\dot y=y$, $\dot z=2z+y^2$ one has $\partial_yz(t)=2tye^{2t}$, which is of the form $e^{\lambda_yt}\,\mathrm{poly}(t)$ only when $y=O(e^{-t})$. Below they are used only on $\phi^{\mathcal N}_{-t}(\supp\mathsf q)$, the support of $\mathsf q\circ\phi^{\mathcal N}_t$.

\begin{lemma}\label{lin:lem:normalfast}
Let $F^{\mathcal N}$ be smooth on a punctured neighbourhood of $0$, with $\mathcal NF^{\mathcal N}=c$ there, and extend it to $\R^d\setminus\{0\}$ by $F^{\mathcal N}(y)=F^{\mathcal N}(\phi^{\mathcal N}_{-t}y)+ct$ for large $t$. Then for fast indices $i,j\in E$,
\[
 (\partial_{ij}F^{\mathcal N})(\phi^{\mathcal N}_{-t}y)=e^{2\lambda_*t}(\partial_{ij}F^{\mathcal N})(y),\qquad \det D\phi^{\mathcal N}_{-t}=e^{-Qt}.
\]
If $F^{\mathcal N}$ is not constant, then some fast pure second derivative of $F^{\mathcal N}$ is not identically zero, and $\Delta_EF^{\mathcal N}$ is a nonzero distribution on $\R^d$. If moreover $E$ is a proper subset or $d\ge3$, then $\Delta_EF^{\mathcal N}\ne0$ on $\R^d\setminus\{0\}$.
\end{lemma}

\begin{proof}
Backward trajectories of $\mathcal N$ converge to $0$, so the extension is well defined and smooth off $0$. The scaling identity follows from differentiating $F^{\mathcal N}(\phi^{\mathcal N}_{-t}x)=F^{\mathcal N}(x)-ct$ twice in fast directions and using \eqref{lin:eq:fast-linear}. The Jacobian identity follows from $\diver\mathcal N=Q$.

Suppose all fast second derivatives vanish. Then $F^{\mathcal N}(y,z)=A(y)\cdot z+B(y)$ for $y\ne0$, and $A,B$ extend smoothly to $y=0$, as in Lemma~\ref{lin:lem:fast}. The terms linear in $z$ in the equation $\mathcal NF^{\mathcal N}=c$ give $\mathcal N_{\rm slow}A+\lambda_*A=0$. Along the backward slow flow, $A$ therefore grows like $e^{\lambda_*t}$ unless $A=0$, which contradicts smoothness at $0$. So $A=0$. Then $\mathcal N_{\rm slow}B=c$, smoothness at $0$ forces $c=0$, and backward convergence makes $B$ constant. Separating one fast coordinate at a time shows in the same way that no fast pure second derivative vanishes identically.

For fixed slow $y\ne0$, the function $z\mapsto F^{\mathcal N}(y,z)$ grows at most logarithmically. Take $t=\lambda_*^{-1}\log|z|+C$ for large $|z|$. Then $\phi^{\mathcal N}_{-t}(y,z)$ lies in a fixed compact annulus near the fast unit sphere: its slow part tends to zero, and its fast part is $e^{-\lambda_*t}z$ plus polynomial terms in $t$ times $e^{-\lambda_*t}$. If $\Delta_EF^{\mathcal N}=0$ off $0$, the harmonic argument of Lemma~\ref{lin:lem:fast}(ii) makes $F^{\mathcal N}$ independent of $z$, and the transport equation then makes it constant. If there are no slow coordinates, $P=0$, and Lemma~\ref{lin:lem:fast} applies directly.
\end{proof}

We next record the forward-perturbation construction. It is also used in the nonlinear lower bounds of Section~\ref{sec:lower-proofs}.

\begin{proposition}[Flow-adapted forward perturbations]\label{lin:prop:construction}
Let Assumption~\ref{ass:base} hold, suppose \condC\ fails, and let $F^{\mathcal N}(y)=F(\mathcal K^{-1}(y))$ be the singular part \eqref{lin:eq:F} in normal-form coordinates. Choose $\mathsf q\in C_c^\infty(\R^d)$ with $\Lambda=\ip{F^{\mathcal N}}{\Delta_E\mathsf q}\ne0$, supported in a fixed flow annulus if $\Teff=4$. For $t\ge t_0$ and $\beta>0$ let
\[
 V_t(x)=\beta\,\mathsf q\bigl(\phi^{\mathcal N}_t(\mathcal K(x))\bigr),\qquad h_t=-S_{u_0}L_{f_0}V_t,
\]
with $V_t$ extended by zero outside the chart. Then $V_t\in C_c^\infty(\cO)$, $A_{f_0}h_t=V_t$, and for every $s\ge0$,
\begin{equation}\label{lin:eq:gen-norms}
\begin{gathered}
 \norm{V_t}_2\asymp\beta e^{-Qt/2},\qquad \norm{h_t}_{H^s(\cO)}\le C\beta e^{(s+2)\lambda_*t-Qt/2},\\
 \ell(h_t)=\beta e^{(2\lambda_*-Q)t}\{f_0(0)\Lambda+o(1)\}
\end{gathered}
\end{equation}
as $t\to\infty$, with constants independent of $t$ and $\beta$.
\end{proposition}

\begin{proof}
For $t_0$ large the support $\mathcal K^{-1}(\phi^{\mathcal N}_{-t}\supp\mathsf q)$ lies in the chart and shrinks to $0$. By Lemma~\ref{lin:lem:S}, $T_{u_0}h_t=-L_{f_0}V_t$, so $A_{f_0}h_t=L_{f_0}^{-1}L_{f_0}V_t=V_t$. The $L^2$ norm follows from $\det D\phi^{\mathcal N}_{-t}=e^{-Qt}$ and the fact that the density of $\mathcal K$ tends to one on the shrinking support.

For the $H^s$ bound, consider integer $s$ first. By the chain rule, $D^k(\mathsf q\circ\phi^{\mathcal N}_t)$ is a sum of products of derivatives of $\mathsf q$ with $m\le k$ factors $D^{j_1}\phi^{\mathcal N}_t,\dots,D^{j_m}\phi^{\mathcal N}_t$, where $j_1+\dots+j_m=k$. On the support $\phi^{\mathcal N}_{-t}(\supp\mathsf q)$, \eqref{lin:eq:orbit-derivs} applies. First derivatives there are at most $Ce^{\lambda_*t}$, with no polynomial loss: fast derivatives are exactly $e^{\lambda_*t}$ by \eqref{lin:eq:fast-linear}, and slow ones are $e^{\lambda_jt}\,\mathrm{poly}(t)\le Ce^{\lambda_*t}$. Higher derivatives are at most $C\,\mathrm{poly}(t)e^{\lambda_*t}$. A product with $m<k$ factors is therefore $O(e^{k\lambda_*t})$. Hence $\norm{V_t}_{H^{k}}\le C\beta e^{k\lambda_*t-Qt/2}$, and $\norm{L_{f_0}V_t}_{H^s}\le C\beta e^{(s+2)\lambda_*t-Qt/2}$ by interpolation. Lemma~\ref{lin:lem:S} gives the bound on $h_t$.

For the pairing we proceed in three steps.

\emph{Step 1: reduction to the singular part.} The argument of the quadratic case gives $\ell(h_t)=\ip{h_t}{b\cdot\nabla w}=\ip{L_{f_0}V_t}w$. Since $V_t$ is compactly supported, $\int L_{f_0}V_t=0$, so we may replace $w$ by $w+\psi_*t-v_\psi(0)$. In normal coordinates and at the shrinking points $y=\phi^{\mathcal N}_{-t}\theta$, $\theta\in\supp\mathsf q$, this function equals
\[
 F^{\mathcal N}(\phi^{\mathcal N}_{-t}\theta)+\psi_*t+v_\psi(\mathcal K^{-1}\phi^{\mathcal N}_{-t}\theta)-v_\psi(0)=F^{\mathcal N}(\theta)+o(1)
\]
uniformly in $\theta$.

\emph{Step 2: change of variables.} Pull back by $\phi^{\mathcal N}_{-t}$, that is, substitute $y=\phi^{\mathcal N}_{-t}\theta$. The Jacobian is $e^{-Qt}$ times the density of $\mathcal K$, which tends to one on the shrinking support.

\emph{Step 3: the surviving principal term.} In normal coordinates, $L_{f_0}$ has principal matrix $f_0(0)I$ at the origin, because $D\mathcal K(0)=I$, and smooth lower-order coefficients. Divide by $e^{2\lambda_*t}$. By \eqref{lin:eq:fast-linear}, fast second derivatives of $\mathsf q\circ\phi^{\mathcal N}_t$ at $\phi^{\mathcal N}_{-t}\theta$ equal $e^{2\lambda_*t}(\partial_{ij}\mathsf q)(\theta)$, so the fast part of $L_{f_0}(\mathsf q\circ\phi^{\mathcal N}_t)$ converges uniformly to $f_0(0)\Delta_E\mathsf q(\theta)$. By \eqref{lin:eq:orbit-derivs}, mixed and slow second derivatives, first-order terms and the variation of the coefficients over the shrinking support carry an extra factor that is a polynomial in $t$ times a decaying exponential; these terms are negligible.

Combining the three steps with $F^{\mathcal N}\in L^2_{\rm loc}$ proves the expansion of $\ell(h_t)$.
\end{proof}

\begin{proof}[Proof of Theorem~\ref{thm:linear}]
The upper bounds and the sufficiency part of the representer criterion are Proposition~\ref{lin:prop:general-upper}. For the parametric case, if \condC\ holds or $\Teff>4$, the representer gives the rate $N^{-1/2}$, and the parametric lower bound of Section~\ref{lin:sec:gauss} matches.

\emph{Necessity of the criterion.} Suppose \condC\ fails, so $F$ is not constant by \eqref{lin:eq:C-F}, and $\Teff\le4$. By Lemma~\ref{lin:lem:normalfast} there is an open set of the small transverse sphere, in normal coordinates, on which a fast pure second derivative of $F^{\mathcal N}$ does not vanish. Along the inward flow its squared integral is proportional to $\int_0^\infty e^{(4\lambda_*-Q)t}\,dt=\infty$. Since smooth changes of coordinates preserve local $H^2$ membership and $v_\psi$ is smooth, $w\notin H^2$, and Proposition~\ref{prop:C-smooth}(ii) excludes an $L^2$ representer.

\emph{Lower bounds, irregular regime ($\Teff<4$).} Apply Proposition~\ref{lin:prop:construction} with $\beta=\varsigma_0e^{-(s+2)\lambda_*t+Qt/2}$ and $e^{(s+2)\lambda_*t}=N^{1/2}$. Then $\norm{h_t}_{H^s}\le C\varsigma_0\le\eta$ for $\varsigma_0$ small, and the divergence is $\frac N2\norm{V_t}_2^2\le C\varsigma_0^2$. The separation is
\[
 |\ell(h_t)|\ge c\varsigma_0e^{-(s+2)\lambda_*t+Qt/2+(2\lambda_*-Q)t}=c\varsigma_0e^{-(s+\Teff/2)\lambda_*t}=c\varsigma_0N^{-(2s+\Teff)/(4s+8)} .
\]
Then \eqref{lin:eq:twopoint} applies.

\emph{Lower bound, critical regime ($\Teff=4$).} By Lemma~\ref{lin:lem:normalfast}, $\Delta_EF^{\mathcal N}\ne0$ off $0$. Here, if $E$ is all of $\{1,\dots,d\}$, then $d=4$; in that case $H$ is isotropic, there are no resonances, and $P=0$. Choose $\mathsf q$ in a fixed flow annulus with $\Lambda>0$, and times $t_j=t_0+j\Delta$ with $\Delta$ so large that the supports of $V_{t_j}$ are disjoint. Normalise $V_{t_j}$ in $L^2$. By \eqref{lin:eq:gen-norms} with $Q=4\lambda_*$, the normalised functions have pairings bounded below by a positive constant and $H^s$ norms of the corresponding $h$ at most $Ce^{(s+2)\lambda_*t_j}$. The superposition argument of Proposition~\ref{lin:prop:rates-quad} then applies verbatim with $m\asymp\log N$ terms and gives the lower bound $\sqrt{(\log N)/N}$.

Combining the upper and lower bounds gives the stated rates, including the case where \condC\ fails and $\Teff>4$, in which $\rate(s,\Teff)=N^{-1/2}$.
\end{proof}

\begin{proof}[Proof of Corollary~\ref{cor:linear}]
In the linearised experiment, the information for $\ell$ is positive exactly when $\psi$ lies in the range of $A_{f_0}^*$ on $L^2$, that is, when an $L^2$ representer exists; see, for example, \cite[Chapter~25]{vandervaart1998asymptotic}. By Theorem~\ref{thm:linear}, this happens if and only if \condC\ holds or $\Teff_0>4$, and these are exactly the cases with rate $N^{-1/2}$. Since $\Teff_0\le d$, the second alternative is impossible for $d\le4$. For $d\ge5$, take the quadratic model with $H=\operatorname{diag}(1,\dots,1,\lambda)$ and $\lambda\ge(d-1)/3$. Then $\Teff=(d-1+\lambda)/\lambda\le4$, and any $\psi$ with $\psi(0)\ne0$ violates \condC.
\end{proof}

\section{Proofs of the nonlinear lower bounds}\label{sec:lower-proofs}

This section proves Theorem~\ref{thm:lower}, the lower bound in Theorem~\ref{thm:isotropic}, Proposition~\ref{prop:instability}, and the statement on interior perturbation classes in Remark~\ref{rem:collar}. All lower bounds come from two-point (or multi-scale two-point) arguments in the nonlinear experiment \eqref{eq:model}, through Lemma~\ref{lin:lem:twopoint}. Each requires three estimates for an alternative $f_1=f_0+h$: membership $\norm h_{H^s}\le\eta$, a divergence bound $N\norm{u_{f_1}-u_{f_0}}_2^2\lesssim1$, and a separation $|\ip h\psi|$ of the order of the claimed rate. In every construction the forward difference $u_{f_1}-u_{f_0}$ is controlled through the exact identity \eqref{eq:exact}, and no Taylor expansion of the forward map is used for the divergence. The constructions are of two kinds.
\begin{enumerate}
\item[(i)] \emph{Conductivity bumps} at the critical point (part (a)). Here the forcing generated by the linear part of $\nabla u_0$ has vanishing mass and first moments, while the $O(|x-x_*|^2)$ remainder is controlled separately in the dual $H^2$ norm. The duality estimate of Lemma~\ref{low:duality} and the Bramble--Hilbert lemma then give $\norm{u_{f_1}-u_{f_0}}_2\lesssim\norm h_\infty\rho^{d/2+2}$ for a bump of radius $\rho$, while $\ip h\psi\asymp\psi(x_*)\norm h_\infty\rho^d$.
\item[(ii)] \emph{Forward-first perturbations}, in which the change $V$ of the potential is prescribed and the conductivity is recovered by the transport inverse (part (b), the isotropic class and the instability results). Then $u_{f_1}-u_{f_0}$ is exactly $V$, or exactly comparable to it by Lemma~\ref{low:transfer}, and only the $H^s$ norm of the recovered conductivity and its pairing with $\psi$ have to be estimated.
\end{enumerate}

\subsection{Tools}\label{low:sec:tools}

Write $\mathcal V=H^2(\cO)\cap H^1_0(\cO)$ with the $H^2$ norm and $\mathcal V^*$ for its dual. For a potential $u$ write $T_uh=\diver(h\nabla u)$.

\begin{lemma}[Uniform duality estimate]\label{low:duality}
Let $f_0$ be as in Assumption~\ref{ass:base} and $s>d/2$. There are $\eta_0,C>0$ such that every $f\in\cF_s(f_0,\eta_0)$ satisfies $\inf_\cO f>0$, $L_f:\mathcal V\to L^2(\cO)$ is an isomorphism with $\norm{L_f^{-1}}_{L^2\to\mathcal V}\le C$, and
\begin{equation}\label{low:eq:duality}
 \norm{L_f^{-1}F}_{2}\le C\norm F_{\mathcal V^*},\qquad F\in\mathcal V^*.
\end{equation}
\end{lemma}

\begin{proof}
Choose $p>d$ with $1+d/2-d/p<s$. Then $H^s\hookrightarrow L^\infty\cap W^{1,p}$, which gives positivity for small $\eta_0$. With $q=2p/(p-2)<2d/(d-2)$ (any $q<\infty$ if $d=2$) we have $H^2\hookrightarrow W^{1,q}$, so for $v\in\mathcal V$ and $h=f-f_0$,
\[
 \norm{\diver(h\nabla v)}_2\le\norm h_\infty\norm{\Delta v}_2+\norm{\nabla h}_p\norm{\nabla v}_q\le C\norm h_{H^s}\norm v_{H^2}.
\]
Thus $L_fv=L_{f_0}v+\diver(h\nabla v)$ is a perturbation of the isomorphism $L_{f_0}:\mathcal V\to L^2$ of norm at most $C\eta_0$, and a Neumann series gives the uniform inverse. Since $L_f$ is symmetric on $\mathcal V$, for $a\in L^2$ we get $\ip{L_f^{-1}F}a=F(L_f^{-1}a)\le\norm F_{\mathcal V^*}\norm{L_f^{-1}a}_{H^2}\le C\norm F_{\mathcal V^*}\norm a_2$.
\end{proof}

For $f_0,f_1$ in this class, the Kullback--Leibler divergence between the laws $P_{f_0},P_{f_1}$ of $Y$ in \eqref{eq:model} is
\begin{equation}\label{low:eq:KL}
 \KL(P_{f_0},P_{f_1})=\frac N2\norm{u_{f_1}-u_{f_0}}_2^2 .
\end{equation}

If this divergence is at most $1/2$, Lemma~\ref{lin:lem:twopoint} with $\theta_i=\Phi(f_i)$ gives $\max_{i=0,1}\E_{f_i}|\widehat\Phi-\Phi(f_i)|\ge\frac18|\Phi(f_1)-\Phi(f_0)|$ for every estimator.

\paragraph*{The transport inverse}
Let $u$ be a smooth potential on $\overline\cO$ with exactly one critical point $x_u$, a nondegenerate minimum, and with $\nabla u$ pointing strictly outward on $\partial\cO$. We use the transport inverse $S_u$ of \eqref{lin:eq:S}, built from the backward flow $X_t$ of $\nabla u$ and its Jacobian $J_t$. The nonlinear constructions need bounds that are uniform over a neighbourhood of potentials. Lemma~\ref{low:S} in Section~\ref{app:lower} provides them: $S_u$ is the unique bounded inverse of $T_u$, that is, $T_uS_uq=q$; it is bounded on $C^k(\overline\cO)$ and on $H^m(\cO)$ for every $k,m\ge0$ (and hence on $H^s$, $s\ge0$, by interpolation); and the bounds are uniform over potentials that are bounded in $C^{k+2}$ (respectively $C^{m+2}$) and uniformly nondegenerate in the sense made precise there. From the formula, $S_uq(x)=0$ whenever the backward orbit of $x$ does not meet $\supp q$; in particular $S_uq$ vanishes near $x_u$ if $\supp q$ does not contain $x_u$, because small sublevel sets of $u$ around $x_u$ are invariant under the backward flow.

\begin{lemma}[Forward-first perturbations]\label{low:transfer}
Let $f_0$ satisfy Assumption~\ref{ass:base}, $s>d/2$, and let $V\in C_c^\infty(\cO)$. Put $h=-S_{u_0}L_{f_0}V$. If $\norm h_{H^s}\le\eta_0$, then $f_0+h\in\cF_s(f_0,\eta_0)$,
\[
 u_{f_0+h}-u_0=L_{f_0+h}^{-1}L_{f_0}V,\qquad \norm{u_{f_0+h}-u_0}_2\le C\norm V_2,
\]
with $C$ as in Lemma~\ref{low:duality}, and $\Phi(f_0+h)-\Phi(f_0)=\ip h\psi$. Moreover $\norm h_{H^s}\le C\norm V_{H^{s+2}}$.
\end{lemma}

\begin{proof}
$T_{u_0}h=-L_{f_0}V$, so \eqref{eq:exact} gives $u_{f_0+h}-u_0=-L_{f_0+h}^{-1}T_{u_0}h=L_{f_0+h}^{-1}L_{f_0}V$. For $a\in L^2$, since $V\in\mathcal V$ and $L_{f_0}$ is symmetric on $\mathcal V$, $\ip{L_{f_0+h}^{-1}L_{f_0}V}a=\ip V{L_{f_0}L_{f_0+h}^{-1}a}\le C\norm V_2\norm a_2$ by Lemma~\ref{low:duality}. The last claim is the $H^s$ bound for $S_{u_0}$.
\end{proof}

\begin{lemma}[Affine approximation]\label{low:BH}
There is $C=C(d)$ such that for every ball $B_r$ and $v\in H^2(B_r)$ there is an affine $\ell$ with $\norm{v-\ell}_{L^2(B_r)}\le Cr^2\norm{D^2v}_{L^2(B_r)}$.
\end{lemma}

This is the Bramble--Hilbert lemma \cite{bramble1970estimation}; a proof by scaling and two Poincar\'e inequalities is recalled in Section~\ref{app:lower}. Finally, for $0<\rho\le1$ and $\varphi\in H^s(\R^d)$, Plancherel gives
\begin{equation}\label{low:eq:scaling}
 \norm{\varphi(\cdot/\rho)}_{H^s(\R^d)}\le\rho^{d/2-s}\norm\varphi_{H^s(\R^d)},
\end{equation}
because $1+|\eta|^2/\rho^2\le\rho^{-2}(1+|\eta|^2)$.

\subsection{Proof of Theorem~\ref{thm:lower}(a) and (c)}\label{low:sec:a}

\begin{proof}[Proof of part (a)]
Translate $x_*$ to the origin, write $H=D^2u_0(0)$ (possibly degenerate) and $\nabla u_0(x)=Hx+r(x)$ with $|r(x)|\le C|x|^2$. Fix an even $\chi_0\in C_c^\infty(B_1)$ with $\chi_0\ge0$, $\chi_0\not\equiv0$, and set
\[
 h_\rho(x)=a_\rho\chi_0(x/\rho),\qquad a_\rho=\gamma\rho^{s-d/2},\qquad 0<\rho\le\rho_0 .
\]
By \eqref{low:eq:scaling}, $\norm{h_\rho}_{H^s}\le C\gamma$, so $f_1=f_0+h_\rho\in\cF_s(f_0,\eta)$ for $\gamma$ small. Split
\begin{gather*}
 T_{u_0}h_\rho=F^{(0)}_\rho+F^{(1)}_\rho,\qquad F^{(1)}_\rho=\diver(h_\rho r),\\
 F^{(0)}_\rho(x)=a_\rho q_0(x/\rho),\qquad q_0(y)=\diver_y(\chi_0(y)Hy).
\end{gather*}
The source $q_0$ has zero mass, and its first moments vanish exactly: $\int y_iq_0(y)\,dy=-\int\chi_0(y)(Hy)_i\,dy=0$ because $\chi_0$ is even. For $v\in\mathcal V$ take the affine $\ell$ of Lemma~\ref{low:BH} on $B_\rho$. Then
\[
 |\ip{F^{(0)}_\rho}v|=|\ip{F^{(0)}_\rho}{v-\ell}|\le\norm{F_\rho^{(0)}}_2\,C\rho^2\norm v_{H^2}\le Ca_\rho\rho^{d/2+2}\norm v_{H^2},
\]
and $|\ip{F^{(1)}_\rho}v|=|\ip{h_\rho r}{\nabla v}|\le\norm{h_\rho r}_2\norm v_{H^1}\le Ca_\rho\rho^{d/2+2}\norm v_{H^2}$. By \eqref{eq:exact} and Lemma~\ref{low:duality},
\[
 \norm{u_{f_1}-u_0}_2\le C\norm{T_{u_0}h_\rho}_{\mathcal V^*}\le C\gamma\rho^{s+2}.
\]
Since $\chi_0$ is even, the linear term of $\psi$ integrates to zero and
\[
 \Phi(f_1)-\Phi(f_0)=a_\rho\rho^d\int\chi_0(y)\psi(\rho y)\,dy=a_\rho\rho^d\Bigl\{\psi(0)\int\chi_0+O(\rho^2)\Bigr\},
\]
whose modulus is at least $c\gamma\rho^{s+d/2}$ for $\rho$ small. Take $\rho=N^{-1/(2s+4)}$ and $\gamma$ so small that \eqref{low:eq:KL} is at most $1/2$. Lemma~\ref{lin:lem:twopoint} gives $\Risk\ge c\rho^{s+d/2}=cN^{-(2s+d)/(4s+8)}$. Only $\nabla u_0(0)=0$ was used, and $h_\rho$ is supported in $B_\rho(x_*)$.
\end{proof}

\begin{proof}[Proof of part (c)]
Put $f_1=f_0+t\psi$ with $t=\gamma N^{-1/2}$. By \eqref{eq:exact} and Lemma~\ref{low:duality}, $\norm{u_{f_1}-u_0}_2\le C|t|\norm{\psi\nabla u_0}_2$, because $|\ip{\diver(\psi\nabla u_0)}v|\le\norm{\psi\nabla u_0}_2\norm v_{H^1}$. The separation is $|t|\norm\psi_2^2$. Apply Lemma~\ref{lin:lem:twopoint} with $\gamma$ small.
\end{proof}

\subsection{Proof of Theorem~\ref{thm:lower}(b)}\label{low:sec:b}

The alternatives are the flow-adapted forward perturbations of Proposition~\ref{lin:prop:construction}; the only new ingredient is that Lemma~\ref{low:transfer} controls the nonlinear forward difference exactly. We also record an integration-by-parts identity used repeatedly below. Let $w$ be the characteristic solution \eqref{lin:eq:w} and $h\in C^1(\overline\cO)$. Since $w=0$ on $\partial\cO$, $|w(x)|\le C(1+|\log|x-x_*||)$ and $|b(x)|\le C|x-x_*|$, integrating over $\cO\setminus B_\epsilon(x_*)$ and letting $\epsilon\to0$ gives
\begin{equation}\label{low:eq:ibp}
 \ip h\psi=\ip h{b\cdot\nabla w}=-\ip{T_{u_0}h}w .
\end{equation}
For $h=-S_{u_0}L_{f_0}V$ with $V\in C_c^\infty(\cO)$ this reads $\ip h\psi=\ip{L_{f_0}V}w$.

\begin{proof}[Proof of Theorem~\ref{thm:lower}(b)]
Write $Q=\tr H$ and $\lambda_*=\lmax(H)$, so that $Q/\lambda_*=\Teff_0$. Since \condC\ fails, Lemma~\ref{lin:lem:normalfast} provides $\mathsf q\in C_c^\infty(\R^d)$ with $\Lambda=\ip{F^{\mathcal N}}{\Delta_E\mathsf q}\ne0$, and, if $\Teff_0=4$, with $\supp\mathsf q$ in a fixed flow annulus and $\Lambda>0$. Let $V_t,h_t$ be as in Proposition~\ref{lin:prop:construction}. Since $A_{f_0}h_t=V_t$ means $T_{u_0}h_t=-L_{f_0}V_t$, Lemma~\ref{low:transfer} applies: as soon as $\norm{h_t}_{H^s}\le\eta$, the conductivity $f_0+h_t$ lies in $\cF_s(f_0,\eta)$, and by \eqref{low:eq:KL}
\[
 \KL(P_{f_0},P_{f_0+h_t})\le CN\norm{V_t}_2^2,\qquad \Phi(f_0+h_t)-\Phi(f_0)=\ip{h_t}\psi .
\]

\emph{Irregular regime ($\Teff_0<4$): a single scale.} Choose $\beta=\gamma e^{-(s+2)\lambda_*t+Qt/2}$ and $t$ with $e^{(s+2)\lambda_*t}=N^{1/2}$. By \eqref{lin:eq:gen-norms}, $\norm{h_t}_{H^s}\le C\gamma$, $\KL\le C\gamma^2$, and
\[
 |\ip{h_t}\psi|\ge c\gamma e^{-(s+2)\lambda_*t+Qt/2+(2\lambda_*-Q)t}=c\gamma e^{-(s+\Teff_0/2)\lambda_*t}=c\gamma N^{-(2s+\Teff_0)/(4s+8)} .
\]
For $\gamma$ small, Lemma~\ref{lin:lem:twopoint} gives the claim.

\emph{Critical regime ($\Teff_0=4$).} Here $Q=4\lambda_*$. Take times $t_j=t_0+j\Delta$ with $\Delta$ so large that the supports of the $V_{t_j}$ are pairwise disjoint, normalise $V_j=V_{t_j}/\norm{V_{t_j}}_2$, and let $H_j=-S_{u_0}L_{f_0}V_j$. By \eqref{lin:eq:gen-norms} with $\beta\asymp e^{Qt_j/2}$, we have $\ip{H_j}\psi\ge c_*>0$ and $\norm{H_j}_{H^s}\le Ce^{(s+2)\lambda_*t_j}$ for $t_0$ large. With $\epsilon=N^{-1/2}$, put
\[
 m=\Bigl\lfloor\frac{\log(1/\epsilon)}{2(s+2)\lambda_*\Delta}\Bigr\rfloor,\qquad V=\frac{\gamma\epsilon}{\sqrt m}\sum_{j=1}^mV_j,\qquad h=-S_{u_0}L_{f_0}V .
\]
Then $\norm V_2=\gamma\epsilon$, so the divergence is at most $C\gamma^2$. By geometric summation, $\norm h_{H^s}\le C\gamma\epsilon m^{-1/2}e^{(s+2)\lambda_*t_m}\le C\gamma\epsilon^{1/2}m^{-1/2}\to0$, and $\ip h\psi\ge c_*\gamma\epsilon\sqrt m\asymp\sqrt{\log N/N}$. Lemmas~\ref{low:transfer} and~\ref{lin:lem:twopoint} complete the proof.
\end{proof}

\begin{remark}\label{low:rem:b}
The perturbations $h_t=-S_{u_0}L_{f_0}V_t$ are not localised: they carry the transport tail $S_{u_0}$ along the flow lines through $\supp V_t$ up to $\partial\cO$. This tail is what distinguishes them from conductivity bumps in flow-adapted boxes, which only give the exponent of part (a). For the canonical model, part (b) is contained in the isotropic-class bounds of Section~\ref{low:sec:iso}, which keep the critical point and Hessian fixed.
\end{remark}

\subsection{Lower bounds on the isotropic class}\label{low:sec:iso}

Throughout this subsection $\cO=B_1$, $f_0=g=1$, $u_0=(|x|^2-1)/(2d)$, and $T_0h=T_{u_0}h=d^{-1}(x\cdot\nabla h+dh)$. The backward flow is $X_tx=e^{-t/d}x$, $J_t=e^{-t}$, and $S_0:=S_{u_0}$ is $S_0q(x)=\int_0^\infty e^{-t}q(e^{-t/d}x)\,dt$. For $\psi\in C_c^\infty(B_1)$ with $c=\psi(0)$, the transport solution of \eqref{eq:transport} is, in polar coordinates $x=r\theta$,
\begin{equation}\label{low:eq:w-ball}
\begin{gathered}
 w(r\theta)=dc\log r+\Psi(\theta)+B(x),\\
 \Psi(\theta)=-d\int_0^1\frac{\psi(t\theta)-c}t\,dt,\qquad B(x)=d\int_0^1\frac{\psi(tx)-c}t\,dt,
\end{gathered}
\end{equation}
with $B$ smooth. If $c=0$, $\Psi$ is the flow profile \eqref{eq:profile} in the parametrisation of $\cL$ by $S^{d-1}$. Condition \condC\ at $f_0=1$ means $c=0$ and $\Psi$ constant.

\begin{remark}[On the class $\cI_s$]\label{low:rem:Is}
Since $\tr D^2u_f(x_*)=g(x_*)/f(x_*)$ at a critical point, every $f\in\cI_s(\eta)$ satisfies $f(0)=1$. For $s>d/2+1$ elliptic regularity gives $u_f\in H^{s+1}\subset C^2$, so the conditions defining $\cI_s(\eta)$ are classical. All alternatives constructed below are smooth. The lower bounds of this subsection are therefore valid for every $s>d/2$ over the smooth subclass
\[
 \cI^\infty_s(\eta)=\bigl\{f\in\cF_s(1,\eta)\cap C^\infty(\overline{B_1}):\ x_*(f)=0,\ D^2u_f(0)=I_d/d\bigr\},
\]
which coincides with $\cI_s(\eta)\cap C^\infty$ when $s>d/2+1$. Lower bounds over $\cI^\infty_s(\eta)$ hold over every larger class, in particular over $\cI_s(\eta)$, over $\cF_s(1,\eta)$, and over the class in Theorem~\ref{thm:isotropic} in which the critical point and the scalar Hessian may vary.
\end{remark}

All alternatives are produced by prescribing the potential and recovering the conductivity by transport. For a smooth potential $v$ with $\diver(F\nabla v)=1$, $F$ is $S_v1$.

\begin{lemma}[Conductivity from a prescribed potential]\label{low:potential}
Let $v_0(x)=|x|^2/(2d)$ on a fixed ball $\mathcal B\supset\overline{B_1}$ and $v=v_0+p$ with $p\in C^\infty$, $\nabla p=0$ near $\partial\mathcal B$, $\nabla p(0)=0$ and $\norm p_{C^{m+2}}\le\delta_m$ sufficiently small. Then $v$ is uniformly strictly convex with unique critical point $0$, $\mathsf f_p:=S_v1$ is smooth and positive, $\diver(\mathsf f_p\nabla v)=1$, and
\begin{align}
 \norm{\mathsf f_p-1}_{C^m}&\le C_m\norm p_{C^{m+2}},\label{low:eq:pot1}\\
 \norm{\mathsf f_p-1+S_0\Delta p}_{C^0}&\le C\norm p_{C^3}^2 .\label{low:eq:pot2}
\end{align}
If $p$ is constant near $0$, then $D^2v(0)=I_d/d$ and $\mathsf f_p=1$ on the set of points whose backward orbits under $\nabla v$ stay in the region where $p$ is constant.
\end{lemma}

The proof is in Section~\ref{app:lower}. If $p$ is compactly supported in $B_1$, then $u:=v_0-1/(2d)+p$ vanishes on $\partial B_1$ and, by uniqueness for the Dirichlet problem, $u_{\mathsf f_p}=u$ on $B_1$. If $p$ is constant near $0$, then $x_*(\mathsf f_p)=0$ and $D^2u_{\mathsf f_p}(0)=I_d/d$.

\begin{proposition}[Single-scale bound inside $\cI_s$]\label{low:iso-single}
Let $s>d/2$. If $\psi(0)\ne0$, or if $\psi(0)=0$ and $\Psi$ is nonconstant, then $\Risk(\cI_s(\eta))\ge cN^{-(2s+d)/(4s+8)}$.
\end{proposition}

\begin{proof}
Choose $q\in C_c^\infty(B_2)$, constant near $0$, specified below. For $\rho$ small let
\begin{equation}\label{low:eq:urho}
 u_\rho(x)=u_0(x)+a\rho^2q(x/\rho),\qquad a=\gamma\rho^{s-d/2}.
\end{equation}
In the variable $y=x/\rho$ let $\mathsf f_a$ be the bounded positive solution on $\R^d$ of $\diver_y\{\mathsf f_a(y)(y/d+a\nabla q(y))\}=1$, and set $f_\rho(x)=\mathsf f_a(x/\rho)$. Then $\diver(f_\rho\nabla u_\rho)=1$ on $B_1$, and $u_\rho=0$ on $\partial B_1$, so $u_{f_\rho}=u_\rho$. Outside $B_2$ the field is radial and the transport equation $d^{-1}(y\cdot\nabla+d)(\mathsf f_a-1)=0$ gives
\begin{equation}\label{low:eq:tail}
 \mathsf f_a(r\theta)-1=r^{-d}B_a(\theta),\qquad r\ge2,
\end{equation}
where $B_a=2^d(\mathsf f_a-1)(2\,\cdot)$ is $O(a)$ in every $C^k(S^{d-1})$ by Lemma~\ref{low:potential} applied on a ball containing $B_3$. The tail \eqref{low:eq:tail} and all its derivatives are square integrable at infinity, so $\norm{\mathsf f_a-1}_{H^s(\R^d)}\le Ca$, and by \eqref{low:eq:scaling}
\[
 \norm{f_\rho-1}_{H^s(B_1)}\le Ca\rho^{d/2-s}=C\gamma .
\]
For $\gamma$ small, $u_\rho$ is uniformly strictly convex with unique critical point $0$ and $D^2u_\rho(0)=I_d/d$ because $q$ is constant near $0$. Hence $f_\rho\in\cI_s(\eta)$. The forward difference is exactly
\begin{equation}\label{low:eq:fwd}
 \norm{u_{f_\rho}-u_0}_2=a\rho^{d/2+2}\norm q_2=\gamma\rho^{s+2}\norm q_2 .
\end{equation}

\emph{Case $\psi(0)\ne0$.} Take $q$ radial, nonnegative, equal to a positive constant near $0$. The radial equation $r^{1-d}\{r^{d-1}\mathsf f_a(r/d+aq'(r))\}'=1$ gives
\[
 \mathsf f_a(r)=\frac1{1+da\,q'(r)/r},
\]
so $\mathsf f_a-1$ is supported in a fixed annulus, and $\mathsf f_a-1=aH_q+O(a^2)$ uniformly with $H_q=-dq'(r)/r$. Integration by parts gives $\int_{\R^2}H_q=4\pi q(0)$ for $d=2$ and $\int_{\R^d}H_q=d(d-2)|S^{d-1}|\int_0^\infty q(r)r^{d-3}\,dr$ for $d\ge3$; both are positive. By radial symmetry the linear term of $\psi$ integrates to zero, and
\[
 \ip{f_\rho-1}\psi=\rho^d\int(\mathsf f_a(y)-1)\psi(\rho y)\,dy=a\rho^d\Bigl\{\psi(0)\int H_q+O(\rho^2)+O(a)\Bigr\}.
\]

\emph{Case $\psi(0)=0$, $\Psi$ nonconstant.} Take $q(r\theta)=\zeta(r)A(\theta)$ with $A=\Delta_{S^{d-1}}\Psi$ and $0\le\zeta\in C_c^\infty((1,2))$, $\zeta\not\equiv0$. Then $q=0$ near $0$. Writing $\mathsf f_a$ as $S_{v}1$ for $v=v_0+aq$, Lemma~\ref{low:potential} gives $\mathsf f_a-1=aH_q+O(a^2)$ in every $C^k$ on a fixed ball, with $H_q=-S_0\Delta q$. For $|y|=R>2$, the substitution $r=e^{-t/d}R$ gives $H_q(R\theta)=-dR^{-d}\int_0^\infty r^{d-1}\Delta q(r\theta)\,dr$. Since $\Delta q=r^{1-d}(r^{d-1}\zeta')'A+r^{-2}\zeta\Delta_{S^{d-1}}A$ and the first part integrates to zero,
\[
 B_a=aB_1+O_{C^k}(a^2),\qquad B_1=-dC_\zeta\Delta_{S^{d-1}}A,\qquad C_\zeta=\int_0^\infty r^{d-3}\zeta(r)\,dr>0 .
\]
Split $\ip{f_\rho-1}\psi=\rho^d\int(\mathsf f_a(y)-1)\psi(\rho y)\,dy$ at $|y|=2$. Since $|\psi(x)|\le C|x|$ and $\mathsf f_a-1=O(a)$, the inner part is $O(a\rho^{d+1})$. By \eqref{low:eq:tail} the outer part equals
\[
 \rho^d\int_{S^{d-1}}B_a(\theta)\int_{2\rho}^\infty\frac{\psi(t\theta)}t\,dt\,d\theta=\rho^d\int_{S^{d-1}}B_a(\theta)\Bigl\{-\frac{\Psi(\theta)}d+O(\rho)\Bigr\}d\theta .
\]
Using self-adjointness of $\Delta_{S^{d-1}}$,
\[
 \ip{f_\rho-1}\psi=a\rho^d\Bigl\{C_\zeta\norm{\Delta_{S^{d-1}}\Psi}_{L^2(S^{d-1})}^2+O(\rho)+O(a)\Bigr\},
\]
and the leading coefficient is positive because $\Psi$ is not constant.

In both cases the separation is at least $ca\rho^d=c\gamma\rho^{s+d/2}$. Take $\rho=N^{-1/(2s+4)}$, use \eqref{low:eq:fwd}, \eqref{low:eq:KL} and Lemma~\ref{lin:lem:twopoint}.
\end{proof}

\begin{proposition}[The logarithm in $d=4$ inside $\cI_s$]\label{low:iso-log}
Let $d=4$, $s>2$, and suppose \condC\ fails at $f_0=1$. Then $\Risk(\cI_s(\eta))\ge c\sqrt{\log N/N}$.
\end{proposition}

\begin{proof}
By \eqref{low:eq:w-ball} and $\Delta\log r=2r^{-2}$ in $d=4$,
\begin{equation}\label{low:eq:Dw4}
 \Delta w=r^{-2}\mathsf a(\theta)+\Delta B,\qquad \mathsf a=8c+\Delta_{S^3}\Psi\quad\text{on }B_1\setminus\{0\}.
\end{equation}
If $\mathsf a\equiv0$, its spherical mean gives $c=0$ and then $\Delta_{S^3}\Psi=0$, so $\Psi$ is constant and \condC\ holds. Hence $\mathsf a\not\equiv0$. Let $0\le\zeta\in C_c^\infty((1,2))$, $\zeta\not\equiv0$, $v(y)=\zeta(|y|)\mathsf a(y/|y|)$, and $r_j=r_04^{-j}$, $v_j(x)=r_j^{-2}v(x/r_j)$. The $v_j$ have disjoint supports, $\norm{v_j}_2=\norm v_2$, and by \eqref{low:eq:Dw4}
\[
 \int v_j\Delta w=K+O(r_j^2),\qquad K=\int_{\R^4}\zeta(|y|)\,|y|^{-2}\mathsf a(y/|y|)^2\,dy>0 .
\]
Fix $r_0$ so that $\int v_j\Delta w\ge K/2$ for all $j$. Let $\epsilon=N^{-1/2}$, fix an integer $m\ge s$ and $0<c_0<\min\{1/(m+4),1/10\}$, and put
\[
 J=\Bigl\lfloor\frac{c_0\log(1/\epsilon)}{\log4}\Bigr\rfloor,\qquad\delta=\frac{\gamma\epsilon}{\sqrt J},\qquad V=\delta\sum_{j=1}^Jv_j .
\]
For every integer $k$, $\norm V_{C^k}\le C_k\delta r_J^{-k-2}$ and $r_J^{-1}\le C\epsilon^{-c_0}$, so $\norm V_{C^{m+2}}\le C\gamma\epsilon^{1-c_0(m+4)}\to0$. The potential $u=u_0+V$ vanishes on $\partial B_1$, equals $u_0$ near $0$, and, by Lemma~\ref{low:potential}, $f=S_u1$ satisfies $u_f=u$, $x_*(f)=0$, $D^2u_f(0)=I_4/4$, and $\norm{f-1}_{H^s}\le C\norm{f-1}_{C^m}\le C\norm V_{C^{m+2}}\to0$. Thus $f\in\cI_s(\eta)$ for $N$ large, and $\KL(P_1,P_f)=\frac N2\norm V_2^2=\frac{\gamma^2}2\norm v_2^2$.

Let $h_1=-S_0\Delta V$. It vanishes near $0$, so \eqref{low:eq:ibp} applies and $\ip{h_1}\psi=\ip{\Delta V}w=\int V\Delta w\ge(K/2)\delta J$. By \eqref{low:eq:pot2},
\[
 |\ip{f-1-h_1}\psi|\le C\norm V_{C^3}^2\le C\delta^2r_J^{-10},
\]
and the ratio of this to $\delta J$ is at most $C\gamma\epsilon^{1-10c_0}J^{-3/2}\to0$. The separation is therefore at least $c\delta J\asymp\gamma\epsilon\sqrt J\asymp\sqrt{\log N/N}$, and Lemma~\ref{lin:lem:twopoint} applies.
\end{proof}

\begin{proposition}[Parametric bound inside $\cI_s$]\label{low:iso-param}
For every nonzero $\psi\in C_c^\infty(B_1)$ and $s>d/2$, $\Risk(\cI_s(\eta))\ge cN^{-1/2}$.
\end{proposition}

\begin{proof}
Since $\psi$ has compact support, $w=0$ near $\partial B_1$. If $\Delta w\equiv0$ on the connected set $B_1\setminus\{0\}$, unique continuation from that collar gives $w\equiv0$ and hence $\psi=b\cdot\nabla w=0$. So there is $q\in C_c^\infty(B_1\setminus\{0\})$ with $\int q\Delta w\ne0$. Let $u_t=u_0+tq$ and $f_t=S_{u_t}1$. For $|t|$ small, $f_t\in\cI_s(\eta)$ and $u_{f_t}-u_0=tq$. By \eqref{low:eq:pot2} and \eqref{low:eq:ibp}, $\Phi(f_t)-\Phi(1)=-t\ip{S_0\Delta q}\psi+O(t^2)=t\int q\Delta w+O(t^2)$. Take $t=\gamma N^{-1/2}$.
\end{proof}

The remaining case is a target satisfying \condC\ at the reference. We first perturb the reference so that its flow profile becomes nonconstant, while preserving the critical point and Hessian (Proposition~\ref{low:activate}). We then transplant the preceding lower-bound construction to this nearby conductivity (Lemma~\ref{low:transplant}). Statistically, this means that compatibility at one conductivity does not protect a uniform risk bound: the least favourable alternatives sit next to a nearby incompatible conductivity, and this yields the uniform lower bound for every nonzero target.

\begin{lemma}[Transplanted lower bounds]\label{low:transplant}
Let $\bar u$ be a smooth, uniformly strictly convex potential on $B_1$ with $\bar u=0$ on $\partial B_1$, $\nabla\bar u$ outward on $\partial B_1$, and $\bar u=u_0$ on a ball $B_{\mathsf b}$ and near $\partial B_1$. Let $\bar f=S_{\bar u}1$, so that $u_{\bar f}=\bar u$ and $\bar f=1$ on $B_{\mathsf b}$. Let $\bar w$ solve $\nabla\bar u\cdot\nabla\bar w=\psi$, $\bar w|_{\partial B_1}=0$, so that on $B_{\mathsf b}\setminus\{0\}$
\[
 \bar w=d\psi(0)\log r+\bar\Psi(\theta)+\bar B(x),\qquad \bar B\text{ smooth}.
\]
If $\psi(0)\ne0$ or $\bar\Psi$ is nonconstant, then for every $\eta'>0$ the conclusions of Propositions~\ref{low:iso-single} and (for $d=4$) \ref{low:iso-log} hold with $\cI_s(\eta)$ replaced by $\{f\in\cF_s(\bar f,\eta'):x_*(f)=0,\ D^2u_f(0)=I_d/d\}$.
\end{lemma}

\begin{proof}
The expansion of $\bar w$ holds because $\nabla\bar u=x/d$ on $B_{\mathsf b}$, and integrating along rays from the sphere $|x|=\mathsf b/2$ gives \eqref{low:eq:w-ball} with $\Psi$ replaced by $\bar\Psi(\theta)=\bar w(\tfrac{\mathsf b}2\theta)-d\psi(0)\log\frac{\mathsf b}2-\bar B(\frac{\mathsf b}2\theta)$ for a smooth $\bar B$.

\emph{Single scale.} Use $u_\rho=\bar u+a\rho^2q(x/\rho)$ with $4\rho<\mathsf b$ and $q$ as in the proof of Proposition~\ref{low:iso-single} (with $\Psi$ replaced by $\bar\Psi$). Backward orbits of $\nabla u_\rho$ starting in $B_{\mathsf b}$ stay in $B_{\mathsf b}$, so the conductivity $f_\rho=S_{u_\rho}1$ equals $\mathsf f_a(x/\rho)$ on $B_{\mathsf b}$, with $\mathsf f_a$ as before. Outside $B_{\mathsf b/2}$, $h=f_\rho-\bar f$ solves $T_{\bar u}h=0$, so it is obtained from its trace $h(\tfrac{\mathsf b}2\theta)=(2\rho/\mathsf b)^dB_a(\theta)$ by transport along the forward flow of $\nabla\bar u$, which crosses the region between $\partial B_{\mathsf b/2}$ and $\partial B_1$ in a uniformly bounded time with bounded derivatives. This defines a fixed linear extension operator, bounded from $C^k(S^{d-1})$ into $C^k$ of the outer region for every $k$. Hence the outer part of $h$ is $O(a\rho^d)$ in every $C^k$, and its deviation from the term linear in $a$ is $O(a^2\rho^d)$. Together with \eqref{low:eq:scaling} in the core, $\norm{f_\rho-\bar f}_{H^s}\le Ca\rho^{d/2-s}=C\gamma$. The forward difference is exactly $a\rho^2q(x/\rho)$, as in \eqref{low:eq:fwd}.

If $\psi(0)\ne0$, $q$ is radial and $\mathsf f_a-1$ is supported in an annulus of radius $<2$, so $h$ vanishes outside $B_{2\rho}$ and the computation of the separation is identical to the case $\psi(0)\ne0$ of Proposition~\ref{low:iso-single}. If $\psi(0)=0$, the linear part $h^{(1)}$ of $h$ in $a$ is $h^{(1)}=-S_{\bar u}L_{\bar f}(a\rho^2q(\cdot/\rho))$, and by \eqref{low:eq:ibp} for the base $\bar f$ and $\bar f=1$ on the support of $q(\cdot/\rho)$,
\[
 \ip{h^{(1)}}\psi=a\rho^2\int q(x/\rho)\Delta\bar w(x)\,dx=a\rho^d\Bigl\{C_\zeta\norm{\Delta_{S^{d-1}}\bar\Psi}_2^2+O(\rho^2)\Bigr\},
\]
using $\Delta\bar w=r^{-2}\Delta_{S^{d-1}}\bar\Psi+\Delta\bar B$ on $B_{\mathsf b}\setminus\{0\}$. The nonlinear remainder $h-h^{(1)}$ contributes $O(a^2\rho^d)$ to the functional: outside $B_{\mathsf b/2}$ by the extension estimate, and inside by \eqref{low:eq:tail} and $|\psi(x)|\le C|x|$, exactly as in the proof of Proposition~\ref{low:iso-single}. The two-point argument concludes as before.

\emph{Dimension four.} Place the annular potentials $v_j$ of Proposition~\ref{low:iso-log} inside $B_{\mathsf b}$, built from $\bar{\mathsf a}=8\psi(0)+\Delta_{S^3}\bar\Psi\not\equiv0$, put $u=\bar u+V$ and $f=S_u1$. Since $\bar f=1$ on $\supp V$, $T_u(f-\bar f)=-\Delta V$, so with $h_1=-S_{\bar u}\Delta V$ we have $f-\bar f-h_1=-S_u\diver(h_1\nabla V)$. Lemma~\ref{low:S} for the fixed strictly convex reference $\bar u$ gives the same $H^s$ control and the same $O(\norm V_{C^3}^2)$ functional remainder as in Lemma~\ref{low:potential}. The forward difference is exactly $V$ and $\ip{h_1}\psi=\int V\Delta\bar w$. The choices of Proposition~\ref{low:iso-log} apply verbatim.
\end{proof}

\begin{proposition}[Activation of the flow profile]\label{low:activate}
Let $\psi\ne0$ satisfy \condC\ at $f_0=1$. For every $m$ and every $\delta>0$ there is a smooth, uniformly strictly convex potential $\bar u$ on $B_1$, equal to $u_0$ near $0$ and near $\partial B_1$, with $\norm{\bar u-u_0}_{C^m}<\delta$, whose flow profile for $\psi$ is nonconstant. Its conductivity $\bar f=S_{\bar u}1$ equals $1$ near $0$, satisfies $\norm{\bar f-1}_{C^{m-2}}\le C\delta$, has critical point $0$ and Hessian $I_d/d$, and \condC\ fails at $\bar f$.
\end{proposition}

\begin{proof}
By \condC, $\psi(0)=0$ and $w_0:=w$ is smooth. Since $\psi\ne0$, there are $\theta_0$ and $r$ with $\psi(r\theta_0)\ne0$. The function $\mathsf F(r)=\partial_rw_0(r\theta_0)/r=d\psi(r\theta_0)/r^2$ is smooth on $(0,1)$, vanishes near $r=1$ and is not identically zero, hence not constant; pick $0\le\zeta\in C_c^\infty((r_1,r_2))$, $0<r_1<r_2<1$, with $\int\zeta'\mathsf F\,dr=-\int\zeta\mathsf F'\,dr\ne0$. Let $\chi\in C^\infty(S^{d-1})$ equal $1$ near $\theta_0$ and vanish on an open set of directions, and $q(r\theta)=\zeta(r)\chi(\theta)$. Put $u_t=u_0+tq$ and let $w_t$ solve $\nabla u_t\cdot\nabla w_t=\psi$, $w_t|_{\partial B_1}=0$. Away from the critical point the characteristics of $\nabla u_t$ between the spheres $|x|=r_1/2$ and $|x|=1$ are transversal to the spheres and have uniformly bounded travel times, so $w_t$ on $|x|=r_1/2$ depends smoothly on $t$, and so does the profile $\Psi_t$ (parametrised by the crossing point with this sphere, where the field is still radial). Differentiating the transport equation at $t=0$,
\[
 \tfrac1d\,x\cdot\nabla\dot w=-\nabla q\cdot\nabla w_0,\qquad \dot w|_{\partial B_1}=0,
 \qquad\text{so}\qquad \dot\Psi(\theta)=d\int_0^1\frac{(\nabla q\cdot\nabla w_0)(r\theta)}r\,dr .
\]
At $\theta_0$, $\nabla_{S^{d-1}}\chi=0$ and $\chi=1$, so $\nabla q\cdot\nabla w_0=\zeta'(r)\partial_rw_0$ and $\dot\Psi(\theta_0)=d\int\zeta'\mathsf F\,dr\ne0$, while $\dot\Psi=0$ on directions where $\chi$ vanishes. Thus $\dot\Psi$ is nonconstant. Since $\Psi_0$ is constant, $\Psi_t=\Psi_0+t\dot\Psi+O(t^2)$ is nonconstant for all small $t\ne0$. For such $t$ put $\bar u=u_t$; it is uniformly strictly convex, equal to $u_0$ near $0$ and near $\partial B_1$, and Lemma~\ref{low:potential} gives the claims on $\bar f$. Near $0$ the only bounded solution of $T_{u_0}\bar f=1$ along backward orbits is $\bar f=1$.
\end{proof}

\begin{corollary}[Lower bounds on the full isotropic class]\label{low:iso-uniform}
For every nonzero $\psi\in C_c^\infty(B_1)$, $s>d/2$ and sufficiently small $\eta>0$,
\[
 \Risk(\cI_s(\eta))\ge c\,\rate(s,d)=c\begin{cases}N^{-(2s+d)/(4s+8)},&d=2,3,\\ \sqrt{\log N/N},&d=4,\\ N^{-1/2},&d\ge5,\end{cases}
\]
whether or not \condC\ holds at $f_0=1$. The same lower bound holds over every class containing $\cI^\infty_s(\eta)$, in particular over $\cI_s(\eta)$ and over the class of Theorem~\ref{thm:isotropic} in which the critical point and the scalar Hessian may vary (Remark~\ref{low:rem:Is}). This is the lower bound in Theorem~\ref{thm:isotropic}; the upper bound is proved in Section~\ref{sec:upper-proofs}.
\end{corollary}

\begin{proof}
For $d\ge5$ use Proposition~\ref{low:iso-param}. For $d\le4$, if \condC\ fails at $f_0=1$, use Propositions~\ref{low:iso-single} ($d=2,3$) and~\ref{low:iso-log} ($d=4$). If \condC\ holds, Proposition~\ref{low:activate} with $m>s+2$ provides $\bar f$ with $\norm{\bar f-1}_{H^s}\le\eta/2$, critical point $0$, Hessian $I_d/d$ and nonconstant profile. Lemma~\ref{low:transplant} with $\eta'=\eta/2$ produces alternatives in $\{f\in\cF_s(\bar f,\eta/2):x_*(f)=0,\ D^2u_f(0)=I_d/d\}\subset\cI_s(\eta)$, and the two-point bounds over this set bound $\Risk(\cI_s(\eta))$ from below.
\end{proof}

\subsection{Proof of Proposition~\ref{prop:instability}}\label{low:sec:instability}

\begin{proof}[Proof of part (a)]
Let $\chi\in C_c^\infty(B_1)$ be radial and equal to $1$ on $B_{1/2}$, and $v_\varepsilon=u_0+\varepsilon\chi(x)(x_1^2-x_2^2)$. For $\varepsilon$ small, $v_\varepsilon$ is uniformly strictly convex, vanishes on $\partial B_1$, has unique critical point $0$, and $\nabla v_\varepsilon$ points outward on $\partial B_1$. Let $f_\varepsilon=S_{v_\varepsilon}1$. Then $T_{v_\varepsilon}(f_\varepsilon-1)=1-\Delta v_\varepsilon=-\varepsilon\Delta\{\chi(x_1^2-x_2^2)\}$, whose right-hand side is $O(\varepsilon)$ in every $C^k$ and vanishes on $B_{1/2}$ because $x_1^2-x_2^2$ is harmonic. By Lemma~\ref{low:S}, uniformly in small $\varepsilon$, $\norm{f_\varepsilon-1}_{H^s}\le C_s\varepsilon$ for every $s$; in particular $f_\varepsilon$ is positive. Backward orbits starting in $B_{1/2}$ stay there, so $f_\varepsilon=1$ on $B_{1/2}$. By uniqueness for the Dirichlet problem, $u_{f_\varepsilon}=v_\varepsilon$, so the critical point is $0$ and
\[
 D^2u_{f_\varepsilon}(0)=\tfrac1dI_d+2\varepsilon\operatorname{diag}(1,-1,0,\dots,0),\qquad \Teff\bigl(D^2u_{f_\varepsilon}(0)\bigr)=\frac{1}{1/d+2\varepsilon}=\frac d{1+2d\varepsilon}.
\qedhere
\]
\end{proof}

\begin{proof}[Proof of part (b)]
This is Proposition~\ref{low:activate}: the conductivity $\bar f$ can be made arbitrarily close to $1$ in every $H^s$ norm, has the same critical point $0$ and Hessian $I_d/d$, and \condC\ fails at $\bar f$ because its flow profile is nonconstant.
\end{proof}

\begin{proof}[Proof of the consequences]
Since $\cI_s(\eta)$ is contained in $\cF_s(1,\eta)$, Corollary~\ref{low:iso-uniform} gives $\Risk(\cF_s(1,\eta))\ge c\,\rate(s,d)$ for every nonzero $\psi$, and $\rate(s,d)/N^{-1/2}\to\infty$ for $d\le4$. If $d\in\{2,3\}$ and $\psi(0)\ne0$, fix $\varepsilon>0$ with $\norm{f_\varepsilon-1}_{H^s}\le\eta/2$. The base $f_\varepsilon$ satisfies Assumption~\ref{ass:base} (with $g=1$), \condC\ fails at $f_\varepsilon$ because $\psi(x_*(f_\varepsilon))=\psi(0)\ne0$, and $\Teff_\varepsilon=d/(1+2d\varepsilon)<d\le4$. Theorem~\ref{thm:lower}(b) at the base $f_\varepsilon$ and $\cF_s(f_\varepsilon,\eta/2)\subset\cF_s(1,\eta)$ give
\[
 \Risk(\cF_s(1,\eta))\ge cN^{-(2s+\Teff_\varepsilon)/(4s+8)},
\]
which is polynomially larger than $\rate(s,d)=N^{-(2s+d)/(4s+8)}$. The same argument in $d=4$ shows that, for $\psi(0)\ne0$, the risk over $\cF_s(1,\eta)$ is polynomially larger than $\sqrt{\log N/N}$.
\end{proof}

\begin{remark}[Moving the critical point]\label{low:rem:shift}
Condition \condC\ at the centre can also be destroyed by moving the critical point, with perturbations supported away from $\partial\cO$. On the canonical model let $w_0=x_1\chi(|x|)$ with a radial cutoff $\chi$ equal to $1$ near $0$, and $\psi=d^{-1}x\cdot\nabla w_0$. Then $\psi\in C_c^\infty(B_1)$, $w_0$ is the transport solution and \condC\ holds. For $f_t=1+t\chi_1$ with $\chi_1\in C_c^\infty(B_1)$ and $\chi_1(x)=x_1$ on $B_{1-\delta}$, the implicit function theorem gives $x_*(f_t)=-t\mathsf k_de_1+O(t^2)$. For $\chi_1=x_1$ on all of $B_1$ one computes $A_1x_1=\frac{d+1}{2d(d+2)}x_1(1-|x|^2)$ and $\mathsf k_d=\frac{d+1}{2(d+2)}$; by continuity $\mathsf k_d\ne0$ also for $\delta$ small. Hence $\psi(x_*(f_t))=-t\mathsf k_d/d+O(t^2)\ne0$ for small $t\ne0$. Theorem~\ref{thm:lower}(a) at the base $f_t$ gives the lower bound $N^{-(2s+d)/(4s+8)}$ over every fixed neighbourhood of $1$, also when $f$ is known near $\partial B_1$.
\end{remark}

\subsection{Interior perturbation classes}\label{low:sec:collar}

For a compact $K\subset\cO$ consider the linearised experiment $Y=A_{f_0}h+N^{-1/2}\mathbb W$ over
\[
 \cF^K_s(\eta)=\{h\in H^s(\cO):\ \supp h\subset K,\ \norm h_{H^s}\le\eta\}.
\]

\begin{proposition}[Only $\psi(x_*)$ matters for interior perturbations]\label{low:collar}
Let Assumption~\ref{ass:base} hold and $K\subset\cO$ be compact.
\begin{enumerate}
\item[(i)] If $\psi(x_*)=0$, there is $z\in C^\infty(\overline\cO)$ with $\ip{A_{f_0}h}z=\ip h\psi$ for all $h\in L^2$ with $\supp h\subset K$. The estimator $\ip Yz$ is unbiased over $\cF^K_s(\eta)$ with standard deviation $N^{-1/2}\norm z_2$, whatever $\Psi$, $d$ and $H$ are.
\item[(ii)] If $\psi(x_*)\ne0$ and $x_*\in\operatorname{int}K$, the minimax risk over $\cF^K_s(\eta)$ is at least $cN^{-(2s+d)/(4s+8)}$.
\end{enumerate}
\end{proposition}

\begin{proof}
(i) Since $\psi(x_*)=0$, the function $v_\psi$ of \eqref{lin:eq:vpsi} is smooth on $\overline\cO$ and satisfies $b\cdot\nabla v_\psi=\psi$. Choose $\chi_\partial\in C^\infty(\overline\cO)$ with $\chi_\partial=0$ on a neighbourhood of $K$ and $\chi_\partial=1$ near $\partial\cO$, and set $w_K=(1-\chi_\partial)v_\psi$ and $z=L_{f_0}w_K$. Then $w_K\in\mathcal V$ and $b\cdot\nabla w_K=\psi$ near $K$. For $\supp h\subset K$,
\[
 \ip{A_{f_0}h}z=-\ip{L_{f_0}^{-1}\diver(hb)}{L_{f_0}w_K}=-\ip{\diver(hb)}{w_K}=\ip h{b\cdot\nabla w_K}=\ip h\psi .
\]
(ii) The perturbations in the proof of Theorem~\ref{thm:lower}(a) are supported in $B_\rho(x_*)\subset K$; the linearised version of that argument is identical, with $A_{f_0}h_\rho=-L_{f_0}^{-1}T_{u_0}h_\rho$ in place of \eqref{eq:exact}.
\end{proof}

\begin{remark}\label{low:rem:collar}
Part (i) shows that the flow-profile part of \condC\ is created by perturbations that reach the boundary: the nonconstant profile can be removed by subtracting from $w$ a function constant along flow lines, and this is possible exactly when the boundary values on those flow lines are not seen by $h$. It is not enough to know the boundary trace of $f$, as in \cite{nickl2023information}; what matters is whether perturbations may come arbitrarily close to $\partial\cO$. In the nonlinear experiment over a fixed neighbourhood the critical point moves with $f$ (Remark~\ref{low:rem:shift}). In $d=2,3$, part (ii) applied at each member of the class shows that uniform parametric estimation over an interior class containing the local alternatives of part (ii) around each of its members requires $\psi(x_*(f))=0$ at every critical point $x_*(f)$ of the class. This is a condition at the admissible critical points, which is weaker than vanishing on an open neighbourhood of them. In $d=4$ the bound of part (ii) equals $N^{-1/2}$, and in $d\ge5$ it is weaker than the parametric bound, so it gives no such conclusion there. For $\psi(x_*)\ne0$ the radial constructions of Propositions~\ref{low:iso-single} and~\ref{low:iso-log} (with $\mathsf a$ replaced by its radial part $8\psi(0)$) are supported in a small ball around $0$, so they also apply to the isotropic class intersected with an interior class.
\end{remark}

\section{Proofs of the nonlinear upper bounds}\label{sec:upper-proofs}

All upper bounds use a prediction pilot and a one-step correction. Given a pilot $p$ from the first half of the observation and $w\in H^2(\cO)\cap H^1_0(\cO)$ determined by $p$, the second half yields
\[
 \widehat\Phi=\ip p\psi+Y_2(L_pw)-\ip{u_p}{L_pw}.
\]
Writing $h=f-p$, $v=u_f-u_p$ and $r=\psi-\nabla u_p\cdot\nabla w$, Lemma~\ref{up:lem:identity} gives
\[
 \E_f|\widehat\Phi-\Phi(f)|\lesssim\E_f|\ip hr|+N^{-1/2}\E_f\norm{L_pw}_2+\E_f\Bigl|\int_\cO h\nabla v\cdot\nabla w\Bigr| .
\]
These terms are the regularisation bias, the noise contribution and the nonlinear remainder. The pilot (Lemma~\ref{up:lem:pilot}) and the transport stability estimate (Lemma~\ref{up:lem:stability}) control the last term, together with a fractional Hardy inequality (Lemma~\ref{up:lem:hardy}). The geometry enters through the choice of $w$: on isotropic classes we regularise at the pilot's critical point and cancel residual moments (Lemma~\ref{up:lem:reg}); when the exact representer is uniformly square integrable, we use it directly (Lemma~\ref{up:lem:exact-iso}, Proposition~\ref{up:prop:uniformC} and Section~\ref{up:sec:full-regular}). Auxiliary results used repeatedly are collected in Section~\ref{up:app}.

\subsection{Standing setting and uniform geometry}\label{up:sec:setting}

Throughout this section Assumption~\ref{ass:base} holds and
\begin{equation}\label{up:eq:s}
 s>\frac d2+2 .
\end{equation}
Let $\cF=\cF_s(f_0,\eta)$, and let $\mathcal G\subset\cF$ be a nonempty class; additional conditions on $\mathcal G$ are imposed in each result. For $p\in\cF$ we write $u_p=u_{f}$ with $f=p$, and $x_p$ for the critical point of $u_p$ when it is unique.

\begin{lemma}[Uniform regularity and geometry]\label{up:lem:uniform}
There are $\eta_0>0$ and $C<\infty$ such that the following holds for all $\eta\le\eta_0$.
\begin{enumerate}
\item[(a)] $\sup_{p\in\cF}\bigl(\norm p_{H^s}+\norm{u_p}_{H^{s+1}}\bigr)\le C$ and $\sup_{p\in\cF}\norm{u_p-u_0}_{C^3(\overline\cO)}\le C\eta$; all $p\in\cF$ are bounded below by a fixed positive constant.
\item[(b)] There are regular values $a<b<0$ of $u_0$ with $\supp\psi\subset\{u_0<a\}$, such that for $K=\{u_0\le b\}\Subset\cO$ and every $p\in\cF$: $u_p$ has exactly one critical point $x_p$, which is a nondegenerate minimum lying in the interior of $K$; $\nabla u_p$ points strictly outward on $\partial K$; and on a fixed ball $B_{r_0}(x_p)$ the Hessian satisfies $D^2u_p\ge mI_d$ for a fixed $m>0$.
\item[(c)] Every backward trajectory of $\dot X=-\nabla u_p(X)$ started in $K$ stays in $K$ and enters $B_{r_0/2}(x_p)$ after a time bounded uniformly in $p$ and in the starting point.
\item[(d)] Let $w_p$ be the solution of $\nabla u_p\cdot\nabla w_p=\psi$ in $\cO\setminus\{x_p\}$ with $w_p=0$ on $\partial\cO$, obtained by integrating $\psi$ along characteristics. Then $w_p$ vanishes outside a compact subset of the interior of $K$ that does not depend on $p$.
\end{enumerate}
\end{lemma}

\begin{proof}
(a) Write $p=f_0+\alpha$ and $v=u_p-u_0$. Subtracting the equations for $u_p$ and $u_0$ gives
\[
 v=-L_{f_0}^{-1}\diver(\alpha\nabla u_0)-L_{f_0}^{-1}\diver(\alpha\nabla v).
\]
The Dirichlet inverse $L_{f_0}^{-1}$ maps $H^{s-1}$ boundedly into $H^{s+1}\cap H^1_0$, since $f_0$ and $\partial\cO$ are smooth. Since $s>d/2$, $H^s$ is an algebra, so $v\mapsto L_{f_0}^{-1}\diver(\alpha\nabla v)$ has operator norm at most $C\norm\alpha_{H^s}\le C\eta$ on $H^{s+1}$. A Neumann series gives $\norm v_{H^{s+1}}\le C\eta$. By \eqref{up:eq:s}, $s+1>d/2+3$, so Sobolev embedding gives the $C^3$ bound. The lower bound on $p$ follows from $H^s\hookrightarrow C^0$.

(b) These properties hold for $u_0$: $u_0<0$ in $\cO$, $u_0$ has a unique critical point, which is a nondegenerate minimum, and $u_0$ has no critical values in $[a,0)$. They are stable under $C^2$-small perturbations, which are provided by (a). Since $u_0$ has no other critical point, $|\nabla u_0|$ is bounded below on $K\setminus B_{r_0/4}(x_*)$, and the same holds for $u_p$ with a smaller constant.

(c) Along a backward trajectory, $\frac{d}{dt}u_p(X_t)=-|\nabla u_p(X_t)|^2$. Since $\nabla u_p$ points outward on $\partial K$, the trajectory cannot leave $K$. The quantity $u_p$ is bounded on $K$ and decreases at a rate bounded below outside $B_{r_0/2}(x_p)$, which gives the entrance-time bound.

(d) $u_p$ increases along forward characteristics. Hence a point with $u_p(x)>a+C\eta$ never reaches $\supp\psi\subset\{u_0<a\}\subset\{u_p<a+C\eta\}$ in forward time, and $w_p(x)=0$ there. For $\eta$ small, $\{u_p\le a+C\eta\}$ lies in a fixed compact subset of the interior of $K$.
\end{proof}

From now on $\eta\le\eta_0$, and $K$, $r_0$ and $m$ are as in Lemma~\ref{up:lem:uniform}. All constants below are uniform over $p,f\in\cF$. Since $K$ is closed, Sobolev norms on $K$, such as $\norm\cdot_{H^\sigma(K)}$, are understood as norms on the interior of $K$, which has smooth boundary.

\subsection{The pilot}\label{up:sec:pilot}

Given the observation \eqref{eq:model}, draw an independent white noise $\mathbb W'$ and set $Y_1=Y+N^{-1/2}\mathbb W'$, $Y_2=Y-N^{-1/2}\mathbb W'$. These are independent observations of $u_f$ with noise level $(N/2)^{-1/2}$. The auxiliary randomisation can be removed by conditioning: the conditional expectation given $Y$ of a randomised estimator has no larger absolute or squared risk, by Jensen's inequality.

Let $\mathcal E$ be a bounded linear extension operator from functions on $\cO$ to functions on $\R^d$, with values supported in a fixed compact set $\mathcal Q$, which is bounded simultaneously $L^2(\cO)\to L^2(\R^d)$ and $H^{s+1}(\cO)\to H^{s+1}(\R^d)$ \cite{stein1970singular}. Let $P_J$ be the Fourier multiplier with symbol $\varphi(\xi/J)$, where $\varphi\in C_c^\infty$ equals one near zero.

\begin{lemma}[Prediction pilot]\label{up:lem:pilot}
Let $J=N^{1/(2s+d+2)}$, $b_N=N^{-(s+1)/(2s+d+2)}$ and $\bar u=(P_J\mathcal EY_1)|_\cO$. Let $\{p_j\}_{j\ge1}$ be a countable subset of $\mathcal G$ such that $\{u_{p_j}\}$ is dense in $\{u_p:p\in\mathcal G\}$ for the $L^2$ norm, and let $p$ be the first $p_j$ with
\[
 \norm{u_{p_j}-\bar u}_2<\inf_k\norm{u_{p_k}-\bar u}_2+b_N .
\]
Then $p$ is a measurable function of $Y_1$ with values in $\mathcal G$, and for every $q\ge1$
\begin{equation}\label{up:eq:pilot}
 \sup_{f\in\mathcal G}\E_f\norm{u_p-u_f}_2^q\le C_qb_N^q .
\end{equation}
\end{lemma}

\begin{proof}
A countable family with the stated density property exists because $H^s$ is separable and $p\mapsto u_p$ is continuous from $H^s$ to $L^2$. The map $Y_1\mapsto\bar u$ is measurable, the infimum over a countable family is measurable, and $p$ takes countably many values on measurable events, so it is measurable. Since $u_f$ lies in the $L^2$-closure of $\{u_{p_k}\}$, we have $\inf_k\norm{u_{p_k}-\bar u}_2\le\norm{u_f-\bar u}_2$, and therefore
\[
 \norm{u_p-u_f}_2\le\norm{u_p-\bar u}_2+\norm{\bar u-u_f}_2\le2\norm{\bar u-u_f}_2+b_N .
\]
It remains to bound $\bar u-u_f=(P_J\mathcal Eu_f-\mathcal Eu_f)|_\cO+(2/N)^{1/2}(P_J\mathcal E\mathbb W_1)|_\cO$, where $\mathbb W_1$ is the white noise in $Y_1$. The first term is deterministic. Since $\varphi=1$ near zero,
\[
 \norm{P_J\mathcal Eu_f-\mathcal Eu_f}_{L^2}\le CJ^{-(s+1)}\norm{\mathcal Eu_f}_{H^{s+1}}\le CJ^{-(s+1)}
\]
by Lemma~\ref{up:lem:uniform}(a). The second term is a centred Gaussian element of $L^2$. Its second moment is $(2/N)$ times the squared Hilbert--Schmidt norm of $P_J\mathcal E=P_J\mathbf 1_{\mathcal Q}\mathcal E$, so it is at most $(2/N)\norm{P_J\mathbf 1_{\mathcal Q}}_{\rm HS}^2\norm{\mathcal E}^2\le CJ^d/N$. Moments of order $q$ of the norm of a Gaussian element are bounded by $C_q$ times the $q/2$-th power of its second moment. With the choice of $J$, both terms are of order $b_N$.
\end{proof}

The projection step ensures that the pilot satisfies every defining condition of $\mathcal G$, such as isotropy of the critical Hessian, without requiring this of a generic smoothed estimate. This construction yields a measurable near-minimiser with the prediction bounds of Lemma~\ref{up:lem:pilot}; a computationally efficient implementation is a direction for further work.

\subsection{The one-step estimator and the exact remainder}\label{up:sec:onestep}

Let $p$ be the pilot from Lemma~\ref{up:lem:pilot}, and let $w\in H^2(\cO)\cap H^1_0(\cO)$ be a function determined by $p$ (a regularised adjoint solution, constructed below). Put
\begin{equation}\label{up:eq:onestep}
 z=L_pw,\qquad r=\psi-\nabla u_p\cdot\nabla w,\qquad
 \widehat\Phi=\ip p\psi+Y_2(z)-\ip{u_p}z .
\end{equation}
Here $Y_2(z)=\ip{u_f}z+(2/N)^{1/2}\mathbb W_2(z)$ for the white noise $\mathbb W_2$ in $Y_2$.

\begin{lemma}[Exact conditional bias]\label{up:lem:identity}
Let $f\in\cF$, $h=f-p$ and $v=u_f-u_p$. Then
\[
 \E\bigl[\widehat\Phi-\Phi(f)\bigm|Y_1\bigr]=-\ip hr+\int_\cO h\,\nabla v\cdot\nabla w,
 \qquad
 \operatorname{Var}\bigl(\widehat\Phi\bigm|Y_1\bigr)=\frac2N\norm z_2^2 .
\]
\end{lemma}

\begin{proof}
Subtracting $\diver(p\nabla u_p)=g$ from $\diver(f\nabla u_f)=g$ gives $L_pv=-\diver(h\nabla u_f)$. Since $v,w\in H^2\cap H^1_0$ and $L_p$ is symmetric on this space,
\[
 \ip vz=\ip{L_pv}w=-\int\diver(h\nabla u_f)\,w=\int h\nabla u_f\cdot\nabla w
 =\int h\nabla u_p\cdot\nabla w+\int h\nabla v\cdot\nabla w .
\]
The first term on the right equals $\ip h{\psi-r}$. Hence $\E[\widehat\Phi\mid Y_1]=\ip p\psi+\ip vz=\ip f\psi-\ip hr+\int h\nabla v\cdot\nabla w$. The variance statement is immediate.
\end{proof}

The identity is exact; no Taylor expansion of the forward map is involved. The term $-\ip hr$ is a regularisation bias, which vanishes if $w$ is the exact adjoint solution. The second term is a quadratic remainder. We bound it by combining the following stability estimate with the weighted inequality of Lemma~\ref{up:lem:hardy}.

\begin{lemma}[Local transport stability]\label{up:lem:stability}
For $f,p\in\cF$, with $h=f-p$ and $v=u_f-u_p$,
\[
 \norm h_{L^2(K)}\le C\norm v_{H^2(K)} .
\]
\end{lemma}

\begin{proof}
Subtracting the two equations in the other order gives $\diver(h\nabla u_p)=-L_fv=:F$, i.e.\ $\nabla u_p\cdot\nabla h+(\Delta u_p)h=F$. Let $X_t$ be the backward flow of $\nabla u_p$, so $\dot X_t=-\nabla u_p(X_t)$. For $x\in K$, the function $\eta(t)=h(X_tx)$ satisfies $\eta'=(\Delta u_p)(X_tx)\,\eta-F(X_tx)$. Put $J_t(x)=\exp(-\int_0^t\Delta u_p(X_\tau x)\,d\tau)$. By Lemma~\ref{up:lem:uniform}(b),(c), $\Delta u_p\ge dm$ on $B_{r_0}(x_p)$ and the entrance time into that ball is bounded, so $J_t(x)\le Ce^{-dmt}$. The function $h$ is bounded, since $s>d/2$. Hence the bounded solution of this linear equation is
\[
 h(x)=\int_0^\infty J_t(x)\,F(X_tx)\,dt .
\]
By Liouville's formula $J_t=\det DX_t$, so the substitution $y=X_tx$, which maps $K$ injectively into $K$, gives
\[
 \int_KJ_t(x)^2F(X_tx)^2\,dx\le\sup_KJ_t\int_KF(y)^2\,dy\le Ce^{-dmt}\norm F_{L^2(K)}^2 .
\]
By Minkowski's inequality, $\norm h_{L^2(K)}\le C\norm F_{L^2(K)}$. Finally $\norm{L_fv}_{L^2(K)}\le C\norm v_{H^2(K)}$, since $f$ is uniformly bounded in $W^{1,\infty}$.
\end{proof}

We also record the interpolation bounds used for the remainder. Let $\delta_p=\norm v_2=\norm{u_f-u_p}_2$ denote the prediction error of the pilot. Since $\norm h_{H^s(K)}$ and $\norm v_{H^{s+1}}$ are uniformly bounded, interpolation between $L^2$ and $H^{s+1}$ for $v$, and between $L^2(K)$ and $H^s(K)$ for $h$, together with Lemma~\ref{up:lem:stability}, gives for $0\le\sigma\le s$ and $0\le\tau\le s+1$:
\begin{equation}\label{up:eq:interp}
 \norm h_{L^2(K)}\le C\delta_p^{\frac{s-1}{s+1}},\qquad
 \norm h_{H^\sigma(K)}\le C\delta_p^{\frac{s-1}{s+1}\left(1-\frac\sigma s\right)},\qquad
 \norm v_{H^\tau}\le C\delta_p^{1-\frac{\tau}{s+1}} .
\end{equation}
For any exponent $0<q<2$, Jensen's inequality and \eqref{up:eq:pilot} give
\begin{equation}\label{up:eq:jensen}
 \E_f\delta_p^q\le(\E_f\delta_p^2)^{q/2}\le C\,b_N^{\,q}=C\,N^{-q(s+1)/(2s+d+2)} .
\end{equation}

\subsection{Representer regularisation on isotropic classes}\label{up:sec:isoreg}

In this subsection and the next, $\mathcal G\subset\cF$ is a class on which every critical Hessian is isotropic:
\begin{equation}\label{up:eq:iso}
 D^2u_p(x_p)=\lambda_pI_d,\qquad \lambda_-\le\lambda_p\le\lambda_+,\qquad p\in\mathcal G,
\end{equation}
for fixed $0<\lambda_-\le\lambda_+$. Both the class $\cI_s(\eta)$ and the larger class in Theorem~\ref{thm:isotropic} are of this form.

\begin{lemma}[Radial coordinates]\label{up:lem:coord}
For $p\in\mathcal G$ there is a $C^2$ diffeomorphism $\Xi_p$ from a fixed neighbourhood of $x_p$ onto a neighbourhood of $0$ with $\Xi_p(x_p)=0$, $D\Xi_p(x_p)=I_d$, $C^2$ norms of $\Xi_p$ and $\Xi_p^{-1}$ bounded uniformly in $p$, and
\[
 D\Xi_p(x)\nabla u_p(x)=\lambda_p\Xi_p(x).
\]
In these coordinates
\begin{equation}\label{up:eq:decomp}
 w_p(x)=\frac{\psi(x_p)}{\lambda_p}\log|\Xi_p(x)|+C_p\Bigl(\frac{\Xi_p(x)}{|\Xi_p(x)|}\Bigr)+q_p(\Xi_p(x)),
\end{equation}
where $C_p\in C^2(S^{d-1})$ and $q_p\in C^2$ have norms bounded uniformly in $p$. Consequently $|D^kw_p(x)|\le C|x-x_p|^{-k}$ for $k=1,2$. Under \condC\ at $p$, $w_p$ extends to a $C^2$ function near $x_p$ with uniformly bounded norm.
\end{lemma}

\begin{proof}
Translate $x_p$ to $0$ and write $b(x)=\nabla u_p(x)=\lambda_px+R(x)$ with $R(x)=O(|x|^2)$, $DR(x)=O(|x|)$ and $D^2R=O(1)$, uniformly by Lemma~\ref{up:lem:uniform}(a). Let $X_t$ be the backward flow. On a small fixed ball, $|X_tx|\le Ce^{-\lambda_pt}|x|$. The variational equation $\frac d{dt}DX_t=-Db(X_t)DX_t$ with $Db(X_t)=\lambda_pI+O(|X_t|)$, and $\int_0^\infty|X_t|\,dt<C$, give $\norm{DX_t}\le Ce^{-\lambda_pt}$. The second variational equation has forcing $D^2b(X_t)[DX_t,DX_t]=O(e^{-2\lambda_pt})$, so $\norm{D^2X_t}\le Ce^{-\lambda_pt}$. Hence
\[
 \Xi_p(x)=\lim_{t\to\infty}e^{\lambda_pt}X_tx=x-\int_0^\infty e^{\lambda_pt}R(X_tx)\,dt
\]
converges in $C^2$: the integrand is $O(e^{-\lambda_pt}|x|)$, and its first two derivatives are $O(e^{-\lambda_p t})$. Moreover $D\Xi_p(0)=I_d$. The flow property $\Xi_p(X_sx)=e^{-\lambda_ps}\Xi_p(x)$, differentiated at $s=0$, gives the conjugacy. The linear part $\lambda_pI$ has no resonances, which is why the conjugacy is $C^2$ without further conditions.

Let $\tilde\psi=\psi\circ \Xi_p^{-1}$ and $q_p(y)=\lambda_p^{-1}\int_0^1(\tilde\psi(ty)-\tilde\psi(0))\,t^{-1}\,dt$. Then $q_p$ is uniformly $C^2$ and $\lambda_py\cdot\nabla q_p(y)=\tilde\psi(y)-\tilde\psi(0)$. Also $y\cdot\nabla\log|y|=1$. In the coordinates $y=\Xi_p(x)$ the transport equation reads $\lambda_py\cdot\nabla_y\tilde w=\tilde\psi$. Thus $\tilde w-\lambda_p^{-1}\tilde\psi(0)\log|y|-q_p$ is constant along rays, so it equals $C_p(y/|y|)$. On a fixed small sphere, $C_p$ is determined by $w_p$, which is obtained by integrating $\psi$ along characteristics that stay away from $x_p$ until they exit $\cO$; it is therefore uniformly $C^2$. The derivative bounds follow. Under \condC\ the logarithmic term is absent and $C_p$ is constant, and the last claim follows.
\end{proof}

The regulariser removes the singularity of $w_p$ inside $B_\rho(x_p)$ in two steps. Subtracting the logarithmic term's value at radius $\rho$ keeps the cutoff-annulus residual bounded independently of $\rho$. Small support alone does not give the required negative Sobolev estimate when $s>d/2$, so we also cancel finitely many residual moments. This yields $\norm{r_{p,\rho}}_{H^{-s}}\lesssim\rho^{s+d/2}$ while preserving the derivative bounds that determine the variance.

Let $\chi\in C^\infty([0,\infty))$ vanish on $[0,1]$ and equal one on $[2,\infty)$. Fix an integer $M\ge0$ with
\begin{equation}\label{up:eq:M}
 M+1>s-\frac d2 ,
\end{equation}
and a nonnegative $\chi_0\in C_c^\infty(B_1)$ that is positive on a ball. Let $\nu_\beta(y)=\chi_0(y)y^\beta$ for multi-indices $|\beta|\le M$. Define
\begin{equation}\label{up:eq:Vd}
 V_d(\rho)=\rho^{d-4}\ (d=2,3),\qquad V_4(\rho)=1+\log(1/\rho),\qquad V_d(\rho)=1\ (d\ge5).
\end{equation}

\begin{lemma}[Moment-cancelled regulariser]\label{up:lem:reg}
There are $\rho_0>0$ and $C$ such that for all $\rho\in(0,\rho_0]$ and $p\in\mathcal G$ there is $w_{p,\rho}\in H^2(\cO)\cap H^1_0(\cO)$, supported in $K$, such that $z_{p,\rho}=L_pw_{p,\rho}$ and $r_{p,\rho}=\psi-\nabla u_p\cdot\nabla w_{p,\rho}$ satisfy
\begin{align*}
 &\supp r_{p,\rho}\subset B_{2\rho}(x_p),\qquad \norm{r_{p,\rho}}_{H^{-s}(\R^d)}\le C\rho^{s+d/2},\\
 &|\nabla w_{p,\rho}(x)|\le\frac C{|x-x_p|+\rho},\qquad \norm{z_{p,\rho}}_2^2\le CV_d(\rho).
\end{align*}
\end{lemma}

\begin{proof}
Let $c=\psi(x_p)\lambda_p^{-1}\log\rho$ and define
\[
 w^0(x)=c+\chi(|x-x_p|/\rho)\bigl(w_p(x)-c\bigr),
\]
which equals $w_p$ outside $B_{2\rho}(x_p)$. On the annulus $\rho\le|x-x_p|\le2\rho$, the decomposition \eqref{up:eq:decomp} and $|\Xi_p(x)|\asymp|x-x_p|$ show that $|w_p-c|\le C$, so $|D^kw^0|\le C\rho^{-k}$ there, for $k=1,2$. On $B_\rho(x_p)$, $w^0$ is constant. The preliminary residual is
\[
 r^0=\psi-\nabla u_p\cdot\nabla w^0=(1-\chi)\psi-(w_p-c)\,\nabla u_p\cdot\nabla\bigl[\chi(|x-x_p|/\rho)\bigr].
\]
It is supported in $B_{2\rho}(x_p)$ and bounded by $C$, because $|\nabla u_p(x)|\le C|x-x_p|$.

Set $w_{p,\rho}=w^0+\sum_{|\beta|\le M}c_\beta\,\nu_\beta((x-x_p)/\rho)$. With $R(y)=r_{p,\rho}(x_p+\rho y)$ and $R^0(y)=r^0(x_p+\rho y)$,
\[
 R(y)=R^0(y)-\sum_\beta c_\beta\,b_\rho(y)\cdot\nabla\nu_\beta(y),\qquad b_\rho(y)=\rho^{-1}\nabla u_p(x_p+\rho y).
\]
The vanishing of $\int y^\alpha R(y)\,dy$ for all $|\alpha|\le M$ is the linear system $\sum_\beta M_{\alpha\beta}c_\beta=\int y^\alpha R^0$, where $M_{\alpha\beta}=\int y^\alpha\,b_\rho\cdot\nabla\nu_\beta$. Since $b_\rho(y)=\lambda_py+O(\rho)$ in $C^1(B_1)$ uniformly, an integration by parts gives
\[
 M_{\alpha\beta}=-\lambda_p(d+|\alpha|)\int\chi_0(y)\,y^{\alpha+\beta}\,dy+O(\rho).
\]
The matrix $\bigl(\int\chi_0y^{\alpha+\beta}\bigr)_{\alpha\beta}$ is a positive definite Gram matrix, because the monomials are linearly independent on the set where $\chi_0>0$. Multiplying it on the left by the invertible diagonal matrix $\operatorname{diag}(-\lambda_p(d+|\alpha|))$ gives an invertible matrix with inverse bounded uniformly for $\lambda_p\in[\lambda_-,\lambda_+]$. Hence for $\rho\le\rho_0$ the system has a solution with $|c_\beta|\le C$. The function $R$ is then bounded by $C$, supported in $B_2$, and orthogonal to all polynomials of degree at most $M$. Lemma~\ref{up:lem:negsob} and \eqref{up:eq:M} give $\norm{r_{p,\rho}}_{H^{-s}}\le C\rho^{s+d/2}$.

The gradient bound follows from $|\nabla w_p|\le C|x-x_p|^{-1}$ outside $B_{2\rho}(x_p)$ and $|\nabla w_{p,\rho}|\le C\rho^{-1}$ inside. Since $p$ is uniformly bounded in $C^1$,
\[
 |z_{p,\rho}|\le C\bigl(|D^2w_{p,\rho}|+|\nabla w_{p,\rho}|\bigr)\le
 \begin{cases}C\rho^{-2},&|x-x_p|\le2\rho,\\ C|x-x_p|^{-2},&|x-x_p|>2\rho,\end{cases}
\]
and $z_{p,\rho}$ is bounded away from $x_p$. Integrating the square gives $\norm{z_{p,\rho}}_2^2\le C(\rho^{d-4}+\int_{2\rho}^{1}t^{d-5}\,dt+1)\le CV_d(\rho)$. In $d=2$ this accounts for the regularisation of the point mass in the distributional Laplacian of $\log|x-x_p|$: the bound $\rho^{-2}$ comes from the inner ball and is not an artefact. Finally, $w_{p,\rho}$ is $C^2$ near $x_p$, vanishes outside a compact subset of $K$ by Lemma~\ref{up:lem:uniform}(d), and is therefore in $H^2\cap H^1_0$.
\end{proof}

For $d\ge5$ no regularisation is needed:

\begin{lemma}[Exact representer for $d\ge5$]\label{up:lem:exact-iso}
If $d\ge5$, then $w_p\in H^2(\cO)\cap H^1_0(\cO)$ with $\norm{w_p}_{H^2}\le C$ uniformly in $p\in\mathcal G$. Thus $z_p=L_pw_p$ satisfies $A_p^*z_p=\psi$ and $\norm{z_p}_2\le C$.
\end{lemma}

\begin{proof}
By Lemma~\ref{up:lem:coord}, $|w_p|\le C(1+|\log|x-x_p||)$ and $|D^kw_p|\le C|x-x_p|^{-k}$ for $k=1,2$, and $\int_{B_1}|x|^{-4}\,dx<\infty$ for $d\ge5$. Lemma~\ref{up:lem:H2} shows that the pointwise derivatives are weak derivatives.
\end{proof}

\subsection{The upper bound in Theorem~\ref{thm:isotropic}}\label{up:sec:iso}

Let
\begin{equation}\label{up:eq:gamma}
 \gamma(s,d)=\frac{(2s-1)^2}{2s(2s+d+2)} .
\end{equation}

\begin{theorem}[Upper bound on isotropic classes]\label{up:thm:iso}
Let \eqref{up:eq:s} and \eqref{up:eq:iso} hold. For $d\le4$ let $\widehat\Phi_\rho$ be the one-step estimator \eqref{up:eq:onestep} with $w=w_{p,\rho}$ from Lemma~\ref{up:lem:reg}; for $d\ge5$ let it use $w=w_p$ from Lemma~\ref{up:lem:exact-iso}, and set $\rho=0$ in the bias term below. Then, uniformly in $f\in\mathcal G$,
\begin{equation}\label{up:eq:iso-risk}
 \E_f|\widehat\Phi_\rho-\Phi(f)|\le C\Bigl\{\rho^{s+d/2}+N^{-1/2}V_d(\rho)^{1/2}+N^{-\gamma(s,d)}\Bigr\}.
\end{equation}
With $\rho=N^{-1/(2s+4)}$, the right-hand side is of order $\rate(s,d)$ provided
\begin{equation}\label{up:eq:gamma-cond}
 \gamma(s,d)\ge\frac{2s+d}{4s+8}\quad(d=2,3),\qquad \gamma(s,d)\ge\frac12\quad(d\ge4).
\end{equation}
Condition \eqref{up:eq:gamma-cond} together with \eqref{up:eq:s} is equivalent to the weaker sufficient thresholds in Remark~\ref{rem:thresholds}. These conditions are implied by the simpler assumption $s>d/2+3$ in Theorem~\ref{thm:isotropic}.
\end{theorem}

\begin{proof}
Condition on $Y_1$ and use Lemma~\ref{up:lem:identity}, with $h=f-p$ and $v=u_f-u_p$.

\emph{Bias.} By interpolation, the extension operator $\mathcal E$ of Section~\ref{up:sec:pilot} is also bounded from $H^s(\cO)$ to $H^s(\R^d)$. Since $r_{p,\rho}$ is supported in $\cO$, Lemma~\ref{up:lem:reg} gives $|\ip h{r_{p,\rho}}|=|\ip{\mathcal Eh}{r_{p,\rho}}|\le\norm{\mathcal Eh}_{H^s}\norm{r_{p,\rho}}_{H^{-s}}\le C\rho^{s+d/2}$. For $d\ge5$ this term vanishes.

\emph{Variance.} The conditional standard deviation is $(2/N)^{1/2}\norm{z}_2\le CN^{-1/2}V_d(\rho)^{1/2}$.

\emph{Remainder.} The function $w$ vanishes outside $K$ and $|\nabla w(x)|\le C|x-x_p|^{-1}$. By the Cauchy--Schwarz inequality and Lemma~\ref{up:lem:hardy} with $\kappa=1$, applied to $h$ and to each component of $\nabla v$,
\[
 \Bigl|\int h\nabla v\cdot\nabla w\Bigr|\le C\Bigl(\int_K\frac{|h|^2}{|x-x_p|}\Bigr)^{1/2}\Bigl(\int_K\frac{|\nabla v|^2}{|x-x_p|}\Bigr)^{1/2}\le C\norm h_{H^{1/2}(K)}\norm v_{H^{3/2}}.
\]
By \eqref{up:eq:interp} with $\sigma=1/2$ and $\tau=3/2$, the right-hand side is at most $C\delta_p^q$ with
\[
 q=\frac{s-1}{s+1}\Bigl(1-\frac1{2s}\Bigr)+\frac{s-1/2}{s+1}=\frac{(2s-1)^2}{2s(s+1)} .
\]
Since $(2s-1)^2<4s(s+1)$, we have $q<2$. By \eqref{up:eq:jensen}, the expected remainder is at most $CN^{-q(s+1)/(2s+d+2)}=CN^{-\gamma(s,d)}$.

Taking expectations over $Y_1$ proves \eqref{up:eq:iso-risk}. With $\rho=N^{-1/(2s+4)}$ and $d\le3$, both $\rho^{s+d/2}$ and $N^{-1/2}\rho^{(d-4)/2}$ equal $N^{-(2s+d)/(4s+8)}$. For $d=4$ they are $N^{-1/2}$ and $N^{-1/2}(1+\log N/(2s+4))^{1/2}$. For $d\ge5$ the bound is $C(N^{-1/2}+N^{-\gamma})$.

\emph{Smoothness conditions.} For $d=2$, \eqref{up:eq:gamma-cond} reads $(2s-1)^2\ge2s(s+1)$, i.e.\ $2s^2-6s+1\ge0$, i.e.\ $s\ge(3+\sqrt7)/2\approx2.82$; together with \eqref{up:eq:s} this is $s>3$. For $d=3$ it reads $(2s-1)^2(4s+8)\ge2s(2s+5)(2s+3)$, i.e.\ $4s^3-8s^2-29s+4\ge0$. This cubic has exactly one root above $3.5$, namely $s_3=3.8264216\ldots$, and \eqref{up:eq:s} is implied. For $d\ge4$, $\gamma\ge1/2$ reads $(2s-1)^2\ge s(2s+d+2)$, i.e.\ $2s^2-(d+6)s+1\ge0$. For $d=4$ this means $s\ge(5+\sqrt{23})/2\approx4.898>4$. At $s=d/2+3$ the quadratic equals $1$, and it is increasing for $s\ge(d+6)/4$; hence $s>d/2+3$ is sufficient in every dimension, including $d=2,3$, where $\gamma\ge1/2$ implies \eqref{up:eq:gamma-cond}.
\end{proof}

\begin{proof}[Proof of Theorem~\ref{thm:isotropic}, upper bound]
Apply Theorem~\ref{up:thm:iso} with $\mathcal G=\cI_s(\eta)$, or with the larger isotropic class, which satisfies \eqref{up:eq:iso} for small $\eta$ since $\lambda_p=g(x_p)/(d\,p(x_p))$. The lower bound is proved in Section~\ref{sec:lower-proofs}.
\end{proof}

\begin{remark}[Squared risk and location of the critical point]\label{up:rem:iso}
(i) All moments in \eqref{up:eq:pilot} are controlled, so the same argument bounds the root mean squared error by the right-hand side of \eqref{up:eq:iso-risk}.
(ii) The proof never compares $x_p$ with $x_f$. The residual moments are cancelled about the pilot's own critical point, and the discrepancy between the characteristics of $u_p$ and $u_f$ is contained in the exact remainder. For completeness, $|x_p-x_f|\le C\norm{u_p-u_f}_{C^1}$ by nondegeneracy, and interpolating $C^1$ through $H^{d/2+1+\delta}$ gives $(\E_f|x_p-x_f|^q)^{1/q}\le C_{q,\delta}N^{-(s-d/2-\delta)/(2s+d+2)}$ for every $\delta>0$.
\end{remark}

\subsection{Classes on which condition (C) holds at every member}\label{up:sec:uniformC}

Under \condC\ at every member, the exact adjoint solution is uniformly regular, and a uniform $O(N^{-1/2})$ risk bound is attainable. Isotropy is not needed here.

\begin{proposition}\label{up:prop:uniformC}
Let \eqref{up:eq:s} hold, $\eta\le\eta_0$, and let $\mathcal G\subset\cF$ be a class such that \condC\ holds at every $p\in\mathcal G$. Then the one-step estimator with $w=w_p$ satisfies $\sup_{f\in\mathcal G}\E_f|\widehat\Phi-\Phi(f)|\le CN^{-1/2}$.
\end{proposition}

\begin{proof}
Let $X_t$ be the backward flow of $\nabla u_p$, which is defined on $\overline\cO$ because $u_p$ extends to a $C^3$ function on a neighbourhood of $\overline\cO$. Put
\[
 V_p(x)=\int_0^\infty\psi(X_tx)\,dt .
\]
Here $\psi(x_p)=0$, and $|X_tx-x_p|$, $\norm{DX_t}$ and $\norm{D^2X_t}$ decay like $e^{-mt}$ after the uniformly bounded entrance time into $B_{r_0}(x_p)$. This follows from the variational equations as in the proof of Lemma~\ref{up:lem:coord}, using only $D^2u_p\ge mI$ there. Differentiating under the integral shows that $V_p\in C^2(\overline\cO)$ with norm bounded uniformly in $p$. Since $\frac d{d\tau}V_p(X_\tau x)=-\psi(X_\tau x)$ and $\frac{d}{d\tau}X_\tau=-\nabla u_p(X_\tau)$, we have $\nabla u_p\cdot\nabla V_p=\psi$. For $y\in\partial\cO$ on the flow line $\ell$, $V_p(y)=\int_{-\infty}^{t_\partial(\ell)}\psi(\gamma_\ell)=-\Psi(\ell)$. Condition \condC\ at $p$ means that $V_p$ is constant, say equal to $c_p$, on $\partial\cO$. Then $w_p=V_p-c_p$ solves \eqref{eq:transport} for $p$, is in $H^2\cap H^1_0$, and has $\norm{w_p}_{C^2}\le C$.

With $w=w_p$ there is no bias term, and the conditional variance is $(2/N)\norm{L_pw_p}_2^2\le C/N$. By Lemma~\ref{up:lem:stability} and \eqref{up:eq:interp},
\[
 \Bigl|\int h\nabla v\cdot\nabla w_p\Bigr|\le C\norm h_{L^2(K)}\norm v_{H^1}\le C\delta_p^{\frac{s-1}{s+1}+\frac s{s+1}}=C\delta_p^{\frac{2s-1}{s+1}},
\]
whose expectation is $O(N^{-(2s-1)/(2s+d+2)})=O(N^{-1/2})$ because $s\ge d/2+2$.
\end{proof}

By Proposition~\ref{prop:instability}(b), such classes are strictly smaller than Sobolev neighbourhoods in $d\le4$. Imposing \condC\ at the centre $f_0$ alone does not give the conclusion of Proposition~\ref{up:prop:uniformC}.

\subsection{Proof of Theorem~\ref{thm:full-regular}}\label{up:sec:full-regular}

The lower bound $N^{-1/2}$ is Theorem~\ref{thm:lower}(c). For the upper bound we use the one-step estimator with the exact adjoint solution, and the pilot projected onto the full class $\mathcal G=\cF$. Lemmas~\ref{up:lem:pilot}, \ref{up:lem:identity} and~\ref{up:lem:stability} did not use isotropy.

\emph{General base in the regular regime ($\Teff_0>4$).} Write $H_p=D^2u_p(x_p)$, $Q_p=\tr H_p$ and $\lambda_p^*=\lmax(H_p)$. Since $H_p$ depends continuously on $p$ in the $C^3$ topology, Lemma~\ref{up:lem:uniform}(a) gives, for $\eta$ small, a fixed $\delta>0$ with $Q_p-4\lambda_p^*\ge\delta$ for all $p\in\cF$. Lemma~\ref{up:lem:flowpolar} shows that $w_p\in H^2\cap H^1_0$ with $\norm{w_p}_{H^2}\le C$ uniformly. Hence $z_p=L_pw_p$ is an exact representer with $\norm{z_p}_2\le C$, the bias vanishes, and the conditional standard deviation is $O(N^{-1/2})$.

Since $\Teff_0\le d$, we have $d\ge5$. Put $t=d/2-1$. By Sobolev embedding, $\norm{\nabla w_p}_{L^{2d/(d-2)}}\le C\norm{w_p}_{H^2}\le C$ and $\norm h_{L^d(K)}\le C\norm h_{H^t(K)}$. By H\"older's inequality with exponents $d$, $2$ and $2d/(d-2)$, and by \eqref{up:eq:interp},
\begin{align*}
 \Bigl|\int h\nabla v\cdot\nabla w_p\Bigr|&\le\norm h_{L^d(K)}\norm{\nabla v}_2\norm{\nabla w_p}_{L^{2d/(d-2)}}\le C\delta_p^{q},\\
 q&=\frac{s-1}{s+1}\Bigl(1-\frac ts\Bigr)+\frac s{s+1}=\frac{2s-d/2+(d/2-1)/s}{s+1}.
\end{align*}
Here $q<2$, so by \eqref{up:eq:jensen} the expected remainder is $O(N^{-\gamma_{\rm gen}})$ with
\[
 \gamma_{\rm gen}=\frac{2s-d/2+(d/2-1)/s}{2s+d+2} .
\]
The condition $\gamma_{\rm gen}\ge1/2$ is $4s-d+(d-2)/s\ge2s+d+2$, i.e.\ $2s^2-(2d+2)s+(d-2)\ge0$, i.e.\ $s^2-(d+1)s+(d-2)/2\ge0$. If $s>d+1$ this holds, since $s(s-d-1)>0$ and $d\ge5$.

\emph{Isotropic base.} If $H$ is isotropic, $d=\Teff_0\ge5$. Fix $\kappa>1$ with
\begin{equation}\label{up:eq:kappa}
 4\kappa<d,\qquad \kappa\Bigl(1-\frac1{2s}\Bigr)<s-\frac d2-2 ;
\end{equation}
this is possible when $s>d/2+3$, since then the right-hand side of the second condition exceeds $1>1-1/(2s)$. By continuity of $D^2u_0$ at $x_*$ and Lemma~\ref{up:lem:uniform}(a), there are a radius $r_*$ and $\eta>0$ such that $m'I\le D^2u_p\le\kappa m'I$ on $B_{r_*}(x_p)$ for all $p\in\cF$. Lemma~\ref{up:lem:aniso-deriv} gives $|D^kw_p|\le C|x-x_p|^{-k\kappa}$ and $|w_p|\le C(1+|\log|x-x_p||)$. Since $4\kappa<d$, Lemma~\ref{up:lem:H2} shows $w_p\in H^2\cap H^1_0$ with uniformly bounded norm. The remainder is bounded with Lemma~\ref{up:lem:hardy} with exponent $\kappa$:
\[
 \Bigl|\int h\nabla v\cdot\nabla w_p\Bigr|\le C\norm h_{H^{\kappa/2}(K)}\norm v_{H^{1+\kappa/2}}\le C\delta_p^{q_\kappa},\qquad q_\kappa=\frac{2s-1-\kappa+\kappa/(2s)}{s+1},
\]
by \eqref{up:eq:interp} with $\sigma=\kappa/2$ and $\tau=1+\kappa/2$:
\[
 \frac{s-1}{s+1}\Bigl(1-\frac\kappa{2s}\Bigr)+\frac{s-\kappa/2}{s+1}=\frac{2s-1-\kappa+\kappa/(2s)}{s+1}.
\]
As $q_\kappa<2$, the expected remainder is $O(N^{-\gamma_\kappa})$ with $\gamma_\kappa=(2s-1-\kappa+\kappa/(2s))/(2s+d+2)$. Finally, $\gamma_\kappa\ge1/2$ is equivalent to $2s-d-4\ge2\kappa(1-1/(2s))$, which is the second condition in \eqref{up:eq:kappa}. \qed

\subsection{Proof of Theorem~\ref{thm:full-bracket}}\label{up:sec:bracket}

Let $\mathcal H$ be as in the theorem, with $\eta\le\eta_0$, and let the pilot be projected onto $\mathcal H$. Put $\ell_\rho=1+\log(1/\rho)$. By Lemma~\ref{up:lem:aniso-deriv}, uniformly in $p\in\mathcal H$,
\begin{equation}\label{up:eq:aniso-w}
 |w_p(x)|\le C\bigl(1+|\log|x-x_p||\bigr),\qquad |D^kw_p(x)|\le C|x-x_p|^{-k\kappa},\quad k=1,2 .
\end{equation}

\emph{Regularisation.} Let $w^0=\chi(|x-x_p|/\rho)\,w_p$. By the transport equation,
\[
 r^0=\psi-\nabla u_p\cdot\nabla w^0=(1-\chi)\psi-w_p\,\nabla u_p\cdot\nabla[\chi(|x-x_p|/\rho)],
\]
which is supported in $B_{2\rho}(x_p)$ and bounded by $C\ell_\rho$. Cancel the moments of $r^0$ through degree $M$ as in Lemma~\ref{up:lem:reg}, working in orthonormal coordinates that diagonalise $H_p=D^2u_p(x_p)$. Now the limiting moment matrix is
\[
 M_{\alpha\beta}=-\Bigl(\tr H_p+\sum_i\alpha_i\lambda_i(H_p)\Bigr)\int\chi_0(y)\,y^{\alpha+\beta}\,dy+O(\rho),
\]
since $\int y^\alpha(H_py)\cdot\nabla\nu_\beta=-\int\nu_\beta\diver(y^\alpha H_py)$. The diagonal factor is bounded below by $dm$, so the inverse is uniformly bounded, and the coefficients satisfy $|c_\beta|\le C\ell_\rho$. Call the result $w_{p,\rho}$, and let $z_{p,\rho}=L_pw_{p,\rho}$, $r_{p,\rho}=\psi-\nabla u_p\cdot\nabla w_{p,\rho}$. Lemma~\ref{up:lem:negsob}, applied to the rescaled residual of size $C\ell_\rho$, gives $\norm{r_{p,\rho}}_{H^{-s}}\le C\rho^{s+d/2}\ell_\rho$. Since $\kappa>1$, $\rho^{-1}\ell_\rho\le C_\kappa\rho^{-\kappa}$. On $B_{2\rho}(x_p)$ this gives $|\nabla w_{p,\rho}|\le C_\kappa\rho^{-\kappa}$ and $|D^2w_{p,\rho}|\le C_\kappa\rho^{-2\kappa}$. Together with \eqref{up:eq:aniso-w} outside,
\[
 |\nabla w_{p,\rho}(x)|\le C(|x-x_p|+\rho)^{-\kappa},\qquad
 \norm{z_{p,\rho}}_2^2\le C\Bigl(\rho^{d-4\kappa}+\int_{\rho}^{1}t^{d-1-4\kappa}\,dt+1\Bigr)\le C\rho^{d-4\kappa},
\]
because $d\le4<4\kappa$.

\emph{Risk bound.} By Lemma~\ref{up:lem:identity}, the bias contributes $C\rho^{s+d/2}\ell_\rho$ and the noise $CN^{-1/2}\rho^{d/2-2\kappa}$. The remainder is bounded exactly as in the isotropic case of Section~\ref{up:sec:full-regular}: Lemma~\ref{up:lem:hardy} with exponent $\kappa<d$ gives $C\delta_p^{q_\kappa}$, hence $O(N^{-\gamma_\kappa})$. Thus
\[
 \sup_{f\in\mathcal H}\E_f|\widehat\Phi-\Phi(f)|\le C\Bigl\{\rho^{s+d/2}\ell_\rho+N^{-1/2}\rho^{d/2-2\kappa}+N^{-\gamma_\kappa}\Bigr\}.
\]
Choose $\rho=(N^{-1/2}/\log(eN))^{1/(s+2\kappa)}$, so that $\ell_\rho\le C\log(eN)$ and the first two terms are of order
\[
 \Bigl(\frac{N^{-1/2}}{\log(eN)}\Bigr)^{\frac{s+d/2}{s+2\kappa}}\log(eN)=N^{-\beta_\kappa}\bigl(\log(eN)\bigr)^{\vartheta_\kappa},
 \qquad \beta_\kappa=\frac{s+d/2}{2(s+2\kappa)},\quad \vartheta_\kappa=\frac{2\kappa-d/2}{s+2\kappa}.
\]
For the finitely many $N$ for which $2\rho>r_*$, use the estimator $\Phi(f_0)$ and enlarge $C$. \qed

\begin{remark}
If $\psi(x_p)=0$ for every $p\in\mathcal H$, then $\norm{w_p}_\infty\le C$, since $|\psi(x)|\le C|x-x_p|$ and the local characteristics leave $B_{r_*}(x_p)$ at an exponential rate. The logarithmic factors then disappear, and $\rho=N^{-1/(2(s+2\kappa))}$ gives $C(N^{-\beta_\kappa}+N^{-\gamma_\kappa})$.
\end{remark}

\subsection{Proof of Corollary~\ref{cor:local-exponent}}\label{up:sec:local}

Let $d\in\{2,3,4\}$ and let $s$ satisfy \eqref{up:eq:s} and $\gamma_1>(s+d/2)/(2s+4)$, where $\gamma_1=\gamma_\kappa|_{\kappa=1}=\gamma(s,d)$. Since $\cF_s(1,\eta)$ increases with $\eta$, the exponent $e(\cF_s(1,\eta))$ is nonincreasing in $\eta$, and the limit as $\eta\downarrow0$ exists.

\emph{Upper bound on the exponent.} Since $\cI_s(\eta)\subset\cF_s(1,\eta)$, Theorem~\ref{thm:isotropic} (lower bound, valid for $s>d/2$) gives $\Risk(\cF_s(1,\eta))\gtrsim\rate(s,d)$. Hence $e(\cF_s(1,\eta))\le(2s+d)/(4s+8)$ for $d=2,3$ and $e(\cF_s(1,\eta))\le1/2$ for $d=4$, for every $\eta>0$.

\emph{Lower bound on the exponent.} Let $\kappa\in(1,d)$. The Hessian of $u_0=(|x|^2-1)/(2d)$ is $I_d/d$. By Lemma~\ref{up:lem:uniform}(a), there are $r_*>0$ and $\eta(\kappa)>0$ such that $mI\le D^2u_f\le\kappa mI$ on $B_{r_*}(x_f)$ for all $f\in\cF_s(1,\eta(\kappa))$, with $m=1/(d\sqrt\kappa)$. Theorem~\ref{thm:full-bracket} with $\mathcal H=\cF_s(1,\eta(\kappa))$ gives
\[
 e\bigl(\cF_s(1,\eta(\kappa))\bigr)\ge\min\{\beta_\kappa,\gamma_\kappa\}.
\]
As $\kappa\downarrow1$, $\beta_\kappa\to(s+d/2)/(2s+4)$ and $\gamma_\kappa\to\gamma_1>(s+d/2)/(2s+4)$. Hence
\[
 \lim_{\eta\downarrow0}e(\cF_s(1,\eta))\ge\frac{s+d/2}{2s+4}.
\]
The right-hand side equals $(2s+d)/(4s+8)$, which is $1/2$ for $d=4$. The constants in Theorem~\ref{thm:full-bracket} depend on $\kappa$, but the exponent statement is for each fixed $\kappa$ and does not require uniformity.

\emph{Strict inequality at fixed $\eta$.} Let $d\in\{2,3\}$, $\psi(0)\ne0$ and $\eta>0$. By Proposition~\ref{prop:instability}(a), for $\varepsilon>0$ small there is $f_\varepsilon$ with $\norm{f_\varepsilon-1}_{H^s}\le\eta/2$. Its potential has a unique nondegenerate minimum at $0$, with $\Teff_\varepsilon=d/(1+2d\varepsilon)<d$, and $f_\varepsilon$ is smooth and positive, so Assumption~\ref{ass:base} holds at $f_\varepsilon$. Since $\psi(0)\ne0$, \condC\ fails at $f_\varepsilon$, and $\Teff_\varepsilon<4$. Theorem~\ref{thm:lower}(b) at the base $f_\varepsilon$, with a radius $\eta'\le\eta/2$, gives $\Risk(\cF_s(f_\varepsilon,\eta'))\gtrsim N^{-(2s+\Teff_\varepsilon)/(4s+8)}$. Since $\cF_s(f_\varepsilon,\eta')\subset\cF_s(1,\eta)$, we obtain $e(\cF_s(1,\eta))\le(2s+\Teff_\varepsilon)/(4s+8)<(2s+d)/(4s+8)$. \qed

\section{Flow estimates and the normal form}\label{lin:app}

\begin{proof}[Proof of Lemma~\ref{lin:lem:flow}]
Since $\partial_\nu u_0>0$ on $\partial\cO$, the field $-b$ points strictly inward on the boundary, so $X_t(\overline\cO)\subset\overline\cO$ for $t\ge0$. Along $X_tx$ we have $\frac d{dt}u_0(X_tx)=-|b(X_tx)|^2$.

Let $U=\{u_0<u_0(x_*)+\delta\}$ for $\delta>0$ small. By nondegeneracy, $U$ is a small convex neighbourhood of $x_*$ on which $Db\ge mI$ for some $m>0$. Outside $U$, $|b|\ge c_\delta>0$, because $x_*$ is the only critical point. Hence every backward trajectory enters $U$ within the time $\max(-u_0)/c_\delta^2$ and stays there. Inside $U$, $\frac d{dt}|X_tx-x_*|^2\le-2m|X_tx-x_*|^2$, which gives exponential convergence.

The variational equation $\frac d{dt}DX_t=-Db(X_t)DX_t$ has a contracting coefficient after the bounded entrance time. So $\norm{DX_t}\le Ce^{-\mu t}$ for any $\mu<m$. The $k$-th variational equation is linear in $D^kX_t$, with the same coefficient and a forcing term that is a sum of products of at least two lower-order derivatives. Induction gives $\norm{D^kX_t}\le C_ke^{-\mu t}$. This proves (i).

For (ii), $u_0(\phi_tx)$ increases at rate $|b|^2\ge c_\delta^2$ once the trajectory has left $U$. Since $u_0\le0$, it reaches $\partial\cO$ in finite time, transversally because $\partial_\nu u_0>0$. Smoothness of $\tau$ and $\vartheta$ follows from the implicit function theorem applied to $u_0(\phi_{\tau}x)=0$. Since the time spent in $U$ tends to infinity as $x\to x_*$, so does $\tau$.

For (iii), Liouville's formula gives $\frac d{dt}\log\det DX_t=-\Delta u_0(X_t)$. Since $\Delta u_0(x_*)=\tr H>0$, we have $\Delta u_0\ge c>0$ on $U$ after shrinking $\delta$, which yields $J_t\le Ce^{-ct}$. Differentiating $\log J_t=-\int_0^t\Delta u_0(X_r)\,dr$ in $x$ gives integrands that contain at least one factor $D^jX_r$, $j\ge1$. These decay exponentially by (i), so all derivatives of $\log J_t$ are bounded uniformly in $t$, and $|D^kJ_t|\le C_kJ_t$.
\end{proof}

\begin{proof}[Proof of Lemma~\ref{lin:lem:S}]
Let $k_t(x)=J_t(x)q(X_tx)$. The flow property gives $J_{t+r}(x)=J_r(x)J_t(X_rx)$ and $X_{t+r}=X_t\circ X_r$. Differentiating $k_{t+r}(x)=J_r(x)k_t(X_rx)$ at $r=0$ yields
\[
 \partial_tk_t=-\Delta u\,k_t-\nabla u\cdot\nabla k_t=-T_uk_t
\]
in the sense of distributions. Integrating in $t$ from $0$ to $\infty$ and using $k_t\to0$ in $L^2$ gives $T_uS_uq=q$. The support statement holds because $k_t(x)\ne0$ requires $X_tx\in\supp q$.

For the Sobolev bound, let $m\ge0$ be an integer. By the chain rule, $D^m k_t$ is a sum of terms $(D^\alpha J_t)\,(D^\gamma q)(X_t)\,\prod D^{j_i}X_t$ with $|\gamma|\le m$. By Lemma~\ref{lin:lem:flow}, these are bounded by $CJ_t\sum_{|\gamma|\le m}|(D^\gamma q)(X_t)|$. Changing variables $y=X_tx$, with $dy=J_t(x)\,dx$, gives
\[
 \int_\cO J_t(x)^2|(D^\gamma q)(X_tx)|^2\,dx\le\sup J_t\int_\cO|D^\gamma q(y)|^2\,dy\le Ce^{-ct}\norm q_{H^m}^2 .
\]
Hence $\norm{k_t}_{H^m}\le Ce^{-ct/2}\norm q_{H^m}$, and by Minkowski's inequality $\norm{S_uq}_{H^m}\le C\norm q_{H^m}$. Interpolation gives the bound for real $s\ge0$.
\end{proof}

\begin{proof}[Proof of Lemma~\ref{lin:lem:normalform}]
\emph{Step 1: finite-order normal form.} Fix $M\ge2$, to be chosen large in Step~2, and write the Taylor expansion of $b$ at $0$ to order $M$ with remainder, $b(x)=Hx+\sum_{k=2}^Mb_k(x)+O(|x|^{M+1})$, where $b_k$ is homogeneous of degree $k$. A polynomial change of variables $x=y+\pi_k(y)$, with $\pi_k$ homogeneous of degree $k$, changes the degree-$k$ part of the field by $[Hy,\pi_k]$ and leaves lower degrees unchanged. On the monomial $y^\alpha e_i$, this homological operator acts as multiplication by $\alpha\cdot\lambda-\lambda_i$. For $k=2,\dots,M$ in turn, we remove all nonresonant monomials of degree $k$ and keep the resonant ones. This is a finite composition of polynomial local diffeomorphisms tangent to the identity. Since resonances have degree at most $\lambda_*/\lambda_1$, choosing $M$ larger than this bound transforms $b$ into $\tilde b=\mathcal N+\mathcal R$, where $\mathcal N=Hy+P$ is as stated and $\mathcal R$ vanishes to order $M+1$ at $0$. The structure of resonances gives $\diver P=0$, since $P_i$ does not depend on $y_i$, and $[Hy,P]=0$.

\emph{Step 2: removing the high-order remainder.} Let $\phi^{\mathcal N}_t$ be the flow of $\mathcal N$; it has the form $\phi^{\mathcal N}_t(y)=e^{tH}G_t(y)$ with $G_t$ polynomial in $t$ and $y$ (see Section~\ref{sec:normalform}), and its derivatives of each fixed order on bounded sets are bounded by $C_k(1+t)^{L_k}e^{\lambda_*t}$ for $t\ge0$. On a small ball the symmetric part of $D\tilde b$ is positive definite. Let $X_t$ be the backward flow of $\tilde b$, and fix $0<m$ below its contraction rate, so that $|X_ty|\le e^{-mt}|y|$ and $\norm{D^kX_t}\le C_ke^{-mt}$ for $k\ge1$, as in the proof of Lemma~\ref{lin:lem:flow}. Choose $M$ so large that $(M+1)m>\lambda_*$. Since $\mathcal R$ vanishes to order $M+1$, the chain rule gives $\norm{\mathcal R\circ X_t}_{C^k}\le C_ke^{-(M+1)mt}$ for every fixed $k$. Indeed, $D^k(\mathcal R\circ X_t)$ is a sum of terms $D^j\mathcal R(X_t)[D^{\alpha_1}X_t,\dots,D^{\alpha_j}X_t]$ with $1\le j\le k$ and $\alpha_1+\dots+\alpha_j=k$. Each such term is bounded by $C|X_t|^{\max(M+1-j,0)}e^{-jmt}\le Ce^{-(M+1)mt}$, because each of the $j$ factors $D^{\alpha_i}X_t$ contributes $e^{-mt}$. Thus a single $M$ serves for all $k$. Put $K_t=\phi^{\mathcal N}_t\circ X_t$. Then
\[
 \frac d{dt}K_t(y)=-D\phi^{\mathcal N}_t(X_ty)\,\mathcal R(X_ty).
\]
The right-hand side has $C^k$ norm at most a polynomial in $t$ times $e^{(\lambda_*-(M+1)m)t}$, which is integrable. So $K_t$ converges in every $C^k$ norm on a small ball to a smooth map $K$ with $DK(0)=I$. Let $\phi^{\tilde b}_s=X_{-s}$ be the forward flow of $\tilde b$. Then $K_t\circ\phi^{\tilde b}_s=\phi^{\mathcal N}_t\circ X_{t-s}=\phi^{\mathcal N}_s\circ K_{t-s}$ for small $s\ge0$ and large $t$. Letting $t\to\infty$ gives $K\circ\phi^{\tilde b}_s=\phi^{\mathcal N}_s\circ K$, so $K$ conjugates $\tilde b$ to $\mathcal N$. Composing with Step~1 gives $\mathcal K$.
\end{proof}

\section{Auxiliary results for the lower bounds}\label{app:lower}

\subsection{The transport inverse}

Let $u^\natural\in C^\infty(\overline\cO)$ have exactly one critical point $x^\natural\in\cO$, a nondegenerate minimum, and let $\nabla u^\natural$ point strictly outward on $\partial\cO$. For a potential $u$ we use the backward flow $X_t$ of $\nabla u$, its Jacobian $J_t$, $T_u$ and the transport inverse $S_u$ as in \eqref{lin:eq:S}, with $u$ in place of $u_0$.

This is a version of Lemma~\ref{lin:lem:S} that is uniform over a neighbourhood of potentials, as needed for the nonlinear constructions.

\begin{lemma}\label{low:S}
For every integer $k\ge0$ there are a $C^{k+2}$-neighbourhood $\mathcal U_k$ of $u^\natural$ and constants $\mu,C_k>0$ such that, for every $u\in\mathcal U_k$ (and with $x_u$ its unique critical point):
\begin{enumerate}
\item[(i)] $X_t(\overline\cO)\subset\overline\cO$ for $t\ge0$, $|X_tx-x_u|\le C_0e^{-\mu t}$, $|D^jX_t(x)|\le C_ke^{-\mu t}$ for $1\le j\le k+1$, $J_t\le C_0e^{-\mu t}$ and $|D^jJ_t|\le C_kJ_t$ for $j\le k$;
\item[(ii)] for $q\in C^k(\overline\cO)$, $S_uq\in C^k(\overline\cO)$ with $\norm{S_uq}_{C^k}\le C_k\norm q_{C^k}$, and $T_uS_uq=q$; every bounded solution $h$ of $T_uh=q$ equals $S_uq$;
\item[(iii)] for $q\in H^k(\cO)$, $\norm{S_uq}_{H^k}\le C_k\norm q_{H^k}$.
\end{enumerate}
By interpolation, $S_u$ is bounded on $H^s(\cO)$ for every $s\ge0$, uniformly over $u\in\mathcal U_{\lceil s\rceil}$.
\end{lemma}

\begin{proof}
(i) If $u$ is $C^2$-close to $u^\natural$, then $-\nabla u$ points strictly inward on $\partial\cO$, so $\overline\cO$ is forward invariant, and $u$ has a unique critical point $x_u$, close to $x^\natural$, with $D^2u\ge2mI$ on a ball $B$ around it. Choose a small sublevel set $W=\{u<u(x_u)+\epsilon_W\}\subset B$; it is convex and forward invariant because $\frac d{dt}u(X_t)=-|\nabla u(X_t)|^2\le0$. Outside $W$, $|\nabla u|\ge c_0>0$, so every orbit enters $W$ before the time $T_0=\operatorname{osc}(u)/c_0^2$. In $W$, $\frac d{dt}|X_t-x_u|^2=-2(X_t-x_u)\cdot\nabla u(X_t)\le-2m|X_t-x_u|^2$. The variational equation $\frac d{dt}DX_t=-D^2u(X_t)DX_t$ gives $|DX_t|\le e^{\norm u_{C^2}T_0}$ for $t\le T_0$ and contraction at rate $m$ afterwards. For $j\ge2$, $\frac d{dt}D^jX_t=-D^2u(X_t)D^jX_t+R_j$, where $R_j$ is a sum of products of derivatives $D^iu(X_t)$, $3\le i\le j+1$, with at least two factors $D^{i'}X_t$ of order $i'<j$. By induction $R_j=O(e^{-2mt})$, and Duhamel's formula gives $D^jX_t=O(e^{-mt})$. Since $\Delta u\ge2dm$ on $W$, $J_t\le Ce^{-2dmt}$. Finally $D\log J_t=-\int_0^tD(\Delta u)(X_r)DX_r\,dr$ is bounded, and the same holds for higher derivatives, so $|D^jJ_t|\le C_jJ_t$. Take $\mu=m$. All constants depend only on $\norm u_{C^{k+2}}$, $m$, $c_0$ and $T_0$, which are uniform on a small neighbourhood.

(ii) The bounds follow from (i) by differentiating under the integral. With $G=S_uq$, the cocycle identity $J_t(X_sx)J_s(x)=J_{t+s}(x)$ gives $J_s(x)G(X_sx)=\int_s^\infty J_\tau(x)q(X_\tau x)\,d\tau$. Differentiating at $s=0$ gives $-\Delta u\,G-\nabla u\cdot\nabla G=-q$, that is, $T_uG=q$. If $h$ is bounded and $T_uh=0$, then $s\mapsto J_s(x)h(X_sx)$ is constant, and letting $s\to\infty$ gives $h(x)=0$.

(iii) For $|\alpha|\le k$, $D^\alpha[J_t\,q\circ X_t]$ is a sum of terms $(D^{\alpha'}J_t)\,(D^\beta q)\circ X_t\prod D^{\gamma_i}X_t$ with $|\beta|\le|\alpha|$. By (i) each is bounded by $CJ_t|(D^\beta q)\circ X_t|$. By Liouville's formula $\det DX_t=J_t$, so the substitution $y=X_tx$ gives
\[
 \int_\cO J_t(x)^2|D^\beta q(X_tx)|^2\,dx=\int_{X_t(\cO)}J_t(X_t^{-1}y)\,|D^\beta q(y)|^2\,dy\le Ce^{-\mu t}\norm q_{H^k}^2 .
\]
Hence $\norm{J_t\,q\circ X_t}_{H^k}\le Ce^{-\mu t/2}\norm q_{H^k}$, and integration in $t$ proves (iii) for smooth $q$; density gives the general case.
\end{proof}

\subsection{Proofs of Lemmas~\ref{low:BH} and~\ref{low:potential}}

\begin{proof}[Proof of Lemma~\ref{low:BH}]
By translation and scaling, $r=1$. Let $\bar g$ denote the mean over $B_1$, $A=\overline{\nabla v}$ and $g(x)=v(x)-A\cdot x$. Then $\nabla g=\nabla v-A$ has mean zero, so by the Poincar\'e inequality $\norm{\nabla g}_{L^2(B_1)}\le C\norm{D^2v}_{L^2(B_1)}$, and again $\norm{g-\bar g}_{L^2(B_1)}\le C\norm{\nabla g}_{L^2(B_1)}$. Take $\ell(x)=A\cdot x+\bar g$.
\end{proof}

\begin{proof}[Proof of Lemma~\ref{low:potential}]
For $\delta_m$ small, $D^2v\ge I/(2d)$, so $v$ is uniformly strictly convex and, since $\nabla v(0)=0$, has the unique critical point $0$; $\nabla v=x/d$ near $\partial\mathcal B$ points outward. Thus $v$ belongs to the neighbourhood $\mathcal U_m$ of $v_0$ in Lemma~\ref{low:S}. By Lemma~\ref{low:S}(ii), $\mathsf f_p=S_v1>0$ is smooth and $T_v\mathsf f_p=1$. Since $T_v1=\Delta v=1+\Delta p$, we have $T_v(\mathsf f_p-1)=-\Delta p$, hence $\mathsf f_p-1=-S_v\Delta p$ by uniqueness, and \eqref{low:eq:pot1} follows from the $C^m$ bound. Let $h_1=-S_0\Delta p$, where $S_0=S_{v_0}$. Then $T_vh_1=T_{v_0}h_1+\diver(h_1\nabla p)=-\Delta p+\diver(h_1\nabla p)$, so $T_v(\mathsf f_p-1-h_1)=-\diver(h_1\nabla p)$ and
\[
 \norm{\mathsf f_p-1-h_1}_{C^0}=\norm{S_v\diver(h_1\nabla p)}_{C^0}\le C\norm{h_1}_{C^1}\norm p_{C^2}\le C\norm p_{C^3}\norm p_{C^2},
\]
which is \eqref{low:eq:pot2}. If $p$ is constant on a neighbourhood $U_0$ of $0$, then $D^2v(0)=I_d/d$, and $S_v\Delta p(x)=0$ whenever the backward orbit of $x$ stays in $U_0$, so $\mathsf f_p(x)=1$ there.
\end{proof}

\section{Auxiliary results for the upper bounds}\label{up:app}

The notation is that of Section~\ref{sec:upper-proofs}. In particular $\eta\le\eta_0$, $K$ is the compact set of Lemma~\ref{up:lem:uniform}, and constants are uniform over the class under consideration.

\subsection{A fractional Hardy inequality}

\begin{lemma}\label{up:lem:hardy}
Let $d\ge2$ and $0<\kappa<d$. There is $C$ such that for all $x_0\in K$ and $a\in H^{\kappa/2}(K)$,
\[
 \int_K\frac{|a(x)|^2}{|x-x_0|^\kappa}\,dx\le C\norm a_{H^{\kappa/2}(K)}^2 .
\]
In particular, for $\kappa=1$ there is no logarithmic loss in $d=2$.
\end{lemma}

\begin{proof}
On $\R^d$ the inequality $\int|x|^{-\kappa}|a|^2\le C\norm{|D|^{\kappa/2}a}_2^2$ for $0<\kappa<d$ is classical (Kato--Herbst inequality, \cite{herbst1977spectral}). We include a short proof. Write $a=c_{d,\kappa}I_{\kappa/2}g$ with $g=|D|^{\kappa/2}a$, where $I_\sigma g(x)=\int|x-y|^{-d+\sigma}g(y)\,dy$ is the Riesz potential. It suffices to show that the operator with kernel $k(x,y)=|x|^{-\kappa/2}|x-y|^{-d+\kappa/2}$ is bounded on $L^2(\R^d)$. We use Schur's test with weight $\omega(x)=|x|^{-d/2}$. The substitution $y=|x|u$ gives
\[
 \int k(x,y)\,\omega(y)\,dy=|x|^{-d/2}\int|e-u|^{-d+\kappa/2}|u|^{-d/2}\,du,
\]
where $e$ is a unit vector. The integral is finite: near $u=e$ because $\kappa>0$, near $u=0$ because $d/2<d$, and at infinity because $d-\kappa/2+d/2>d$, i.e.\ $\kappa<d$. Similarly, $x=|y|u$ gives
\[
 \int k(x,y)\,\omega(x)\,dx=|y|^{-d/2}\int|u|^{-(\kappa+d)/2}|u-e|^{-d+\kappa/2}\,du,
\]
which is finite because $(\kappa+d)/2<d$ near $u=0$ and $3d/2>d$ at infinity. Schur's test proves the bound on $\R^d$. For the local statement, extend $a$ by a bounded extension operator $H^{\kappa/2}(K)\to H^{\kappa/2}(\R^d)$, which exists since $\partial K$ is smooth, and translate $x_0$ to the origin.
\end{proof}

\subsection{Negative Sobolev norm of moment-cancelled residuals}

\begin{lemma}\label{up:lem:negsob}
Let $s>0$, $M$ an integer with $M+1>s-d/2$, and $A\ge1$. Let $R\in L^\infty(\R^d)$ be supported in $B_2$, with $\norm R_\infty\le A$, and with $\int y^\alpha R(y)\,dy=0$ for all $|\alpha|\le M$. For $x_0\in\R^d$ and $\rho\in(0,1]$ let $r(x)=R((x-x_0)/\rho)$. Then $\norm r_{H^{-s}(\R^d)}\le CA\rho^{s+d/2}$, where $C$ depends only on $d,s,M$.
\end{lemma}

\begin{proof}
$\widehat r(\xi)=\rho^de^{-ix_0\cdot\xi}\widehat R(\rho\xi)$, so the substitution $\zeta=\rho\xi$ gives
\[
 \norm r_{H^{-s}}^2=\rho^{d}\int(1+|\zeta|^2/\rho^2)^{-s}|\widehat R(\zeta)|^2\,d\zeta
 \le\rho^{d+2s}\int|\zeta|^{-2s}|\widehat R(\zeta)|^2\,d\zeta .
\]
The vanishing moments and Taylor's formula for $e^{-iy\cdot\zeta}$ give $|\widehat R(\zeta)|\le C\norm R_{L^1}|\zeta|^{M+1}\le CA|\zeta|^{M+1}$ for $|\zeta|\le1$. Since $2(M+1)-2s>-d$, the integral over $|\zeta|\le1$ is at most $CA^2$. The integral over $|\zeta|>1$ is at most $\norm R_2^2\le CA^2$.
\end{proof}

\subsection{Weak second derivatives}

\begin{lemma}\label{up:lem:H2}
Let $B$ be a ball centred at $x_0$, $d\ge2$, and let $w\in C^2(\overline B\setminus\{x_0\})$. If $w$, $\nabla w$ and $D^2w$, computed pointwise on $B\setminus\{x_0\}$, are in $L^2(B)$, then $w\in H^2(B)$ and its weak derivatives are the pointwise ones.
\end{lemma}

\begin{proof}
Fix $\eta_0\in C^\infty(\R)$ with $\eta_0=0$ on $(-\infty,1]$, $\eta_0=1$ on $[2,\infty)$, and put $\eta_\varepsilon(x)=\eta_0(|x-x_0|/\varepsilon)$, so that $|\nabla\eta_\varepsilon|\le C\varepsilon^{-1}\mathbf 1_{B_{2\varepsilon}(x_0)}$. Let $v\in C^1(\overline B\setminus\{x_0\})$ with $v,\nabla v\in L^2(B)$, and let $\varphi\in C_c^\infty(B)$. Since $\eta_\varepsilon\varphi\in C_c^1(B\setminus\{x_0\})$,
\[
 \int_B v\,\partial_i\varphi
 =-\int_B\partial_iv\,\eta_\varepsilon\varphi
 +\int_B v(1-\eta_\varepsilon)\partial_i\varphi
 -\int_B v\varphi\,\partial_i\eta_\varepsilon .
\]
As $\varepsilon\to0$, the first term tends to $-\int_B\partial_iv\,\varphi$ and the second to zero, by dominated convergence. By Cauchy--Schwarz the third is at most
\[
 C\norm\varphi_\infty\varepsilon^{-1}|B_{2\varepsilon}|^{1/2}\norm v_{L^2(B_{2\varepsilon}(x_0))}
 \le C\varepsilon^{d/2-1}\norm v_{L^2(B_{2\varepsilon}(x_0))}\to0,
\]
because $d\ge2$ and $\norm v_{L^2(B_{2\varepsilon}(x_0))}\to0$. Thus the weak derivative of $v$ is the pointwise one. Applying this with $v=w$, and then with $v=\partial_jw$, proves the claim. A point is therefore removable for $H^2$ under $L^2$ hypotheses alone, and no pointwise growth bound is required.
\end{proof}

\subsection{Derivative bounds for the adjoint solution}

For $p\in\cF$ let $w_p$ be the solution of \eqref{eq:transport} for $u_p$, as in Lemma~\ref{up:lem:uniform}(d). Away from a fixed neighbourhood of $x_p$, $w_p$ is $C^2$ with uniformly bounded norm. Indeed it is obtained by integrating $\psi$ along characteristics that reach $\partial\cO$, or a region where $w_p=0$, in uniformly bounded time, and the flow depends smoothly on the starting point.

\begin{lemma}[Controlled Hessian ratio]\label{up:lem:aniso-deriv}
Suppose that $mI\le D^2u_p\le MI$ on $B_{r_*}(x_p)$ for all $p$ in a class $\mathcal H\subset\cF$, and set $\kappa=M/m$. Then, uniformly in $p\in\mathcal H$ and $x\in K\setminus\{x_p\}$,
\[
 |w_p(x)|\le C\bigl(1+|\log|x-x_p||\bigr),\qquad |D^kw_p(x)|\le C|x-x_p|^{-k\kappa},\quad k=1,2 .
\]
\end{lemma}

\begin{proof}
Let $\rho_1=r_*/2$, $S=\partial B_{\rho_1}(x_p)$, and let $x\in B_{\rho_1}(x_p)\setminus\{x_p\}$, $r=|x-x_p|$. Let $\phi_t$ be the forward flow of $\nabla u_p$. Since $\nabla u_p(y)=\int_0^1D^2u_p(x_p+\theta(y-x_p))\,d\theta\,(y-x_p)$, we have $\frac d{dt}|\phi_tx-x_p|^2\ge2m|\phi_tx-x_p|^2$ in the ball. Hence the exit time $\tau(x)$ from $B_{\rho_1}(x_p)$ satisfies $\tau(x)\le m^{-1}\log(\rho_1/r)$. The first variational equation $\frac d{dt}D\phi_t=D^2u_p(\phi_t)D\phi_t$ and $\norm{D^2u_p}\le M$ give $\norm{D\phi_t}\le e^{Mt}$. The second variational equation has forcing $D^3u_p[D\phi_t,D\phi_t]=O(e^{2Mt})$, which gives $\norm{D^2\phi_t}\le Ce^{2Mt}$. On $S$ the vector field is transverse: $(y-x_p)\cdot\nabla u_p(y)\ge m\rho_1^2$. Implicit differentiation of $|\phi_{\tau(x)}x-x_p|^2=\rho_1^2$ therefore gives $|D\tau|\le Ce^{M\tau}$ and $|D^2\tau|\le Ce^{2M\tau}$. Since $\frac d{dt}w_p(\phi_tx)=\psi(\phi_tx)$,
\[
 w_p(x)=w_p(\phi_{\tau(x)}x)-\int_0^{\tau(x)}\psi(\phi_tx)\,dt .
\]
Differentiating once and twice, and using the uniform $C^2$ bounds for $w_p$ near $S$, gives $|w_p|\le C(1+\tau)$ and $|D^kw_p(x)|\le Ce^{kM\tau(x)}\le C(\rho_1/r)^{kM/m}$.
\end{proof}

\begin{lemma}[Effective dimension above four]\label{up:lem:flowpolar}
Suppose that $\tr H_p-4\lmax(H_p)\ge\delta>0$ for all $p$ in a class $\mathcal H\subset\cF$, where $H_p=D^2u_p(x_p)$. Then $w_p\in H^2(\cO)\cap H^1_0(\cO)$ and $\sup_{p\in\mathcal H}\norm{w_p}_{H^2}<\infty$.
\end{lemma}

\begin{proof}
Write $Q_p=\tr H_p$ and $\lambda^*_p=\lmax(H_p)$. Let $S=\partial B_{\rho_1}(x_p)$ with $\rho_1\le r_0$ small, and let $X_t$ be the backward flow. The map $(t,\theta)\mapsto X_t\theta$, for $t>0$ and $\theta\in S$, is a diffeomorphism onto $B_{\rho_1}(x_p)\setminus\{x_p\}$. By Liouville's formula its Jacobian is $(\nabla u_p\cdot n)(\theta)\,J_t(\theta)$ with $J_t(\theta)=\exp(-\int_0^t\Delta u_p(X_\sigma\theta)\,d\sigma)$. Since $|X_\sigma\theta-x_p|\le Ce^{-m\sigma}$ and $\Delta u_p(y)=Q_p+O(|y-x_p|)$, we get $J_t\asymp e^{-Q_pt}$.

Let $y=X_t\theta$. The forward orbit from $y$ reaches $S$ at time $t$, and $\norm{D^2u_p(\phi_\sigma y)-H_p}=O(e^{-m(t-\sigma)})$ is integrable along it. Grönwall's inequality therefore gives $\norm{D\phi_\sigma(y)}\le Ce^{\lambda^*_p\sigma}$ and, from the second variational equation, $\norm{D^2\phi_\sigma(y)}\le Ce^{2\lambda^*_p\sigma}$ for $0\le\sigma\le t$. As in Lemma~\ref{up:lem:aniso-deriv}, this yields
\[
 |D^kw_p(X_t\theta)|\le Ce^{k\lambda^*_pt},\quad k=1,2,\qquad |w_p(X_t\theta)|\le C(1+t).
\]
These are bounds in the flow-polar variable $t$. They should not be converted into bounds of the form $|y-x_p|^{-k}$: for anisotropic $H_p$, a point at Euclidean distance $r$ from $x_p$ can have $t$ as large as $\lmin(H_p)^{-1}\log(1/r)$, so only $|D^kw_p(y)|\le C|y-x_p|^{-k\lambda_p^*/\lmin(H_p)}$ holds pointwise. We use the $L^2$ bounds directly. In flow-polar coordinates,
\[
 \int_{B_{\rho_1}(x_p)}|D^2w_p|^2\,dx\le C\int_0^\infty e^{(4\lambda^*_p-Q_p)t}\,dt\le\frac C\delta ,
\]
and similarly
\[
 \int_{B_{\rho_1}(x_p)}(|w_p|^2+|\nabla w_p|^2)\le C\int_0^\infty\{(1+t)^2+e^{2\lambda_p^*t}\}e^{-Q_pt}\,dt\le C/\delta .
\] Lemma~\ref{up:lem:H2}, which needs only these $L^2$ bounds, shows that the pointwise derivatives are weak derivatives. Outside $B_{\rho_1}(x_p)$, $w_p$ is uniformly $C^2$ and vanishes near $\partial\cO$.
\end{proof}

\section{Details of the numerical computations}\label{num:app}

This section gives the details behind Section~\ref{sec:numerics} of the main paper. The computations are deterministic: they evaluate finite-dimensional moduli of continuity \eqref{num:eq:modulus} in the linearised experiment at a fixed reference. They do not evaluate the nonlinear one-step estimator, and the parameter sets $\Theta$ are spectral, polynomial, or finite-difference Hilbert balls rather than the ordinary Sobolev balls of Theorem~\ref{thm:linear}.

\paragraph*{Discretisations}
In each case $\Theta$ is the unit ball of a Hilbert norm.
\begin{enumerate}
\item[(i)] \emph{Balls, isotropic Hessian.} The canonical model on $B_1\subset\R^d$, $d=2,\dots,6$, so that $\Teff=d$. The problem decouples into angular sectors. We use the Dirichlet eigenbasis of each sector ($K=150$ radial modes) and the spectral class $\sum_kj_k^{2s}a_k^2\le1$.
\item[(ii)] \emph{Ellipsoids, anisotropic Hessian.} $\cO=\{\sum_ix_i^2/a_i^2<1\}$ with $f_0=g=1$, so $u_0$ is quadratic. We take $d-1$ equal semi-axes, so the eigenvalues of $H$ are proportional to $(\mu,\dots,\mu,1)$ and $\Teff=((d-1)\mu+1)/\max\{\mu,1\}$. On polynomials of degree at most $N_{\rm b}\le180$ the Galerkin forward map is exact. The class has weight $(1+n(n+d))^{s/2}$ on orthogonal polynomials of degree $n$.
\item[(iii)] \emph{A square model with zero boundary trace.} On $(-1,1)^2$ let $Ah=\Delta_D^{-1}\diver(h\,Hx)$ with $H=\operatorname{diag}(1,\kappa)$, so that $\Teff=1+1/\kappa$. Here $h=0$ on the boundary, and the norm is $\norm{\Delta h}_2$ for $s=2$ and $(\norm{\nabla\Delta h}_2^2+\norm{\Delta h}_2^2)^{1/2}$ for $s=3$. We use finite differences on $N_{\rm g}\times N_{\rm g}$ grids, $N_{\rm g}\le511$.
\end{enumerate}
The targets are of three kinds: $\psi(x_*)\ne0$; $\psi(x_*)=0$ with nonconstant flow profile; and a control for which the flow profile is constant, exactly in (i) and up to a relative variation $6\times10^{-3}$ in (ii). The two nonconstant kinds are the irregular targets. The control in (i) is regular. In (ii), its flow profile is approximately constant and it exhibits near-parametric behaviour at the resolved scales; it does not satisfy \condC\ exactly. Failure of \condC\ does not preclude regularity when $\Teff>4$.

\begin{figure}[ht]
\centering
\includegraphics[width=0.55\textwidth]{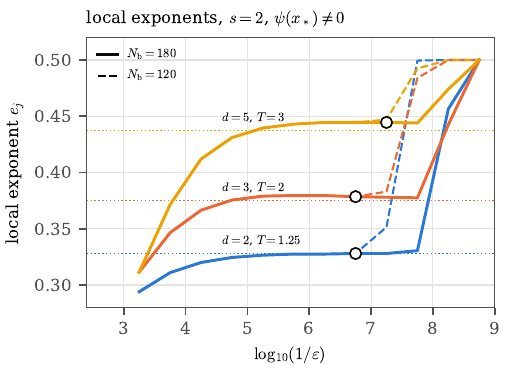}
\caption{Local exponents $e_j$ against $\log_{10}(1/\varepsilon)$ in discretisation (ii) for $\psi(x_*)\ne0$ and $s=2$, at basis degrees $N_{\rm b}=180$ (solid) and $N_{\rm b}=120$ (dashed). Open circles mark the selected exponent, dotted lines the prediction \eqref{num:eq:pred}. The rise to $1/2$ on the right is the truncation artefact.}
\label{num:fig:slopes}
\end{figure}

\paragraph*{Local exponents}
In discretisations (i) and (ii), the modulus is computed at $\varepsilon_j=10^{-3-j/2}$, $j=0,\dots,12$. The local exponent
\[
 e_j=\frac{\log\omega(\varepsilon_{j+1})-\log\omega(\varepsilon_j)}{2\{\log\varepsilon_{j+1}-\log\varepsilon_j\}}
\]
is the absolute-error exponent read off one step of the $\varepsilon$-grid. In discretisations (i) and (ii), the fixed finite-dimensional forward map $A$ is injective, so $\omega(\varepsilon)\asymp\varepsilon$ for $\varepsilon$ below the smallest resolved singular value, and $e_j$ eventually rises to $1/2$. This is a truncation artefact. In discretisation (iii), the nodal centred-flux scheme on the odd grids has a nullspace detected by the Gaussian target. Its fixed-grid modulus therefore has a positive limit as $\varepsilon\downarrow0$; the selected slopes describe the resolved finite-scale range rather than that fixed-grid limit. We therefore report one \emph{selected} local exponent per case, chosen by the rules below, which were fixed in advance and are applied by script. Figure~\ref{num:fig:slopes} shows local exponents at two resolutions together with the selected values; the selected values for (i) and (ii) are listed in Table~\ref{num:tab:exponents}, while the square-model values are reported under \emph{Further observations} below.

\paragraph*{Further observations}
Over the resolved ranges the computations are consistent with Theorem~\ref{thm:linear}.
\begin{enumerate}
\item \emph{Effective dimension.} For example, $d=2$, $\mu=0.125$ gives $0.320$ against the prediction $0.320$, while the value with $d$ in place of $\Teff$ is $0.375$; $d=3$, $\mu=4$ gives $0.390$ against $0.391$ (value with $d$: $0.438$).
\item \emph{Threshold.} In $d=5$, for $s=2$ we obtain $0.416$, $0.445$ and $0.468$ at $\Teff=2.6,3,3.4$, against $0.412$, $0.438$ and $0.463$. At $\Teff=4.5$ the selected value is $0.487$, still increasing at the resolved scales.
\item \emph{Critical regime.} For the four-dimensional ball, $(\omega/\varepsilon)^2$ increases by an almost constant amount per unit of $\log(1/\varepsilon)$ over five decades. This is consistent with the factor $\sqrt{\log N}$; the selected exponents near $0.48$ for $\Teff=4$ in Table~\ref{num:tab:exponents} reflect this factor.
\item \emph{The two irregular targets.} The selected exponents for $\psi(x_*)\ne0$ and for a nonconstant profile agree within $0.014$, the largest differences occurring for the isotropic ball in $d=3$ and for moderate anisotropy.
\item \emph{Square model.} The selected local exponents decrease with anisotropy and are close to the predictions over the resolved range: for $s=2$ and $\kappa=1,1.5,2,3,6,10$ they are $0.374$, $0.362$, $0.351$, $0.339$, $0.326$, $0.321$, against $0.375$, $0.354$, $0.344$, $0.333$, $0.323$, $0.319$. This experiment imposes $h=0$ on the square boundary. It does not require $h$ to vanish in a boundary collar, and its zero extension need not belong to the ordinary $H^s$ class on a larger domain. It therefore illustrates the anisotropic trend in a separate local model.
\end{enumerate}

\paragraph*{Limitations}
In discretisations (i) and (ii), the parameter sets are spectral or polynomial Hilbert balls, not the balls $\{\norm h_{H^s}\le\eta\}$ of Theorem~\ref{thm:linear}. For interior-localised functions they are equivalent to such balls, but near $\partial\cO$ they impose different boundary behaviour, and the flow-profile obstruction involves perturbations that reach the boundary (Remark~\ref{rem:collar}). The selected exponents are one-step slopes at finite resolution, not limits. For moderate anisotropy ($\mu\in[0.5,0.7]$) and $\psi(x_*)\ne0$ they exceed the prediction by up to $0.01$, for example $0.409$ against $0.400$ for $d=3$, $\mu=0.7$, $s=2$, and this does not change with the basis degree. We attribute it to slowly decaying corrections from the slow directions. Their relative weight at fast scale $r$ is $r^{2(1-\mu)}$, which pulls the slope towards the isotropic value at accessible scales. A surrogate model with these terms reproduces the sign and size of the deviation.

\paragraph*{The discrete modulus and its computation}
Each discretisation fixes a finite-dimensional space $V$ of perturbations $h$, a matrix representing $A=A_{f_0}$ on $V$ with values measured in $L^2$, and a positive definite Gram matrix $Q_2$ of the Hilbert norm $\norm\cdot_\Theta$. Put $Q_1=A^*A/\varepsilon^2$. The discrete modulus $\omega_V(\varepsilon)$ is the support function in direction $\psi$ of $\{h\in V:\ip h{Q_1h}\le1,\ \ip h{Q_2h}\le1\}$.

\emph{Exact identity.} This is a convex problem with the Slater point $h=0$. Lagrangian duality therefore gives
\begin{equation}\label{num:eq:dual}
 \omega_V(\varepsilon)=\inf_{t\in(0,1)}\ip\psi{(tQ_1+(1-t)Q_2)^{-1}\psi}^{1/2}.
\end{equation}
With $C=AQ_2^{-1/2}=U\Sigma V^*$ and $q=V^*Q_2^{-1/2}\psi$,
\[
 \ip\psi{(tQ_1+(1-t)Q_2)^{-1}\psi}=\sum_i\frac{q_i^2}{t\sigma_i^2/\varepsilon^2+1-t}.
\]
Independently of \eqref{num:eq:dual}, the inclusions $\{\ip h{(Q_1+Q_2)h}\le1\}\subset\{Q_1,Q_2\text{-balls}\}\subset\{\ip h{(tQ_1+(1-t)Q_2)h}\le1\}$ show that the value of the right-hand side at $t=1/2$ lies within a factor $\sqrt2$ of $\omega_V(\varepsilon)$ \cite{donoho1994statistical}.

\emph{Numerical errors.} In discretisations (i) and (ii) the minimum in \eqref{num:eq:dual} is taken over a finite logarithmic grid of $800$--$1000$ values of $t$, clustered at both ends of $(0,1)$. The computed value is therefore an upper bound for $\omega_V(\varepsilon)$. Refining the grid eightfold changed $\omega_V$ by at most $1.3\times10^{-4}$ relative and the local exponents by at most $6\times10^{-5}$ in the cases tested ($d=3$ and $d=2$ in (i); $d=2,\mu=0.25$ and $d=3,\mu=0.5$ in (ii) at $N_{\rm b}=80$). Further errors come from quadrature and from the singular value decomposition, both at the level of machine precision. The dominant error is the truncation to $V$, which is discussed under \emph{Selection rules} below.

In discretisation (iii), $\psi$ is not available in a diagonalising basis. For $s=2$, the $N_{\rm g}=127$, $\kappa=1$ run uses $m=10^{2+2i}$, $i=0,\dots,6$, while the runs at $N_{\rm g}=255$ and $511$ use $m=10^{4+i}$, $i=0,\dots,14$. For $s=3$, the runs at $N_{\rm g}=255$ use $m=10^{6+i}$, $i=0,\dots,14$. For each multiplier $m$, we compute the maximiser $h_m$ of $\ip h\psi-\frac12(\norm h_\Theta^2+m\norm{Ah}_2^2)$ from its sparse optimality system. The normalised point $h_m/\norm{h_m}_\Theta$ satisfies the Karush--Kuhn--Tucker conditions of \eqref{num:eq:modulus} on $V$ at $\varepsilon=\norm{Ah_m}_2/\norm{h_m}_\Theta$. Thus $(\varepsilon,\ip{h_m}\psi/\norm{h_m}_\Theta)$ is an exact point of the graph of $\omega_V$, up to the error of the sparse linear solver. For $s=3$, the sparse-solver error becomes visible at the two largest multipliers; the corresponding last two adjacent-interval slopes are excluded from exponent selection.

\paragraph*{Discretisation (i)}
On $B_1$ with $f_0=g=1$, $\nabla u_0=x/d$. On the sector of angular degree $\ell$, $\diver(h\nabla u_0)=(rh'+dh)/d$ for $h=h(r)Y_\ell$, where $Y_\ell$ is an $L^2(S^{d-1})$-normalised spherical harmonic. The basis is $\phi_k(r)=c_kr^{1-d/2}J_\nu(j_kr)$, with $\nu=\ell+d/2-1$, $J_\nu(j_k)=0$, and $c_k$ chosen so that $\int_0^1\phi_k^2r^{d-1}dr=1$. The matrix $\ip{\phi_m}{(r\phi_k'+d\phi_k)/d}$ is computed with $2\times10^4$ Gauss--Legendre nodes, and $\Delta_D^{-1}$ acts diagonally through $-j_m^{-2}$. The class is $\sum_kj_k^{2s}a_k^2\le1$ for $h=\sum_{k\le K}a_k\phi_k$, with $K=150$; $K=120$ and $K=250$ give the same local exponents to three digits. Each target lies in a single sector, which is invariant under the forward map. Targets enter through their coefficients $\ip{\phi_k}\psi$, that is, through their $L^2$ projection onto the span of the basis. The targets are:
\begin{itemize}
\item $\psi(x)=e^{-|x|^2/0.16}$ ($\ell=0$, $\psi(0)\ne0$);
\item $\psi(x)=r\,b(r)Y_1(x/r)$ with $b(r)=e^{-((r-0.45)/0.15)^2}$ ($\ell=1$, $\psi(0)=0$). For the flow $\dot x=x/d$ the profile is proportional to $Y_1\int_0^1\psi_1(r)r^{-1}dr$, which is nonconstant;
\item the control $\psi=r\{b(r)-a\,b_2(r)\}Y_1$, with $b_2(r)=e^{-((r-0.25)/0.08)^2}$ and $a$ chosen so that $\int_0^1\psi_1(r)r^{-1}dr=0$. Then \condC\ holds exactly.
\end{itemize}
The degree-one radial-Gaussian targets in (i) are not $C^\infty$ at the origin. The control nevertheless admits a boundary-zero $H^2$ transport solution and hence an $L^2$ representer.

\paragraph*{Discretisation (ii)}
In the coordinates $\xi_i=x_i/a_i$ on the unit ball, $L=\sum_ia_i^{-2}\partial_{\xi_i}^2$ and $\diver(h\nabla u_0)=2c\sum_ia_i^{-2}\partial_{\xi_i}(\xi_ih)$. Write $\xi=(\xi',\xi_d)$, $r=|\xi|$ and $z=\xi_d$. The basis consists of
\[
 \Phi_{jk}=r^kP_j^{(0,k+d/2-1)}(2r^2-1)\,C_k^{((d-2)/2)}(z/r),\qquad k\text{ even},\quad n=2j+k\le N_{\rm b},
\]
normalised numerically in $L^2(B_1,d\xi)$; orthonormality holds to $3\times10^{-13}$. For $d=2$ we use the limiting convention $C_k^{(0)}=T_k$, the Chebyshev polynomial of the first kind. Thus $C_k^{(0)}(\cos\theta)=\cos k\theta$, which agrees with $\lim_{\lambda\to0}C_k^{(\lambda)}/\lambda$ up to a factor removed by the normalisation. These functions are radial in $\xi'$ and even in $z$, a subspace invariant under the forward map. Solutions are sought as $(1-r^2)q$ with $\deg q\le N_{\rm b}$. Since the Dirichlet problem on an ellipsoid maps polynomial data of degree $n$ to such solutions, the Galerkin forward map is exact on $V$. $\norm{Ah}_2$ is computed through the Cholesky factor of the Gram matrix of $\{(1-r^2)\Phi_{jk}\}$. All integrals use Gauss--Jacobi quadrature in $r^2$ and in $z/r$ (weight $(1-t^2)^{(d-3)/2}$), and are exact for the polynomial integrands. The class weight is $(1+n(n+d))^{s/2}$. We use $N_{\rm b}\in\{120,150,180\}$, giving up to $4186$ basis functions. Targets enter through their coefficients $\ip{\Phi_{jk}}\psi$. In $\xi$-coordinates they are:
\begin{itemize}
\item $\psi=e^{-|\xi|^2/0.16}$ ($\psi(0)\ne0$);
\item $\psi=10\,(z^2-|\xi'|^2/(d-1))\,e^{-|\xi|^2/0.2025}$, a degree-two harmonic times a Gaussian ($\psi(0)=0$, nonconstant profile);
\item the control $\psi=\nabla u_0\cdot\nabla w$ with $w=e^{-|\xi|^2/0.16}(1+3z^2)$. The boundary-zero solution of \eqref{eq:transport} is $w$ minus the function that is constant along flow lines and equals $w$ on $\partial\cO$. On $\partial\cO$, $w$ ranges over $e^{-6.25}[1,4]$, so the flow profile is constant up to a variation $3e^{-6.25}\approx5.8\times10^{-3}$, compared with $w(0)=1$. Hence \condC\ holds only approximately; at the resolved scales this target behaves like a regular one.
\end{itemize}
None of these targets is compactly supported.

\paragraph*{Discretisation (iii)}
Centred differences on a uniform grid with spacing $2/(N_{\rm g}+1)$ are used for $\Delta_D$ and for $\diver(h\,Hx)$, together with the norms $\norm{\Delta h}_2$ and $(\norm{\nabla\Delta h}_2^2+\norm{\Delta h}_2^2)^{1/2}$. The target is $\psi(x)=e^{-|x|^2/0.18}$. Each optimality system is a sparse saddle-point system of size $3N_{\rm g}^2$. The $\varepsilon$-values are those produced by the multipliers; here the local exponent is defined as above between consecutive frontier points.

\paragraph*{Selection rules}
The rules were fixed before the values were tabulated, and are applied by the script that produces Table~\ref{num:tab:exponents}. Write $\gamma_j=2e_j$ for the local slope of $\log\omega$.
\begin{itemize}
\item In (i), the last local exponent ($\varepsilon$ between $10^{-8.5}$ and $10^{-9}$) is used. At $K=150$ there is no truncation artefact in this range.
\item In (ii), if two basis degrees $N_1<N_2$ are available, let $j^*$ be the largest index such that $|\gamma_i^{(N_1)}-\gamma_i^{(N_2)}|\le0.005$ for all $3\le i\le j^*$; the selected value is $e^{(N_2)}_{j^*}$. The tolerance $0.005$ thus applies to the local slope of $\log\omega$, which is twice the exponent $e_j$. Otherwise, the onset of the truncation artefact is the first index $i>5$ (so $\varepsilon_i\le10^{-5.5}$) with $\gamma_i-\gamma_{i-1}>0.04$, again in units of $\gamma$. The selected value is $e_{i-2}$, or the last $e_j$ if there is no such index.
\item In (iii), the largest local exponent is used, since mesh effects lower the slopes at small $\varepsilon$.
\end{itemize}

\paragraph*{Resolution checks}
For $\mu=1$ and the common target with $\psi(0)\ne0$, (ii) reproduces the exponents of (i) to within $0.001$ for $d=2$: $0.374$ and $0.3995$ against $0.375$ and $0.3995$ for $s=2,3$. In (ii), increasing $N_{\rm b}$ from $120$ to $180$ leaves the local exponents unchanged to three digits up to the onset of the truncation artefact. It only moves that onset to smaller $\varepsilon$, by about half a decade (Figure~\ref{num:fig:slopes}). In (iii), $\kappa=1$ and $s=2$ give $e=0.369$, $0.372$ and $0.374$ for $N_{\rm g}=127$, $255$ and $511$, respectively, compared with the prediction $0.375$. These values use different multiplier and $\varepsilon$ windows, so this is a descriptive cross-resolution comparison rather than a convergence test on a common window. For $s=5$ the accessible range of $\varepsilon$ is pre-asymptotic in all discretisations, and such runs are not reported.

\begin{table}[ht]
\caption{Selected local exponents $e$ in the linearised experiment, for $\psi(x_*)\ne0$ (``loc.'') and for $\psi(x_*)=0$ with nonconstant flow profile (``prof.''). ``pred.'' is $e_{\rm pred}(s,\Teff)$ from \eqref{num:eq:pred}; a dagger marks $\Teff=4$, where the rate carries an additional logarithmic factor. ``$d$'' is $\min\{(2s+d)/(4s+8),1/2\}$, the value with $d$ in place of $\Teff$. Rows ``ball'' are discretisation (i); rows with $\mu$ are discretisation (ii), with Hessian eigenvalues proportional to $(\mu,\dots,\mu,1)$. ``res.'' is the resolution of the selected run. All control targets give $e\ge0.4965$.}
\label{num:tab:exponents}
\centering\footnotesize
\setlength{\tabcolsep}{2.6pt}
\begin{tabular}{@{}rlrl cccc cccc@{}}
\toprule
& & & & \multicolumn{4}{c}{$s=2$} & \multicolumn{4}{c}{$s=3$}\\
\cmidrule(lr){5-8}\cmidrule(l){9-12}
$d$ & geometry & $\Teff$ & res. & loc. & prof. & pred. & $d$ & loc. & prof. & pred. & $d$\\
\midrule
2 & ball & 2 & $K=150$ & 0.375 & 0.376 & 0.375 & 0.375 & 0.400 & 0.403 & 0.400 & 0.400 \\
3 & ball & 3 & $K=150$ & 0.435 & 0.442 & 0.438 & 0.438 & 0.445 & 0.459 & 0.450 & 0.450 \\
4 & ball & 4 & $K=150$ & 0.482 & 0.491 & 0.500$^\dagger$ & 0.500 & 0.479 & 0.492 & 0.500$^\dagger$ & 0.500 \\
5 & ball & 5 & $K=150$ & 0.498 & 0.500 & 0.500 & 0.500 & 0.495 & 0.499 & 0.500 & 0.500 \\
6 & ball & 6 & $K=150$ & 0.500 & 0.500 & 0.500 & 0.500 & 0.499 & 0.500 & 0.500 & 0.500 \\
\addlinespace
2 & $\mu=0.5$ & 1.5 & $N_{\rm b}=180$ & 0.345 & 0.346 & 0.344 & 0.375 & 0.376 & 0.377 & 0.375 & 0.400 \\
2 & $\mu=0.25$ & 1.25 & $N_{\rm b}=180$ & 0.328 & 0.328 & 0.328 & 0.375 & 0.362 & 0.362 & 0.362 & 0.400 \\
2 & $\mu=0.125$ & 1.125 & $N_{\rm b}=180$ & 0.320 & 0.320 & 0.320 & 0.375 & 0.356 & 0.356 & 0.356 & 0.400 \\
\addlinespace
3 & $\mu=4$ & 2.25 & $N_{\rm b}=180$ & 0.390 & 0.390 & 0.391 & 0.438 & 0.411 & 0.411 & 0.412 & 0.450 \\
3 & $\mu=2$ & 2.5 & $N_{\rm b}=120$ & 0.405 & 0.407 & 0.406 & 0.438 & 0.424 & 0.426 & 0.425 & 0.450 \\
3 & $\mu=0.7$ & 2.4 & $N_{\rm b}=150$ & 0.409 & 0.405 & 0.400 & 0.438 & 0.427 & 0.423 & 0.420 & 0.450 \\
3 & $\mu=0.5$ & 2 & $N_{\rm b}=180$ & 0.379 & 0.377 & 0.375 & 0.438 & 0.403 & 0.402 & 0.400 & 0.450 \\
3 & $\mu=0.4$ & 1.8 & $N_{\rm b}=150$ & 0.364 & 0.363 & 0.362 & 0.438 & 0.391 & 0.391 & 0.390 & 0.450 \\
3 & $\mu=0.25$ & 1.5 & $N_{\rm b}=120$ & 0.344 & 0.344 & 0.344 & 0.438 & 0.375 & 0.375 & 0.375 & 0.450 \\
\addlinespace
5 & $\mu=2$ & 4.5 & $N_{\rm b}=180$ & 0.487 & 0.488 & 0.500 & 0.500 & 0.488 & 0.490 & 0.500 & 0.500 \\
5 & $\mu=0.75$ & 4 & $N_{\rm b}=180$ & 0.486 & 0.478 & 0.500$^\dagger$ & 0.500 & 0.490 & 0.483 & 0.500$^\dagger$ & 0.500 \\
5 & $\mu=0.6$ & 3.4 & $N_{\rm b}=150$ & 0.468 & 0.458 & 0.463 & 0.500 & 0.474 & 0.466 & 0.470 & 0.500 \\
5 & $\mu=0.5$ & 3 & $N_{\rm b}=180$ & 0.445 & 0.436 & 0.438 & 0.500 & 0.456 & 0.449 & 0.450 & 0.500 \\
5 & $\mu=0.4$ & 2.6 & $N_{\rm b}=150$ & 0.416 & 0.412 & 0.412 & 0.500 & 0.433 & 0.429 & 0.430 & 0.500 \\
5 & $\mu=0.25$ & 2 & $N_{\rm b}=120$ & 0.376 & 0.374 & 0.375 & 0.500 & 0.401 & 0.400 & 0.400 & 0.500 \\
\bottomrule
\end{tabular}
\end{table}

\clearpage
\bibliographystyle{imsart-number}
\bibliography{refs}

\end{document}